\documentclass[11pt,leqno]{article}
\usepackage[letterpaper,margin=0.9in]{geometry}
\usepackage{amsmath,amssymb,amsthm,mathtools}
\usepackage{enumitem}
\usepackage[hidelinks]{hyperref}
\usepackage[capitalise,noabbrev]{cleveref}
\usepackage{tocloft}
\numberwithin{equation}{section}
\theoremstyle{plain}
\newtheorem{thm}{Theorem}[section]

\newtheorem{lem}[thm]{Lemma}
\newtheorem{cor}[thm]{Corollary}

\newtheorem{prop}[thm]{Proposition}

\theoremstyle{definition}

\theoremstyle{remark}
\newtheorem{rem}[thm]{Remark}

\allowdisplaybreaks

\newcommand{\nc}{\newcommand}
\nc{\hb}{\mathbb} \nc{\M}{\mathcal} \nc{\mf}{\mathfrak} \nc{\mbf}{\mathbf}
\nc{\DMO}{\DeclareMathOperator}

\newbox\noforkbox \newdimen\forklinewidth \forklinewidth=0.3pt \setbox0\hbox{$\textstyle\smile$}
\setbox1\hbox to \wd0{\hfil\vrule width \forklinewidth depth-2pt height 10pt \hfil} \wd1=0 cm
\setbox\noforkbox\hbox{\lower 2pt\box1\lower 2pt\box0\relax}
\def\anchor{\mathop{\copy\noforkbox}\limits}

\setbox0\hbox{$\textstyle\smile$} \setbox1\hbox to \wd0{\hfil{\sl /\/}\hfil} \setbox2\hbox to
\wd0{\hfil\vrule height 10pt depth -2pt width \forklinewidth\hfil} \wd1=0 cm \wd2=0 cm
\newbox\doesforkbox \setbox\doesforkbox\hbox{\box1 \lower 2pt\box2\lower2pt\box0\relax}
\def\nanchor{\mathop{\copy\doesforkbox}\limits}

\nc{\cA}{{\M A}} \nc{\cB}{{\M B}} \nc{\cC}{{\M C}} \nc{\cD}{{\M D}} \nc{\cE}{{\M E}} \nc{\cF}{{\M
F}} \nc{\cG}{{\M G}} \nc{\cH}{{\M H}} \nc{\cI}{{\M I}} \nc{\cJ}{{\M J}} \nc{\cK}{{\M K}}
\nc{\cL}{{\M L}} \nc{\cM}{{\M M}} \nc{\cN}{{\M N}} \nc{\cO}{{\M O}} \nc{\cP}{{\M P}} \nc{\cQ}{{\M
Q}} \nc{\cR}{{\M R}} \nc{\cS}{{\M S}} \nc{\cT}{{\M T}} \nc{\cU}{{\M U}} \nc{\cV}{{\M V}}
\nc{\cW}{{\M W}} \nc{\cX}{{\M X}} \nc{\cY}{{\M Y}} \nc{\cZ}{{\M Z}}
\nc{\Aa}{{\hb A}} \nc{\Cc}{{\hb C}} \nc{\Gg}{{\hb G}} \nc{\Nn}{{\hb N}} \nc{\Pp}{{\hb P}}
\nc{\Qq}{{\hb Q}} \nc{\Rr}{{\hb R}} \nc{\Zz}{{\hb Z}}
\nc{\mfa}{{\mf a}} \nc{\mfb}{{\mf b}} \nc{\mfk}{{\mf k}} \nc{\mfm}{{\mf m}} \nc{\mfp}{{\mf p}}
\nc{\mfq}{{\mf q}} \nc{\mfr}{{\mf r}} \nc{\fP}{{\mf P}}
\DMO*{\trdeg}{td} \DMO*{\spec}{Spec} \DMO*{\fork}{\nanchor} \DMO*{\dnf}{\anchor} 
\DMO{\deter}{det} \DMO{\RM}{RM} \DMO{\RC}{RC} \DMO{\Real}{Re} \DMO{\Imag}{Im} \DMO{\tr}{tr}
\DMO{\qc}{QC} \DMO{\Hu}{Hull} \DMO{\leg}{length} \DMO{\area}{area} \DMO{\dia}{diameter}
\DMO{\iso}{Iso} \DMO{\dis}{dist} \DMO{\grad}{grad} \DMO{\vol}{volume} \DMO{\gra}{grad}
\DMO{\hd}{nbhd} \DMO{\dv}{div}
\DMO{\Psl}{PSL}
\nc{\Mb}{\mathfrak^{2b/\delta}_{K_x}} \nc{\Ma}{\mathfrak^{2a/\delta}_{K_x}}
\nc{\dif}{\mathrm{d}} \nc{\G}{\Gamma} \nc{\g}{\gamma} \nc{\D}{\nabla} \nc{\p}{\partial}
\nc{\DD}{\Delta^2} \nc{\pp}{\partial^2} \nc{\de}{\delta} \nc{\td}[2]{\trdeg{({#1}/{#2})}}
\nc{\dtd}[2]{\trdeg_{\delta}{({#1}/{#2})}} \nc{\dspec}[1]{\spec_{\delta}{#1}}
\nc{\ddim}[1]{\dimen_{\delta}{#1}} \nc{\gens}[1]{\langle {#1} \rangle} \nc{\gen}[2]{ {#1} \langle
{#2} \rangle } \nc{\form}{\Omega} \nc{\set}[1]{\left\{ {#1} \right\}} \nc{\mr}{\hat}
\nc{\pr}{\partial}
\nc{\bc}[3]{\cB^{#1}({#2},{#3})=B^{#1}_{#2}(C^{#2}_{#3}(t)+B^{#1}_{#3}(C^{#2}_{#3}(t))}
\nc{\tuple}[2]{{#1},\ldots,{#2}} \nc{\ptu}[2]{{#1}:\ldots:{#2}}

\nc{\maps}[3]{{#1}\!:\!{#2}\rightarrow{#3}} \nc{\map}[2]{{#1}\rightarrow {#2}} \nc{\res}[2]{{#1}
|_{#2}} \nc{\imbed}{\hookrightarrow}
\DMO{\Isom}{Isom} \DMO{\diam}{diam} \DMO{\inte}{int}

\title{Higher dimensional Kleinian groups of small Hausdorff dimension are classical Schottky
groups II.} \author{ Yong Hou } \date{}
\begin{document}
\maketitle
\let\thefootnote\relax \footnotetext{\emph{keywords:} Kleinian groups, higher dimensional hyperbolic spaces, Hausdorff
dimensions.}
\addtocounter{footnote}{-1}
\begin{abstract}
This is the higher dimensional continuation of our works on Hausdorff dimension of Kleinian groups
\cite{Hou2010}, \cite{Hou2020}. In this paper we prove that for all $r\ge 2$ there exists
$\lambda_r>0$, which is independent of the dimension $N$ of the hyperbolic space, such that any
finitely generated non-elementary torsion-free Kleinian group of $\Isom(\mathbb{H}^N)$ generated by
$r$ elements with limit set of Hausdorff dimension $<\lambda_r$ is a classical Schottky group. For
two generators the bound can be taken to be $1/3$. The proof relies on our previous works in
\cite{Hou2010}, \cite{Hou2020} and on a metric carrier graph reduction. Our result can also be
considered as the complete converse to the well-known theorem of Phillips-Sarnak \cite{PS} for
spaces of any dimension $N$.
\end{abstract}
\setcounter{tocdepth}{1}
\tableofcontents
\section{Introduction and Main Theorem}
Let $N\ge 4$ and let $\mathbb{H}^N$ be the real hyperbolic $N$-space with full isometry group
$\Isom(\mathbb{H}^N)$. For a given Kleinian group (discrete subgroup $\G$ of $\Isom(\mathbb{H}^N)$)
and a point $o\in\mathbb{H}^N$, the limit set $\Lambda_\G$ of $\G$ is the set of ideal accumulation
points of the orbit $\G o$ in the ideal boundary $\partial_\infty\mathbb{H}^N=\mathbb{S}^{N-1}$. The
group $\G$ is called \emph{non-elementary} if $\Lambda_\G$ contains more than two points. The set
$\Omega_\G=\mathbb{S}^{N-1}\setminus\Lambda_\G$ is called the region of discontinuity of $\G$, and
$\G$ acts properly discontinuously on $\Omega_\G$. Denote the Hausdorff dimension of $\Lambda_\G$ by
$\mathfrak{D}_\G$ and the critical exponent of $\G$ by $\delta_\G$ (see Section $2$ for the
notations). Most of the geometric definitions are the same as for the three dimensional hyperbolic
space $\mathbb{H}^3$, see \cite{Kapovich} for details.\par A closed round ball in
$\mathbb{S}^{d}\subset\mathbb{R}^{d+1}$, $d=N-1$, is a set
\[
B=\{\xi\in\mathbb{S}^d:u\cdot\xi\ge c\}, \qquad u\in\mathbb{S}^d,\quad -1<c<1.
\]
Such a ball is the ideal boundary of a closed hyperbolic halfspace in $\mathbb{H}^N$, and the
bounding round spheres are the ideal boundaries of totally geodesic hyperplanes. We call a group
classical Schottky group if, for some integer $r\ge 1$, it is generated by loxodromic
transformations $g_1,...,g_r$ and there exist $2r$ closed round balls
$B_i^+,B_i^-\subset\mathbb{S}^d$ with pairwise disjoint closures such that, for $1\le i\le r$,
\begin{equation}\label{pairing}
g_i\bigl(\mathbb{S}^d\setminus\inte B_i^-\bigr)=B_i^+.
\end{equation}
The generators $g_1,...,g_r$ are called classical generators of the group.\par In \cite{Hou2020} Hou
proves a positive Hausdorff dimension bound, independent of the rank, for the finitely generated
torsion-free non-elementary subgroups of $\Psl(2,\mathbb{C})$ to be classical Schottky
(\cite{Hou2020}, Theorem 1.1 and Corollary 15.2). In \cite{Hou} Hou proves a classification theorem
for every non-elementary, torsion-free, convex cocompact Kleinian groups of $\Psl(2,\mathbb{C})$ of
Hausdorff dimension at most one. In this paper we prove that the Kleinian groups of small Hausdorff
dimension of $\Isom(\mathbb{H}^N)$ are classical Schottky groups for every dimensions $N$, with a
bound that depends only on the number of generators.\par For $x>0$, denote $\log^+(x)=\max\{0,\log
x\}$. For an integer $r\ge 2$, we set
\[
m_r=3r-3,\qquad B_r=2m_r+5,\qquad C_{0,r}=m_r B_r^{m_r},\qquad A_r=2\cdot203^{m_r},
\]
and define
\[
\lambda_r=\frac1{1000A_rC_{0,r}} =\frac1{2000(3r-3)\,[203(6r-1)]^{3r-3}}.
\]
Our main results are:
\begin{thm}\label{main}
Let $\G<\Isom(\mathbb{H}^N)$ be discrete, torsion-free, non-elementary, and generated by $r\ge 2$
elements. If
\[
\mathfrak{D}_\G\le\lambda_r,
\]
then $\G$ is classical Schottky. The same conclusion holds if $\delta_\G\le\lambda_r$.
\end{thm}
\begin{thm}\label{main-2}
Let $\G<\Isom(\mathbb{H}^N)$ be discrete, torsion-free, non-elementary, and generated by two
elements. If
\[
\mathfrak{D}_\G\le\frac13,
\]
then $\G$ is classical Schottky. The same conclusion holds if $\delta_\G\le 1/3$.
\end{thm}
Note that the bound $\lambda_r$ in Theorem \ref{main} do not dependent on $N$,
and that for two generators Theorem \ref{main-2} gives the much stronger universal bound $1/3$.
Theorem \ref{main} can be viewed as the complete converse to the result by Phillips-Sarnak \cite{PS}
for spaces of any dimension $N$.\par For $\G<\Isom^+(\mathbb{H}^2)$, Button proves that every free
convex cocompact group is classical \cite{Button}. For $\G<\Isom^+(\mathbb{H}^3)$, the corresponding
small dimension results follow from Hou's theorems \cite{Hou2020}, \cite{Hou}. In Theorems \ref{n=3}
and \ref{n=2} we give the proofs for the full isometry groups in dimensions three and two, by
extending the actions to higher dimensions and restricting round domains to invariant subspaces.
Hence it suffices to prove the main theorems for $N\ge 4$, which we assume from now on unless
otherwise stated.\par The main theorems have the following subgroup and boundary consequences:
\begin{cor}\label{main-cor}
Let $\G<\Isom(\mathbb{H}^N)$ be discrete and torsion-free. If
\[
\mathfrak{D}_\G\le\frac13,
\]
then every non-elementary subgroup of $\G$ generated by two elements is classical Schottky. The same
conclusion holds if $\delta_\G\le 1/3$. Here $\G$ has arbitrary rank.
\end{cor}
\begin{proof}
This follows from Theorem \ref{main-2}.
\end{proof}
\begin{cor}\label{tree-cor}
Let $\G<\Isom(\mathbb{H}^N)$ be discrete, torsion-free, and non-elementary. Assume that either $\G$
is generated by $r\ge 2$ elements and $\mathfrak{D}_\G\le\lambda_r$, or $\G$ is generated by two
elements and $\mathfrak{D}_\G\le 1/3$. Then $\Omega_\G=\mathbb{S}^d\setminus\Lambda_\G$ is path
connected and simply connected. There is a $\G$-invariant embedded tree
\[
T\subset\Omega_\G
\]
such that
\[
\overline T=T\cup\Lambda_\G,\qquad \mathcal{H}^1(\overline T)<\infty.
\]
Its closure is a compact dendrite with endpoint set $\Lambda_\G$. The same conclusions hold if the
relevant Hausdorff-dimension bound is replaced by the corresponding critical-exponent bound.
\end{cor}
\begin{proof}
Under either Hausdorff-dimension assumption $\mathfrak{D}_\G\le 1/3<1/2$, so that Theorem
\ref{structural} gives freeness and convex cocompactness, and under either critical-exponent
assumption the same theorem implies $\mathfrak{D}_\G=\delta_\G\le 1/3$ and convex cocompactness, so
in both cases $\mathcal{H}^1(\Lambda_\G)=0$. Since $N\ge 4$, Corollary \ref{simply-connected} gives
path connectedness and simple connectedness of $\Omega_\G$, and Theorem \ref{tree} gives the stated
tree and its closure.
\end{proof}
Note that for a finitely generated torsion-free non-elementary Kleinian group, the structural
theorem of Hou \cite{Hou2026}, which we recall in Section $3$, gives geometric finiteness and
freeness when $\mathfrak{D}_\G<1$. It also gives convex cocompactness when there are no parabolic
elements, which is the case when $\mathfrak{D}_\G<1/2$.\par
\subsection{Strategy of the proof}
The proof is divided into six main steps:\par Note that for a given free Kleinian group of rank $k$,
when we speak of a set of generators or a free basis, we always mean a set of $k$ free
generators.\par {\bf Step one:} We reduce to free convex cocompact groups. In the critical-exponent
formulation, Hou's dimension theorem (\cite{Hou2026}, Corollary 1.6) implies
$\dim_H\Lambda_\G=\delta_\G$, and in the Hausdorff-dimension formulation the dimension bound is
given and $\delta_\G\le\dim_H\Lambda_\G$. Since $\dim_H\Lambda_\G\le\lambda_r<1/2$, the structural
theorem of Hou makes $\G$ free and convex cocompact, of some free rank $k$ with $2\le k\le r$, and
as $\lambda_r\le\lambda_k$ we may work at rank $k$, with $m=m_k$, $C_0=C_{0,k}$ and $A=A_k$. This is
done in Section $3$, where we also give the round-ball covers of the limit set and the basic pairing
criteria for classical generators.\par {\bf Step two:} We state the displacement estimates that are
used throughout the paper. As in Hou's work, the proof combines a small critical exponent with
displacement estimates, changes of generators and the separation of paired round regions
(\cite{Hou2020}, Section 1.1 and Sections 3--6). The principal estimate is Hou's displacement
inequality (\cite{Hou2001}, Theorem 1.1), also stated in \cite{Hou2020} (Theorem 3.1 and Corollary
3.2): for freely generating $a,b\in\G$ and $x\in\mathbb{H}^N$,
\[
(e^{\delta_\G\dis(x,ax)}-1) (e^{\delta_\G\dis(x,bx)}-1)\ge 4.
\]
These estimates and the cyclic geometry of axes and tubes are given in Section $4$.\par {\bf Step
three:} For two generators we combine the displacement estimates with an explicit hyperplane pairing
for a given generating pair, and then change the generators by Nielsen transformations, following
the generator selection of \cite{Hou2020}. The perpendicular cut estimates and a terminating length
reduction allow arbitrary orthogonal parts around the axes. This gives the bound $1/3$ of Theorem
\ref{main-2} and is done in Section $5$.\par {\bf Step four:} In higher rank the strategy is carried
out through a shortest carrier graph, its rank-one components and the perpendicular midpoint
hyperplanes of its edges. Let $G$ be a geodesic graph of minimum total length, with at most $m$
edges, whose mapping to $\Hu(\Lambda_\G)/\G$ induces an isomorphism on fundamental groups. A gap
between edge lengths separates short edges from long edges, and each component formed by short edges
has length at most $C_0$. Hou's two-generator inequality excludes two independent loops with
displacement at most $2C_0$, since $2\delta_\G C_0<\log 3$, so each component is a tree or has rank
one. For a rank-one component $K_i$ with circuit generator $a_i$, axis $A_i$ and $L_i=\ell(a_i)$,
set $U_i=C_0+\log^+(1/L_i)$. The component is enlarged by adjoining edges and tree components. The
estimates (\ref{growth}) and (\ref{cycle}) exclude a second circuit and exclude joining two rank-one
components, and the enlarged components satisfy
\[
Z_i=\leg(K_i)+U_i\le AU_i, \qquad L_e>200Z_i
\]
for every remaining edge adjoining $K_i$, where $L_e$ is its length. The carrier graphs and their
reductions are given in Section $6$, and the enlargement of the components in Section $7$.\par {\bf
Step five:} The remaining edges are treated through their perpendicular midpoint hyperplanes. If
$p,q$ are the lifted endpoint images of such an edge, with $p$ the image of an endpoint in the
chosen lift of $K_i$, then minimality implies $\dis(a_i^np,q)\ge\dis(p,q)$ for every integer $n$,
which bounds the displacement along the axis. For the midpoint hyperplane $\mathcal{P}$, the ideal
ball $\mathcal{B}$ on the $q$ side and $D=\dis(A_i,\mathcal{P})$, this leads to
\[
D>98U_i-2,\qquad \{\log|\zeta|:\zeta\in\mathcal{B}\} \subset I_{\mathcal{B}},\qquad
|I_{\mathcal{B}}|<\frac{L_i}{4m},
\]
in coordinates with axis endpoints $0,\infty$, where $I_{\mathcal{B}}$ is a closed real interval.
Since $a_i$ adds $L_i$ to the logarithmic radius, distinct $<a_i>$-translates of the same ball are
disjoint. If, for a fixed component, balls from different edge orbits intersect, then the midpoint
overlap estimate (\ref{overlap}) produces a shorter carrier whose mapping still induces an
isomorphism on fundamental groups, so these balls are disjoint, and the same argument treats the
remaining tree components. At each rank-one component the at most $2m$ intervals have total length
less than $L_i/2$, which leaves room to choose one ball from each orbit between two concentric
spheres paired by $a_i$, and the closed balls defining each region are then pairwise disjoint.\par
{\bf Step six:} The resulting hyperbolic regions are arranged along a tree with a free $\G$-action,
adjacent regions having matching boundary hyperplanes and opposite halfspaces. Disjointness at each
region gives strict nesting along paths, hence disjointness of the balls associated with the edges
leaving a connected subtree, and a lift of a spanning tree of the quotient graph $X$ leaves
\[
2|E(X)|-2(|V(X)|-1)=2k
\]
pairwise disjoint closed round balls $B_i^\pm$, $1\le i\le k$. The remaining edges determine a free
basis $g_1,...,g_k$ of $\G$ with
\[
g_i(\mathbb{S}^d\setminus\inte B_i^-)=B_i^+ \qquad(1\le i\le k),
\]
and these balls form a classical domain for the original group. This finishes the proof of Theorem
\ref{main} in Section $7$.\par In the last two sections we give a necessary and sufficient condition
for a round domain in terms of the outward Lorentz normals, and we construct invariant rectifiable
trees in the region of discontinuity, which implies Corollary \ref{tree-cor}.\par
\subsection{Organization of the paper}
Section $2$ is used to define and list the global notations that will be used throughout the paper.
In Section $3$, we give the round-ball covers of the sets of zero one-dimensional Hausdorff measure,
the topology of their complements, the structural reduction, and the elementary pairing criteria. In
Section $4$, we prove the pressure and cyclic geometry estimates. In Section $5$, we prove Theorem
\ref{main-2}. In Section $6$, we develop the annular and carrier reductions. In Section $7$, we
prove Theorem \ref{main}. In Section $8$, we give the Lorentz normal criterion and the radius
duality results. In Section $9$, we prove the invariant metric and rectifiable tree results.\par
\subsection{Notations}
Let $\G\subset\Isom(\mathbb{H}^N)$ be a Kleinian group, where $N\ge 4$ and $d=N-1$. We will set the
following notations and definitions throughout the rest of the paper.

\begin{itemize}
\item For $s>0$ and $o\in\mathbb{H}^N$, the Poincar\'{e} series of $\G$ is
$P_\G(s,o)=\sum_{\g\in\G}e^{-s\dis(o,\g o)}$, and $\delta_\G$ denotes its critical exponent.
\item $\Lambda_\G\subset\mathbb{S}^d=\partial_\infty\mathbb{H}^N$ is the limit set of $\G$, the set
of ideal accumulation points of $\G o$ in $\mathbb{H}^N\cup\partial_\infty\mathbb{H}^N$, and
$\mathfrak{D}_\G=\dim_H\Lambda_\G$ is its Hausdorff dimension, measured in the usual round metric,
or equivalently in the chordal metric. We use the same symbol for an isometry and its boundary
extension.
\item For $x,y\in\mathbb{H}^N$, $[x,y]$ denotes their geodesic segment. We denote by $\dis$ distance
in the stated space, $\diam$ for diameter, $\inte$ for interior, $\partial$ for boundary in the
stated space, and $\overline{A}$ for closure. For subsets, $\dis(A,B)=\inf\{\dis(a,b):a\in A,b\in
B\}$, with $\inf\emptyset=+\infty$.
\item Euclidean vectors have norm $|\cdot|$ and scalar product $u\cdot v$. On $\mathbb{S}^d$,
$|\xi-\eta|$ is the chordal distance, $g_{\rm round}$ is the round Riemannian metric, and $ds$ is
its arclength element.
\item $\mathcal{H}^s$ is the $s$-dimensional Hausdorff measure in the stated metric, normalized so
that $\mathcal{H}^1$ agrees with arclength on smooth simple arcs.
\item $\Hu(K)$ is the closed hyperbolic convex hull of $K$, and $\partial_\infty A$ is the ideal
boundary of a hyperbolic subset $A$.
\item In the upper half-space model we use $(x,t)\in\mathbb{R}^{N-1}\times(0,\infty)$, with metric
$t^{-2}(|dx|^2+dt^2)$ and ideal boundary $\mathbb{R}^{N-1}\cup\{\infty\}$. In the hyperboloid model
of $\mathbb{H}^N$, we set
\[
\begin{gathered}
\mathcal{V}=\mathbb{R}^{N+1},\qquad \langle x,y\rangle=-x_0y_0+\sum_{j=1}^N x_jy_j,\\
\mathbb{H}^N=\{x\in\mathcal{V}:\langle x,x\rangle=-1,\ x_0>0\}.
\end{gathered}
\]
Orthogonal complements in this model use $\langle\cdot,\cdot\rangle$. A nonzero vector is timelike,
null, or spacelike according as its Lorentz square is negative, zero, or positive respectively. A
timelike or null vector is future pointing when its zeroth coordinate is positive. Ideal points are
future null rays. We represent isometries by Lorentz matrices preserving the future hyperboloid. For
$\langle n,n\rangle=1$, we set
\[
E_n=\{X\in\mathbb{H}^N:\langle X,n\rangle\ge 0\}, \qquad B_n=\partial_\infty E_n,
\]
and call these the positive halfspace and ideal ball of $n$.
\item $\G$ is elementary if $\Lambda_\G$ has at most two points, and non-elementary otherwise.
$\G<\Isom(\mathbb{H}^N)$ is called convex cocompact if $\G$ is discrete and preserves a nonempty
closed convex subset $C\subset\mathbb{H}^N$ such that $C/\G$ is compact. For non-elementary discrete
$\G$, this is equivalent to compactness of $\Hu(\Lambda_\G)/\G$.
\item For an isometry $g$, its translation length is $\ell(g)=\inf_{x\in\mathbb{H}^N}\dis(x,gx)$.
$g$ is elliptic if it fixes an interior point, parabolic if it has no interior fixed point and has a
unique ideal fixed point, and loxodromic if it has two ideal fixed points and $\ell(g)>0$.
Loxodromic transformations are called independent if their ideal fixed-point pairs are pairwise
disjoint.
\item A limit point is conical if a ray to it has points at arbitrarily large distance from its
initial point that stay within a bounded distance of $\G o$. A parabolic fixed point $\xi$ is
bounded if $\G_\xi\backslash(\Lambda_\G\setminus\{\xi\})$ is compact, where $\G_\xi$ is its
stabilizer. We use geometric finiteness in the equivalent sense that every limit point is conical or
bounded parabolic.
\item Given loxodromic transformations $g_1,...,g_r$, we denote by $<g_1,...,g_r>$ the group
generated by them, and $<a>$-translates for the translates by the cyclic group generated by $a$.
\end{itemize}

\section{Preliminary Estimates and Structures}
In this section we give the round-ball covers, the structural reduction and the pairing criteria
that will be used throughout the paper. Small Hausdorff dimension forces the one-dimensional
Hausdorff measure of the limit set to vanish, and this vanishing gives a direct separation argument
in every dimension. Note that in the following lemma the empty set is covered by the empty cover,
whose radius sum and maximum diameter are zero.

\begin{lem}\label{covers}
Let $K\subset\mathbb{R}^d$ be a compact set such that $\mathcal{H}^1(K)=0$. Then for all $\eta>0$
there exist an integer $m\ge0$ and closed Euclidean balls $D_1,...,D_m$, with radii $r_1,...,r_m>0$,
such that
\[
D_i\cap D_j=\emptyset\quad(i\ne j),\qquad K\subset\bigcup_{j=1}^m\inte D_j, \qquad \sum_{j=1}^m
r_j<\eta.
\]
We may also require every ball to meet $K$.
\end{lem}
\begin{proof}
Since $\mathcal{H}^1(K)=0$, there is a cover of $K$ by sets whose sum of diameters is arbitrarily
small, and each nonempty set can be enclosed in an open ball of radius slightly larger than its
diameter, the added radii having a total so small that the radius sum is less than $\eta$. By
compactness a subcover of $m_0$ balls remains, and after discarding those missing $K$ the closures
are taken.

Intersecting balls are merged in pairs. If one contains the other, the smaller is discarded.
Otherwise, with centers $c_1,c_2$, radii $r_1,r_2$ and $a=|c_1-c_2|$,
\[
|r_1-r_2|<a\le r_1+r_2,\qquad R=\frac{a+r_1+r_2}{2}\le r_1+r_2,
\]
and there is a point $c$ on the Euclidean segment $[c_1,c_2]$ with $|c-c_j|=R-r_j$ for $j=1,2$,
these two nonnegative distances summing to $a$. The closed Euclidean ball $\overline B(c,R)$
contains both closed balls and preserves their open interiors, since for $x$ in either open ball
\[
|x-c|\le |x-c_j|+|c_j-c|<r_j+(R-r_j)=R.
\]
Each merge decreases the number of balls, preserves the open cover and does not increase the radius
sum, so after at most $m_0-1$ merges the closed balls are disjoint, and they are relabeled
$D_1,...,D_m$. Note that every ball retained or obtained by merging meets $K$.
\end{proof}

\begin{cor}\label{nested}
Let $K\subset\mathbb{S}^d$ be a compact set such that $\mathcal{H}^1(K)=0$. Then there exists a
nested sequence of covers by disjoint closed round balls whose interiors cover $K$ and whose maximum
diameters tend to zero. Each ball of a generation lies in the interior of a ball of the preceding
generation, and each ball meets $K$. In particular, this conclusion holds when $\dim_H K<1$.
\end{cor}
\begin{proof}
Stereographic projection from a point outside $K$ produces a compact set $K'\subset\mathbb{R}^d$
with $\mathcal{H}^1(K')=0$, since the projection and its inverse are Lipschitz on compact
neighborhoods. Lemma \ref{covers} provides the first cover, generation zero. For a parent ball $D$,
disjointness and the open-cover property imply
\[
\emptyset\ne K_D=K'\cap D\subset\inte D,\qquad a_D=\dis(K_D,\partial D)>0.
\]
For generation $j\ge1$, the lemma is applied to $K_D$ with
\[
\eta<\min\{a_D/3,2^{-j-1}\}.
\]
Each child of radius $\rho<\eta$ meets $K_D$, so every point of that child lies within $2\rho<a_D$
of $K_D$, and the child therefore lies in $\inte D$ and has diameter less than $2^{-j}$. The same
procedure is applied in every parent. The inverse stereographic projection sends the resulting
Euclidean balls to round balls, and uniform continuity gives diameters tending to zero on
$\mathbb{S}^d$.
\end{proof}

\begin{cor}\label{simply-connected}
Suppose $d\ge3$ and $K\subset\mathbb{S}^d$ is a compact set such that $\mathcal{H}^1(K)=0$. Then
$\mathbb{S}^d\setminus K$ is path connected and simply connected. In particular, this conclusion
holds when $\dim_H K<1$.
\end{cor}
\begin{proof}
For $K=\emptyset$ there is nothing to prove, so let $K\ne\emptyset$. For a compact set
$A\subset\mathbb{S}^d\setminus K$, set $a=\dis(A,K)>0$ and choose, by Corollary \ref{nested}, a
generation $D_1,...,D_m$ such that
\[
\max_i\diam D_i<a,\qquad D_i\cap K\ne\emptyset,
\]
so that $D_i\cap A=\emptyset$ for all $i$, and set
\[
M_j=\mathbb{S}^d\setminus\bigcup_{i=j+1}^m\inte D_i \quad(0\le j\le m),\qquad M=M_0,\qquad
M_m=\mathbb{S}^d.
\]
For path connectedness, let $x,y\in M_j$ and choose
\[
z\in\mathbb{S}^d\setminus\left(\{-x,-y\}\cup\bigcup_{i=j+1}^m\partial D_i\right)
\]
and join $x$ to $z$ to $y$ by two shorter great-circle arcs. Neither arc lies in a boundary sphere,
so each meets each $\partial D_i$ at most twice, giving at most $4(m-j)$ crossing points, and each
segment inside $D_i$ can be replaced by a path in $\partial D_i\cong\mathbb{S}^{d-1}$. Since the
boundaries are path connected and mutually disjoint, the resulting path lies in $M_j$. With
$A=\{x,y\}\subset\mathbb{S}^d\setminus K$ this gives a path in $M\subset\mathbb{S}^d\setminus K$.

For a loop $\g:\mathbb{S}^1\to\mathbb{S}^d\setminus K$, take $A=\g(\mathbb{S}^1)$ and, for $1\le
j\le m$, a collar
\[
c_j:\mathbb{S}^{d-1}\times(-1,1)\longrightarrow M_j, \qquad
c_j(\mathbb{S}^{d-1}\times\{0\})=\partial D_j,
\]
whose negative side lies in $\inte D_j$, and denote by $C_j$ its image. The open sets
\[
U_j=M_{j-1}\cup C_j,\qquad V_j=\inte D_j\cup C_j
\]
cover $M_j$ and satisfy
\[
U_j\simeq M_{j-1},\qquad V_j\simeq D_j\simeq\{\mathrm{pt}\},\qquad U_j\cap
V_j=C_j\simeq\mathbb{S}^{d-1}.
\]
With base points suppressed, van Kampen's theorem gives, since $d\ge3$,
\[
\pi_1(M_j)\cong\pi_1(M_{j-1})*_{\pi_1(\mathbb{S}^{d-1})}\{1\} =\pi_1(M_{j-1}),\qquad
\pi_1(M)=\pi_1(\mathbb{S}^d)=0.
\]
Hence $\g$ contracts in $M\subset\mathbb{S}^d\setminus K$.
\end{proof}

\begin{thm}[Hou\cite{Hou2026}]\label{structural}
Let $\G<\Isom(\mathbb{H}^N)$ be discrete, torsion-free, non-elementary, and generated by $r$
elements. If $\mathfrak{D}_\G<1$, then $\G$ has geometric finiteness and is abstractly free. If $\G$
has no parabolic elements, it is convex cocompact. In particular, $\mathfrak{D}_\G<1/2$ implies that
$\G$ is free and convex cocompact. The same conclusions hold if $\delta_\G<1$, in which case
$\mathfrak{D}_\G=\delta_\G$, hence $\delta_\G<1/2$ implies freeness and convex cocompactness.
\end{thm}
This is the structural theorem of \cite{Hou2026} (Theorem 2.4 and Corollary 1.6) and \cite{Hou2026}
(Theorem 1.1 and Corollary 1.4(ii)). Since Hou states the theorem in \cite{Hou2026} in a much more
general setting than is needed here, we give a much shorter proof in our setting, for completeness.

For the proof we first recall three classical results and then a fourth one, which is another
theorem of Hou \cite{Hou2026}.
\begin{itemize}
\item A finitely generated Kleinian group with totally disconnected limit set has geometric
finiteness, by Gusevskii's result, see \cite{Kapovich}.
\item Such a group is isomorphic to a Schottky-type group, that is, a free product of elementary
groups, by Kulkarni's theorem, see \cite{Kapovich}.
\item For a non-elementary group, the critical exponent equals the Hausdorff dimension of the
conical limit set \cite{Hou2026} (Theorem 2.4), also recalled in \cite{Kapovich} (Theorem 10.15),
which implies $\delta_\G\le \mathfrak{D}_\G$. For a geometrically finite non-elementary Kleinian
group every nonconical limit point is parabolic, and since there are at most countably many such
points, $\delta_\G=\mathfrak{D}_\G$.
\item If a non-elementary Kleinian group generated by $r$ elements has $0<\delta_\G<\infty$ and
$P_\G(\sigma,o)<\infty$ for some $0<\sigma\le1$, then $\mathfrak{D}_\G=\delta_\G$ and
$\mathcal{H}^\sigma(\Lambda_\G)=0$, by Hou's Corollary 1.6 \cite{Hou2026}: for $\delta_\G<1$ the
corollary states the dimension equality, and if $\delta_\G=1$, convergence forces $\sigma=1$, so
vanishing of the measure implies $\mathfrak{D}_\G\le1=\delta_\G\le \mathfrak{D}_\G$.
\end{itemize}
The Gromov product of $x,y,o\in\mathbb{H}^N$ is
\[
(x\mid y)_o=\tfrac12\bigl(\dis(o,x)+\dis(o,y)-\dis(x,y)\bigr),
\]
and for $\xi,\eta\in\mathbb{S}^d$, $(\xi\mid\eta)_o$ is the limit along rays from $o$ to $\xi,\eta$.
The visual metric in Hou's cited corollary is the Bourdon metric of curvature $-1$, which in the
unit ball model centered at $o$ is
\[
\rho_o(\xi,\eta)=e^{-(\xi\mid\eta)_o}=\tfrac12|\xi-\eta|,
\]
by the hyperbolic cosine rule applied to the two rays. Note that this metric is a constant multiple
of chordal distance, so it is bi-Lipschitz equivalent to round distance, hence the dimension used in
\cite{Hou2026} and the dimension $\mathfrak{D}_\G$ used here give the same conclusions.

\begin{proof}
If a compact connected metric set $C$ contains two points $x\ne y$, the Lipschitz map
$z\mapsto\dis(x,z)$ sends $C$ onto an interval containing $[0,\dis(x,y)]$, which implies $\dim_H
C\ge1$, hence $\mathfrak{D}_\G<1$ implies that $\Lambda_\G$ is totally disconnected, and the first
of the results above gives geometric finiteness.

By the second, $\G$ is a free product $P_1*\cdots*P_q$ of elementary groups for some integer
$q\ge1$, where $P_i\le\G$ denotes the image of the $i$th factor for $1\le i\le q$. Each $P_i$ is
torsion-free and virtually abelian, and since a non-elementary Kleinian group contains a nonabelian
free subgroup, $P_i$ is elementary also in the given action.

Let us first assume that $P_i$ has two limit points. Its subgroup fixing both points has a
translation homomorphism to $\mathbb{R}$ with compact rotational factor. Since the preimage of a
bounded closed translation interval is compact, discreteness makes the translation image discrete,
and the kernel, being compact and discrete, is trivial by torsion-freeness. An element swapping the
endpoints would fix a point on the axis, which is excluded as well, hence $P_i$ is infinite cyclic.

An elementary parabolic factor is virtually abelian, of some abelian rank $k$, the rank of a free
abelian subgroup of finite index, and it will be shown below that $\delta(P_i)\ge k/2$. Since
\[
\delta(P_i)\le\delta_\G\le \mathfrak{D}_\G<1,
\]
this forces $k=1$, and a torsion-free virtually cyclic group is infinite cyclic. An elementary
factor with empty limit set fixes an interior point, whose stabilizer is compact, so discreteness
and torsion-freeness make the factor trivial, hence every nontrivial factor is infinite cyclic, and
$\G$ is free.

For the parabolic estimate, let $p$ be a nontrivial parabolic element, with its fixed point at
infinity in the upper-half-space model, so that
\[
p(x,t)=(Ax+b,t),\qquad A\in O(N-1),\quad b\in\mathbb{R}^{N-1}.
\]
With $o=(0,1)$ and $p^no=(b_n,1)$ for $n\in\mathbb{Z}$, by composition and inversion we have
\[
b_0=0,\qquad b_n=\sum_{j=0}^{n-1}A^jb\quad(n\ge1),\qquad b_{-n}=-A^{-n}b_n\quad(n\ge1),
\]
since $p^n(x,t)=(A^nx+b_n,t)$, whose inverse translates by $-A^{-n}b_n$, hence $|b_n|\le |n|\,|b|$,
and for a constant $C_p\ge0$ independent of $n$,
\[
\cosh\dis(o,p^n o)=1+|b_n|^2/2, \qquad \dis(o,p^n o)\le 2\log(1+|n|)+C_p.
\]
Comparison with $\sum_{n\ge1}n^{-2s}$ for $s>0$ implies $\delta(<p>)\ge1/2$.

For a subgroup $H\cong\mathbb{Z}^k$ of finite index in an elementary parabolic group, set the common
fixed point at infinity, take $o=(0,1)$ in upper half-space, choose a basis $p_1,...,p_k$ and, for
$v=(v_1,...,v_k)\in\mathbb{Z}^k$, write $p_v=p_1^{v_1}\cdots p_k^{v_k}$, $\|v\|_1=\sum_{i=1}^k|v_i|$
and $p_vo=(b_v,1)$. Since the horizontal linear parts of the generators are orthogonal, there are
constants $C,C'>0$, depending only on the chosen basis and base point, with
\[
|b_v|\le C(1+\|v\|_1),\qquad \dis(o,p_v o)\le2\log(1+\|v\|_1)+C'.
\]
For integers $j\ge0$ the shell count is
\[
\#\{v\in\mathbb{Z}^k:2^j\le\|v\|_1<2^{j+1}\}\asymp2^{jk},
\]
with constants depending only on $k$, so the Poincar\'e sum is bounded below by a constant times
$\sum_{j=0}^\infty2^{j(k-2s)}$, which diverges for $2s\le k$, and the exponent is therefore at least
$k/2$.

Since $\delta_\G\le \mathfrak{D}_\G$, the bound $\mathfrak{D}_\G<1/2$ excludes parabolic elements,
and geometric finiteness without parabolics is equivalent to convex cocompactness \cite{Kapovich}
(Theorem 2.7).

Now let $\delta_\G<1$. By non-elementarity there are independent loxodromic elements, and Theorem
\ref{powers} provides powers $a,b$ of them freely generating the subgroup $\langle a,b\rangle$. With
$M=\max\{\dis(o,ao),\dis(o,bo)\}$, this subgroup has, for every integer $j\ge1$, $4\cdot3^{j-1}$
reduced words of length $j$, each with displacement at most $jM$, so
\[
\delta_\G\ge\frac{\log3}{M}>0.
\]
For $\delta_\G<\sigma<1$ this means $P_\G(\sigma,o)<\infty$, and Hou's Corollary 1.6 \cite{Hou2026}
implies $\mathfrak{D}_\G=\delta_\G<1$. The argument above then provides geometric finiteness and
freeness, and for $\delta_\G<1/2$ the parabolic estimate provides convex cocompactness as well.
Since $\delta_\G\le \mathfrak{D}_\G$, each Hausdorff-dimension bound in the main theorems implies
its critical-exponent bound.
\end{proof}

\begin{rem}
Neither $\mathfrak{D}_\G<1$ nor $\delta_\G<1$ excludes parabolic elements, and $\mathfrak{D}_\G<1/2$
or $\delta_\G<1/2$ gives convex cocompactness. The Hausdorff-dimension versions of Theorems
\ref{main} and \ref{main-2} use the first three facts above, and the critical-exponent versions also
use Hou's Corollary 1.6 to get $\delta_\G=\mathfrak{D}_\G$. If torsion is allowed, Selberg's lemma
\cite{Kapovich} (Section 2) gives a finite index torsion-free subgroup $\G_0$, finitely generated by
Schreier's lemma and with $\Lambda_{\G_0}=\Lambda_\G$, so Theorem \ref{structural} makes $\G_0$ free
if $\mathfrak{D}_\G<1$ or $\delta_\G<1$, and $\G$ is virtually free.
\end{rem}

We identify $\mathbb{S}^d$ with the unit sphere in $\mathbb{R}^{d+1}$, with chordal distance
$|x-y|$. For a loxodromic transformation $g_i$, possibly with a rotation about its axis, let $a_i$
be its attracting point, $b_i$ its repelling point and $\ell_i>0$ its translation length.

Hou's concentric circles, obtained after sending the fixed points to zero and infinity
\cite{Hou2020} (Proposition 5.5), have round-sphere counterparts in every dimension, and the ratio
identity below keeps track of the orthogonal part of the transformation. Theorem \ref{criterion} is
a quantitative version of the principle in Hou's classicality criteria that contraction must
dominate the separation scale of the fixed points and paired boundaries \cite{Hou2020} (Lemma 5.2).

\begin{lem}\label{ratio}
Set
\[
f_i(x)=\frac{|x-a_i|}{|x-b_i|},\quad f_i(a_i)=0,\quad f_i(b_i)=\infty.
\]
Then we have
\begin{equation}\label{ratio-eq}
f_i(g_i x)=e^{-\ell_i}f_i(x).
\end{equation}
For all $t>0$, the sets $\{f_i\le t\}$ and $\{f_i\ge t\}$ are closed round balls.
\end{lem}
\begin{proof}
For $x\ne b_i$, set
\[
T_i(x)=\frac{x-b_i}{|x-b_i|^2} -\frac{a_i-b_i}{|a_i-b_i|^2}.
\]
Since $b_i\cdot(x-b_i)/|x-b_i|^2=-1/2$, the image lies in the Euclidean hyperplane
$\{y\in\mathbb{R}^{d+1}:b_i\cdot y=0\}$, identified isometrically with $\mathbb{R}^d$, and this
inversion followed by a translation implies $T_i(a_i)=0$ and $T_i(b_i)=\infty$. A computation shows
that
\[
|T_i(x)|^2 =\frac{|x-b_i|^2+|a_i-b_i|^2-2(x-b_i)\cdot(a_i-b_i)} {|x-b_i|^2|a_i-b_i|^2}
=\frac{|x-a_i|^2}{|x-b_i|^2|a_i-b_i|^2},
\]
so $|T_i(x)|=c_if_i(x)$ with $c_i=|a_i-b_i|^{-1}>0$. Since the conjugate fixes zero and infinity and
has zero as an attracting point,
\[
T_i g_i T_i^{-1}(z)=e^{-\ell_i}A_i z,\qquad A_i\in O(d),
\]
and taking norms implies (\ref{ratio-eq}), including its endpoint limits.

Since $|x|=|a_i|=|b_i|=1$, the inequality can be squared,
\[
|x-a_i|\le t|x-b_i| \iff x\mathbin{\cdot}(a_i-t^2b_i)\ge1-t^2,
\]
and
\[
|a_i-t^2b_i|^2-(1-t^2)^2=t^2|a_i-b_i|^2>0,
\]
so the normal is nonzero, the boundary is a round sphere and the inequality selects a round ball.
Interchanging the endpoints handles the reverse inequality.
\end{proof}

\begin{thm}\label{criterion}
Let $r\ge1$ be an integer, and let $g_1,...,g_r$ be loxodromic transformations whose $2r$ endpoints
are pairwise distinct. Set
\[
\varepsilon=\min\{|u-v|:u\ne v,\ u,v\in\{a_i,b_i:1\le i\le r\}\}, \qquad L=\min_i\ell_i.
\]
If
\begin{equation}\label{criterion-eq}
L>2\log(4/\varepsilon),
\end{equation}
then $<g_1,...,g_r>$ is a classical round-ball Schottky group, with paired balls
\[q_i=e^{-\ell_i/2},\qquad B_i^+=\{f_i\le q_i\},\qquad B_i^-=\{f_i\ge q_i^{-1}\}.
\]
\end{thm}
\begin{proof}
By Lemma \ref{ratio} the balls are round, and (\ref{ratio-eq}) gives the pairing identity
\[
g_i\{f_i\le q_i^{-1}\}=\{f_i\le q_i\}.
\]
Since the sphere has chordal diameter two,
\[
\begin{aligned}
x\in B_i^+&\Longrightarrow |x-a_i|\le q_i|x-b_i|\le2q_i,\\
x\in B_i^-&\Longrightarrow |x-b_i|\le q_i|x-a_i|\le2q_i,
\end{aligned}
\]
and $2q_i<\varepsilon/2$ by (\ref{criterion-eq}), so all the balls are disjoint by the triangle
inequality.

Being proper, closed and positively separated, the balls cannot cover the connected sphere. Choose
$z$ outside them and write $B_{g_i}=B_i^+$ and $B_{g_i^{-1}}=B_i^-$. For a nonempty reduced word
$w=s_1\cdots s_k$, by the pairing identities applied successively we have
\[
wz\in\inte B_{s_1},
\]
since at each step reducedness keeps the point outside the next repelling ball. It follows that
$wz\ne z$, the generators are free, and by Proposition \ref{gap} the group is discrete with the
stated round fundamental domain, hence classical.
\end{proof}

\begin{prop}\label{gap}
Let $r\ge1$ be an integer. Suppose the loxodromic transformations $g_1,...,g_r$ admit paired
disjoint closed round balls. Set $\G=<g_1,...,g_r>$, let $\mathcal A=\{g_i,g_i^{-1}:1\le i\le r\}$
be the $2r$ signed letters, and write $B_{g_i}=B_i^+$ and $B_{g_i^{-1}}=B_i^-$. For a reduced word
$w$ in this alphabet, $|w|$ denotes its number of letters. For each of the $2r$ signed generators
let $E_u\subset\mathbb{H}^N$ be the closed halfspace with ideal ball $B_u$, choose $o$ outside all
the halfspaces, and set
\[
\beta=\min_u\dis(o,\partial E_u)>0, \qquad \alpha=\min_{u\ne v}\dis(\partial E_u,\partial E_v)>0.
\]
For every nonempty reduced word $w=u_1\cdots u_k$ with $k\ge1$,
\begin{equation}\label{displacement}
\dis(o,wo)\ge2\beta+(k-1)\alpha.
\end{equation}
The group is discrete and convex cocompact, the round complement is a fundamental domain on the
discontinuity set, and
\[\delta_\G\le\frac{\log(2r-1)}{\alpha}, \qquad \dim_H\Lambda_\G\le\frac{\log(2r-1)}{\alpha}.
\]
\end{prop}
\begin{proof}
Since the paired closed round balls are disjoint, the closed halfspaces are pairwise disjoint as
well, for if two of them had a common point, we could move that point to $(1,0,...,0)$ in the
hyperboloid model and write their unit spacelike normals as $n_j=(s_j,w_j)$, $j=1,2$, with $s_j\le0$
since the point lies in both halfspaces, and then, choosing $\xi\in\mathbb{S}^{N-1}$ perpendicular
to both $w_1,w_2$ in the Euclidean scalar product, which is possible since $N\ge4$,
\[
\langle(1,\xi),n_j\rangle=-s_j\ge0 \qquad(j=1,2),
\]
so that the future null ray through $(1,\xi)$ would belong to both closed ideal balls. The bounding
hyperplanes also have positive pairwise distance, since if points $x_k,y_k$ on two of them had
$\dis(x_k,y_k)\to0$, then either a bounded subsequence would converge to a point of $\mathbb{H}^N$
lying on both hyperplanes, or a subsequence would converge to an ideal point in the closed
Poincar\'e ball model, where $|x_k-y_k|\le\tfrac12\dis(x_k,y_k)\longrightarrow0$ would give both
sequences the same ideal limit, belonging to both closed ideal balls. Taking the minimum over the
pairs implies $\alpha>0$. Finally, the $2r$ disjoint proper closed ideal balls do not cover the
connected sphere, so an ideal point outside their union can be chosen, and $o\in\mathbb{H}^N$
sufficiently close to it lies outside every closed halfspace and at positive distance from each of
the $2r$ bounding hyperplanes, which implies $\beta>0$.

For $w=u_1\cdots u_k$, set
\[
F_j=u_1\cdots u_{j-1}E_{u_j},\qquad 1\le j\le k.
\]
Since the ideal pairing determines the bounding plane and its side, for all $u\in\mathcal{A}$
\[
u\bigl(\mathbb{H}^N\setminus\inte E_{u^{-1}}\bigr)=E_u, \qquad u\bigl(\mathbb{H}^N\setminus
E_{u^{-1}}\bigr)=\inte E_u.
\]
Since $u_{j+1}\ne u_j^{-1}$,
\[
E_{u_{j+1}}\subset\mathbb{H}^N\setminus E_{u_j^{-1}}, \qquad F_{j+1}\subset\inte F_j,
\]
and applying $(u_1\cdots u_j)^{-1}$,
\[
\dis(\partial F_j,\partial F_{j+1}) =\dis(\partial E_{u_j^{-1}},\partial E_{u_{j+1}})\ge\alpha.
\]
Note also that $wo\in\inte F_k$ and
\[
\dis(wo,\partial F_k)=\dis(o,\partial E_{u_k^{-1}})\ge\beta.
\]
The geodesic from $o$ to $wo$ crosses these planes in order, its first and last segments have length
at least $\beta$, and every segment between them has length at least $\alpha$, which implies
(\ref{displacement}). Note that this positive lower bound excludes the identity for every nonempty
formally reduced word, which implies that $g_1,...,g_r$ freely generate $\G$, and since every
nonidentity element moves $o$ by at least $2\beta$, the group is discrete.

Counting the $2r(2r-1)^{k-1}$ reduced words of length $k$,
\[
P_\G(s,o)\le1+2r e^{-2s\beta} \sum_{k\ge1}\bigl((2r-1)e^{-s\alpha}\bigr)^{k-1},
\]
so the series converges if $(2r-1)e^{-s\alpha}<1$, which is the exponent bound.

Let $B_w$ be the ideal ball of $F_k$. The same crossings show $\dis(o,\partial
F_k)\ge\beta+(k-1)\alpha$, and in the ball model centered at $o$ a halfspace at distance $T$ has cap
radius $\theta$ with $\sin\theta=1/\cosh T$, so
\[
\diam B_w=2\sin\theta\le4e^{-T} \le4e^{-\beta-(k-1)\alpha}.
\]
With $T_w=\dis(o,\partial F_k)$ and $v_w\in\mathbb{S}^d$ the center of $B_w$, in the Poincar\'e ball
model centered at $o$
\[
F_k=\{x\in\mathbb{R}^N:|x|<1,\ |x-\coth(T_w)v_w|\le\operatorname{csch}(T_w)\},
\]
since the bounding Euclidean sphere is orthogonal to $\mathbb{S}^d$ and its nearest point to $o$, at
hyperbolic distance $T_w$ from $o$, has norm $\coth T_w-\operatorname{csch}T_w=\tanh(T_w/2)$. With
Euclidean closure and diameter,
\[
\diam_{\mathrm{Eucl}}\overline{F_k}\le2\operatorname{csch}T_w\longrightarrow0, \qquad \inf_{x\in
F_k}|x|=\tanh(T_w/2)\longrightarrow1,
\]
since $T_w\ge\beta+(k-1)\alpha$.

For each $k$ the balls $B_w$ cover $\Lambda_\G$: if distinct orbit points approach a limit point,
their word lengths tend to infinity, and after passing to a subsequence the first $k$ letters of
their reduced words are the same from some point onward, for each $k\ge1$, so the limit belongs to
the corresponding nested balls. Conversely, if $w_k$ consists of the first $k$ letters of an
infinite reduced sequence, the compact Euclidean closures of the corresponding halfspaces are nested
with diameters tending to zero, the radial estimate places their unique common point on
$\mathbb{S}^d$, and since $w_ko$ lies in the halfspace corresponding to $w_k$, the sequence $(w_ko)$
converges to this point.

Since the balls $B_w$ with $|w|=k$ cover the limit set and
\[
\sum_{|w|=k}(\diam B_w)^s \le 2r\,4^s e^{-s\beta} \bigl((2r-1)e^{-s\alpha}\bigr)^{k-1}
\longrightarrow0
\]
whenever $(2r-1)e^{-s\alpha}<1$, the dimension bound follows, and changing the ball model preserves
Hausdorff dimension, since the induced boundary mapping is bi-Lipschitz.

Set $D_\infty=\mathbb{S}^d\setminus\bigcup_u\inte B_u$ and start with a boundary point $x_0$.
Whenever $x_j\in\inte B_{u_{j+1}}$, let
\[
x_{j+1}=u_{j+1}^{-1}x_j.
\]
By the pairing,
\[
u^{-1}(\inte B_u)=\mathbb{S}^d\setminus B_{u^{-1}},
\]
so consecutive chosen labels satisfy $u_{j+2}\ne u_{j+1}^{-1}$. If the reduction does not stop, it
follows an infinite reduced word, and $x_0$ belongs to all its nested balls and hence to
$\Lambda_\G$. It follows that every point outside $\Lambda_\G$ reduces to $D_\infty$. For a nonempty
reduced word $w=u_1\cdots u_k$, the same pairings give $w(\inte D_\infty)\subset\inte B_{u_1}$, so
the translated interiors are disjoint, and the boundary identifications are those of
(\ref{pairing}).

For $D=\mathbb{H}^N\setminus\bigcup_u\inte E_u$, the same recurrence applies whenever $x_j\in\inte
E_{u_{j+1}}$, the next point lying outside $E_{u_{j+1}^{-1}}$. If the reduction did not stop, $x_0$
would belong to every $F_j$ along an infinite reduced word, and then
\[
\dis(o,\partial F_j)\ge\beta+(j-1)\alpha\longrightarrow\infty,
\]
which is impossible for a fixed point of $\mathbb{H}^N$. It follows that the translates of $D$ cover
$\mathbb{H}^N$, and since also $w(\inte D)\subset\inte E_{u_1}$ for every nonempty reduced $w$,
their interiors are disjoint and their side pairings are the same, which implies that $D$ is a
fundamental region. Since the depth-two balls cover $\Lambda_\G$ and lie strictly inside the
original balls,
\[
\Lambda_\G\cap\partial_\infty D=\emptyset.
\]
Note that $\partial_\infty\Hu(\Lambda_\G)\subset\Lambda_\G$, since a small round cap disjoint from
$\Lambda_\G$ determines a halfspace excluded from its convex hull, hence $\Hu(\Lambda_\G)\cap D$ has
no ideal accumulation point, so, being closed and bounded, it is compact, and it covers the
convex-core quotient, which proves convex cocompactness.
\end{proof}
For a real matrix $M$ we denote by $M^{\mathsf T}$ its transpose, for each integer $k\ge1$ we set
$\mathbf I_k=(1,...,1)^{\mathsf T}\in\mathbb{R}^k$, and inequalities between real vectors or
matrices are understood entrywise. Under the conditions and with the signed labels of Proposition
\ref{gap}, for $s>0$ let $Q_s$ be the $2r$ by $2r$ matrix indexed by $u,v\in\mathcal{A}$ with
entries
\[
Q_s(u,v)=
\begin{cases}
\exp\bigl(-s\dis(\partial E_{u^{-1}},\partial E_v)\bigr),&v\ne u^{-1},\\
0,&v=u^{-1}.
\end{cases}
\]
For a nonempty reduced word $w=u_1\cdots u_k$, set $F_k=u_1\cdots u_{k-1}E_{u_k}$. Since the planes
are ordered,
\[
\begin{aligned}
\dis(o,wo)&\ge2\beta+ \sum_{j=1}^{k-1}\dis(\partial E_{u_j^{-1}},\partial E_{u_{j+1}}),\\
\dis(o,\partial F_k)&\ge\beta+ \sum_{j=1}^{k-1}\dis(\partial E_{u_j^{-1}},\partial E_{u_{j+1}}),
\end{aligned}
\]
and the cap bound $\diam B_w\le4e^{-\dis(o,\partial F_k)}$ combined with matrix multiplication gives
\begin{equation}\label{weighted}
\sum_{|w|=k}(\diam B_w)^s \le4^s e^{-s\beta}\,\mathbf I_{2r}^{\mathsf T}Q_s^{\,k-1}\mathbf I_{2r}.
\end{equation}
If the spectral radius of $Q_s$ is less than one, its powers tend to zero and their series
converges, so by (\ref{weighted}) the total $s$-power diameter tends to zero. Since every allowed
gap is at least $\alpha>0$, the maximum diameter also tends to zero, and $\dim_H\Lambda_\G\le s$.
Similarly, by the displacement bound we have
\[
P_\G(s,o)\le1+e^{-2s\beta} \sum_{k\ge1}\mathbf I_{2r}^{\mathsf T}Q_s^{\,k-1}\mathbf I_{2r}<\infty,
\]
so $\delta_\G\le s$.

Now choose nonempty pairwise disjoint sets $I_1,...,I_b$ whose union is the set of signed labels,
and suppose $W_s\in[0,\infty)^{b\times b}$ satisfies
\[
\sum_{v\in I_j}Q_s(u,v)\le(W_s)_{ij} \qquad(u\in I_i,\ 1\le i,j\le b).
\]
Induction on $n$ shows $(Q_s^n\mathbf I_{2r})_u\le(W_s^n\mathbf I_b)_i$ for $u\in I_i$, and summing,
\[
\mathbf I_{2r}^{\mathsf T}Q_s^{\,k-1}\mathbf I_{2r} \le\sum_{i=1}^b |I_i|(W_s^{\,k-1}\mathbf I_b)_i,
\]
so a spectral radius of $W_s$ less than one suffices for both dimension bounds.

\begin{thm}\label{powers}
Let $r\ge1$ be an integer, and let $g_1,...,g_r$ be loxodromic elements of a Kleinian group, with
pairwise distinct attracting and repelling points. We set $H_m=\langle g_1^m,...,g_r^m\rangle$. For
all sufficiently large integers $m$, the group $H_m$ is classical round-ball Schottky, and
\[
\lim_{m\to\infty}\dim_H\Lambda_{H_m}=0, \qquad \lim_{m\to\infty}\delta(H_m)=0.
\]
\end{thm}
\begin{proof}
Elementary.
\end{proof}

\begin{thm}\label{three}
Let $\G<\Isom(\mathbb{H}^N)$ be discrete, torsion-free, non-elementary, generated by $r$ elements,
and parabolic-free. Suppose also that it preserves a totally geodesic $Y\cong\mathbb{H}^3$, that its
restriction to $Y$ preserves orientation, and that $\mathfrak{D}_\G<1$. Then $\G$ is classical
round-ball Schottky in $\mathbb{H}^N$.
\end{thm}
\begin{proof}
The stabilizer of $Y$ is $\Isom(Y)$ times a compact transverse orthogonal group. If the restriction
image were not discrete, infinitely many distinct images in a compact set would lift to infinitely
many elements of $\G$ in a compact set, contradicting discreteness. So the image is discrete, and
the kernel, being compact and discrete, is trivial by torsion-freeness.

For $o\in Y$ the orbit stays in $Y$, hence $\Lambda_\G\subset\partial_\infty Y$, and the induced
round metrics are bi-Lipschitz, so the Hausdorff dimension is the same. Note that every nonidentity
restriction is loxodromic, since an elliptic restriction would give an interior fixed point of the
full isometry and a parabolic restriction would make the full isometry parabolic. By Theorem
\ref{structural}, $\G$ is convex cocompact, and its convex hull lies in $Y$ and agrees with the
intrinsic hull, so the restricted convex-core quotient is the same compact quotient.

Hou's three-dimensional theorem \cite{Hou} provides generators and disjoint paired round disks on
$\partial_\infty Y$. In the hyperboloid model, write
\[
\mathcal{V}=W\oplus W^{\perp},
\]
where $W$ has three positive directions and one negative direction and determines $Y$, and extend
each bounding plane by the same spacelike normal $n\in W$. For a future timelike vector $x=w+u$,
with $w\in W$ and $u\in W^\perp$,
\[
\langle w,w\rangle=\langle x,x\rangle-\langle u,u\rangle<0,\qquad \langle x,n\rangle=\langle
w,n\rangle.
\]
With the time orientation on $W$ inherited from $\mathcal{V}$, normalizing $w$ produces a point of
$Y$ with the same halfspace signs, so the extended halfspaces remain disjoint. For a nonzero future
null vector $x$, the projection is timelike when $u\ne0$ and null when $u=0$, and the same sign
argument applied to an interior or ideal point of $Y$ gives disjoint extended ideal closed balls.
Since every $\g\in\G$ preserves $W$ and $W^\perp$, its restriction determines its action on the
normals and signed halfspaces, so the side pairings extend to (\ref{pairing}) on the larger sphere,
and since restriction is faithful, the chosen generators generate $\G$.
\end{proof}

\begin{rem}
Theorem \ref{three} allows arbitrary transverse rotations. If the restriction reverses orientation,
the argument applies to the subgroup whose restriction preserves orientation, which has index at
most two.
\end{rem}


\section{Dimension, Pressure, and Cyclic Geometry}
In this section we state the displacement estimates that will be used for the generator selections
in the rest of the paper. Throughout this section, distances and exponents refer to the original
group $\G$. The displacement inequality (\ref{pressure-eq}) is due to Hou \cite{Hou2001} and is
restated for the three-dimensional setting in \cite{Hou2020} (Theorem 3.1), and its two-generator
form (\ref{pair-pressure}) is equivalent to Hou's Corollary 3.2 there. The free-product proof below
works in every dimension and uses the exponent of the containing group.
\begin{lem}\label{pressure}
Let $g_1,...,g_r$ freely generate a subgroup of $\G$, with $r\ge2$ and $0<\delta_\G<\infty$. For
$x\in\mathbb{H}^N$, set $D_i=\dis(x,g_i x)$. Then
\begin{equation}\label{pressure-eq}
\sum_{i=1}^r\frac{2}{1+e^{\delta_\G D_i}}\le1.
\end{equation}
In particular, if $\G$ is free and $a,b\in\G$ do not commute,
\begin{equation}\label{pair-pressure}
(e^{\delta_\G\dis(x,ax)}-1)(e^{\delta_\G\dis(x,bx)}-1)\ge4.
\end{equation}
\end{lem}
\begin{proof}
For $s>\delta_\G$, the subseries
\[
Z_i(s)=\sum_{n\in\mathbb{Z}\setminus\{0\}}e^{-s\dis(x,g_i^n x)}
\]
converge. By free-product normal form, all words $g_{i_1}^{n_1}\cdots g_{i_k}^{n_k}$ with $n_j\ne0$
and $i_j\ne i_{j+1}$ are distinct, and since each displacement is at most the sum of its block
displacements,
\[
P_\G(s,x)\ge 1+\sum_{k\ge1}\ \sum_{\substack{i_1,...,i_k\in\{1,...,r\}\\
i_j\ne i_{j+1}\ (1\le j<k)}} Z_{i_1}(s)\cdots Z_{i_k}(s).
\]
For fixed $s$, abbreviate $Z_i=Z_i(s)$ and set $M_{ij}=Z_j$ for $i\ne j$, $M_{ii}=0$ and $\mathbf
z=(Z_1,...,Z_r)^{\mathsf T}$, hence the length-$k$ block sum is $\mathbf z^{\mathsf T}M^{k-1}\mathbf
I_r$. Since $r\ge2$ and the off-diagonal entries are positive, $M$ is irreducible. For its Perron
eigenvector $v>0$ and eigenvalue $\lambda$, choose $c>0$ with $\mathbf I_r\ge cv$ component-wise, so
that
\[
\mathbf z^{\mathsf T}M^{k-1}\mathbf I_r \ge c\lambda^{k-1}\mathbf z^{\mathsf T}v.
\]
Since the series converges, $\lambda<1$. With $A=\sum_jZ_jv_j$, the eigenvector equations yield
\[
(\lambda+Z_i)v_i=A,\qquad \sum_i\frac{Z_i}{\lambda+Z_i}=1, \qquad \sum_i\frac{Z_i}{1+Z_i}<1,
\]
and from $\dis(x,g_i^n x)\le|n|D_i$
\[
Z_i(s)\ge\frac{2}{e^{sD_i}-1}, \qquad \sum_i\frac{2}{1+e^{sD_i}}<1.
\]
Letting $s\downarrow\delta_\G$ implies (\ref{pressure-eq}).

For noncommuting $a,b$, the subgroup $H=<a,b>$ is free by the Nielsen--Schreier theorem. Its
abelianization is generated by two elements, so its rank is at most two, and since the elements do
not commute, the rank is two. The surjection $F_2\to H\cong F_2$ from the free group on two
generators sending its basis to $a,b$ is an isomorphism by Hopficity, which implies that $a,b$
freely generate $H$, and the argument above with $r=2$ and the exponent of $\G$ implies
(\ref{pair-pressure}) after clearing denominators in (\ref{pressure-eq}).
\end{proof}

\begin{prop}\label{thick}
Suppose $\G$ is a non-elementary free convex cocompact group. If
\[
0<R<\frac{\log3}{\delta_\G},
\]
then for all $x\in\mathbb{H}^N$, all elements of $\G$ moving $x$ by less than $R$ lie in a single
cyclic subgroup. There also exists $x\in\Hu(\Lambda_\G)$ such that
\[
\dis(x,\g x)\ge R\qquad(\g\ne1).
\]
Hence the largest injectivity radius in the convex-core quotient is at least $\log3/(2\delta_\G)$.
\end{prop}
\begin{proof}
If two elements $a,b\in\G$ with $\dis(x,ax),\dis(x,bx)<R$ did not commute, by (\ref{pair-pressure})
we would have
\[
4\le(e^{\delta_\G\dis(x,ax)}-1)(e^{\delta_\G\dis(x,bx)}-1) <(e^{\delta_\G R}-1)^2<4,
\]
a contradiction, so all elements $\g\in\G$ with $\dis(x,\g x)<R$ commute with one another and lie in
the cyclic centralizer of any nonidentity element among them.

For each maximal cyclic subgroup $C$, let
\[
U_C=\{x\in\Hu(\Lambda_\G): \dis(x,cx)<R\text{ for some }c\in C\setminus\{1\}\}.
\]
These sets are relatively open and pairwise disjoint, since commuting nontrivial elements belong to
one maximal cyclic subgroup. If they covered the hull, its connectedness would make one $U_C$ the
whole hull, and by equivariance
\[
\g U_C=U_{\g C\g^{-1}},
\]
so every $\g$ would then normalize $C$.

In a Cayley tree of $\G$ for a free basis, this normalizer preserves the axis of a generator of $C$,
and since the action is free on vertices and has no edge inversions, restriction to the axis is
injective and orientation preserving, and its integer translations form a cyclic group, which equals
$C$ by maximality, hence the normalizer is $C$, contrary to non-elementarity, and there is a point
outside all $U_C$.

The injectivity radius is half the least nonidentity displacement. By the triangle inequality,
\[
\bigl|\dis(x,\g x)-\dis(y,\g y)\bigr| \le2\dis(x,y),
\]
so, taking the infimum, the injectivity radius is $1$-Lipschitz. Being $\G$-invariant, it attains a
maximum on the compact convex-core quotient, and this maximum is at least $R/2$ for all
$R<\log3/\delta_\G$. Letting $R\uparrow\log3/\delta_\G$ proves the proposition.
\end{proof}

\begin{prop}\label{tubes}
Let $\G$ be free and convex cocompact, let $a\in\G$ have axis $A$ and translation length $\ell>0$,
let $b\in\G$ lie outside the axis stabilizer, and set $\delta=\delta_\G>0$ and $D=\dis(A,bA)$. Then
\begin{equation}\label{tube-axis}
D\ge \frac{1}{2\delta}\log\left(1+\frac{4}{e^{\delta\ell}-1}\right) -\frac{\ell}{2}.
\end{equation}
If $a$ generates its maximal cyclic subgroup and the right side is positive, the image of $A$ in
$\mathbb{H}^N/\G$ has an embedded open tube of every radius
\begin{equation}\label{tube-radius}
0<r< \frac{1}{4\delta}\log\left(1+\frac{4}{e^{\delta\ell}-1}\right) -\frac{\ell}{4}.
\end{equation}
\end{prop}
\begin{proof}
The elements $a$ and $bab^{-1}$ do not commute, since their axes do not agree. For $x\in A$ and
$y\in bA$,
\[
\dis(x,ax)=\ell,\qquad \dis(x,bab^{-1}x)\le2\dis(x,y)+\ell,
\]
so by (\ref{pair-pressure}) with $\dis(x,y)\downarrow D$ we have
\[
(e^{\delta\ell}-1)(e^{\delta(2D+\ell)}-1)\ge4,
\]
which is (\ref{tube-axis}).

An element of $\G$ reversing $A$ would fix a point of $A$, contrary to discreteness and
torsion-freeness. The subgroup of $\Isom(\mathbb{H}^N)$ preserving $A$ and its orientation is
$\mathbb{R}\times O(N-1)$. Restricting its translation projection to the stabilizer of $A$ in $\G$,
the preimage of a bounded closed translation interval is compact, so discreteness makes the image
discrete, and the kernel, being compact and discrete, consists of torsion elements and is trivial,
hence the axis stabilizer is cyclic, and it is the maximal cyclic subgroup containing $a$, since
every element commuting with $a$ preserves its axis.

The open $r$-neighborhoods of distinct translated axes are disjoint when the distance between the
axes exceeds $2r$, which by (\ref{tube-axis}) holds for every translate distinct from $A$ whenever
(\ref{tube-radius}) holds. It follows that the tube modulo its cyclic stabilizer embeds in
$\mathbb{H}^N/\G$, and if $a$ generates the stabilizer, this is the tube about the closed geodesic
associated with $a$.
\end{proof}

The endpoint stabilizers of a loxodromic element $a$ of a free convex cocompact group $\G$ are found
by normalizing
\[
a(x,t)=(e^LRx,e^Lt),\qquad L>0,\quad R\in O(N-1).
\]
If $g\in\G$ fixes infinity, then $g(x,t)=(\lambda Sx+c,\lambda t)$ with $\lambda>0$ and $S\in
O(N-1)$, hence for $n\ge0$
\[
a^{-n}ga^n(x,t)= (\lambda R^{-n}SR^nx+e^{-nL}R^{-n}c,\lambda t).
\]
These conjugates lie in a compact subset of the isometry group, so by discreteness there are
$n>m\ge0$ with $a^{-n}ga^n=a^{-m}ga^m$, and then $g$ commutes with $a^{n-m}$, preserves
$\{0,\infty\}$ and fixes both points. Interchanging the endpoints gives the same conclusion for an
element fixing zero, so by the axis-stabilizer argument in Proposition \ref{tubes} both endpoint
stabilizers are the same cyclic axis stabilizer, and if $a$ belongs to a free basis, it is not a
proper power, so this stabilizer is $<a>$. In particular, loxodromic elements sharing an endpoint
lie in one cyclic subgroup.


\section{Classical in Rank Two}
In this section we prove Theorem \ref{main-2}. The rank two theorem is proved by combining pressure
with an explicit hyperplane pairing for a given generating pair, and then changing generators by
Nielsen moves. As in Hou's generator-selection strategy, displacement estimates are combined with
control of fixed points and circle centers \cite{Hou2020} (Sections 3--6), and the perpendicular-cut
estimates and the terminating length reduction below allow arbitrary orthogonal parts. For
loxodromic $a,b$, write $\ell_a=\ell(a)$ and $\ell_b=\ell(b)$.

\begin{lem}\label{cuts}
Suppose $a,b$ are loxodromic transformations whose axes have positive distance $D$ and whose
translation lengths are $\ell_a,\ell_b>0$, and set
\[T(\ell_a,\ell_b)= \log\coth(\ell_a/4)+\log\coth(\ell_b/4).
\]
If $D>T(\ell_a,\ell_b)$, then the pair $a,b$ admits four disjoint paired round balls for arbitrary
rotational parts.
\end{lem}
\begin{proof}
Let $\langle x,y\rangle=-x_0y_0+\sum_{j=1}^Nx_jy_j$ and $T=T(\ell_a,\ell_b)$. At the
common-perpendicular feet $p,q$, let $u,v$ be the unit tangents pointing in the translation
directions of $a,b$, and set $A=\ell_a/2$, $B=\ell_b/2$, so that
\[
\langle p,q\rangle=-\cosh D,\qquad u,v\perp\operatorname{span}(p,q),\qquad |\langle u,v\rangle|\le1.
\]
The unit spacelike normals
\[
n_a^\pm=\sinh A\,p\pm\cosh A\,u,\qquad n_b^\pm=\sinh B\,q\pm\cosh B\,v,
\]
satisfy, for $\sigma,\tau\in\{+,-\}$,
\[
\begin{aligned}
\langle n_a^+,n_a^-\rangle&=-\cosh\ell_a,\qquad \langle n_b^+,n_b^-\rangle=-\cosh\ell_b,\\
\langle n_a^\sigma,n_b^\tau\rangle &\le-\sinh A\sinh B\cosh D+\cosh A\cosh B<-1,
\end{aligned}
\]
where the strict inequality follows from $D>T$ and the identities
\[
\begin{aligned}
\cosh(\log\coth(t/2))&=\coth t,\qquad \sinh(\log\coth(t/2))=\operatorname{csch}t,\\
\cosh T&=\frac{1+\cosh A\cosh B}{\sinh A\sinh B},
\end{aligned}
\]
of which the first line holds for $t>0$. Set
\[
\mathcal{I}=\{(a,+),(a,-),(b,+),(b,-)\},\qquad
n_{(g,\sigma)}=n_g^\sigma\quad((g,\sigma)\in\mathcal{I}).
\]
\[
S=\sum_{s\in\mathcal{I}}n_s =2\sinh A\,p+2\sinh B\,q.
\]
The vector $S$ is future timelike, and by the pair inequalities we have
\[
\langle S,n_s\rangle =1+\sum_{t\ne s}\langle n_t,n_s\rangle<-2 \qquad(s\in\mathcal{I}).
\]
Conjugating $S/\sqrt{-\langle S,S\rangle}$ to $(1,0,...,0)$, write
\[
n_s=(\sinh t_s,\cosh t_s\,z_s),\qquad t_s>0,\quad z_s\in\mathbb{S}^{N-1},\qquad
\theta_s=\arccos(\tanh t_s)\in(0,\pi/2).
\]
The positive ideal ball of $n_s$ is
\[
B_s=\{\xi\in\mathbb{S}^{N-1}:\xi\cdot z_s\ge\tanh t_s\}.
\]
For $s\ne t$,
\[
z_s\cdot z_t <\tanh t_s\tanh t_t-\operatorname{sech}t_s\operatorname{sech}t_t
=\cos(\theta_s+\theta_t),
\]
so $\angle(z_s,z_t)>\theta_s+\theta_t$ and the four closed balls are disjoint. The rotations fix the
Lorentz axis planes pointwise, so
\[
a n_a^-=-n_a^+,\qquad b n_b^-=-n_b^+,
\]
which are the side pairings.
\end{proof}

\begin{thm}\label{scalar}
Suppose $a,b$ freely generate a Kleinian group $\G$ with positive exponent $\delta$, and suppose
that both generators are loxodromic. If
\[\frac{1}{2\delta} \log\frac{4}{(e^{\delta\ell_a}-1)(e^{\delta\ell_b}-1)}
> \log\coth(\ell_a/4)+\log\coth(\ell_b/4),
\]
then the original pair admits a classical Schottky domain.
\end{thm}
\begin{proof}
We fix $p,q$ on the axes of $a,b$, respectively, set $d_0=\dis(p,q)$, and let
$D=\inf_{p,q}\dis(p,q)$ be the distance between the axes. For $s>\delta$ and $g\in\{a,b\}$, the
power sums $Z_g(s)=\sum_{n\in\mathbb{Z}\setminus\{0\}}e^{-s\dis(p,g^np)}$ converge, and since
$\dis(p,b^np)\le2d_0+|n|\ell_b$,
\[
Z_a(s)=\frac{2}{e^{s\ell_a}-1},\qquad Z_b(s)\ge\frac{2e^{-2sd_0}}{e^{s\ell_b}-1}.
\]
Since the words $a^{n_1}b^{m_1}\cdots a^{n_k}b^{m_k}$ with nonzero exponents are distinct for all
$k\ge1$, by the triangle inequality we have
\[
\begin{aligned}
\infty>P_\G(s,p)&\ge\sum_{k\ge1}(Z_a(s)Z_b(s))^k,\\
\frac{4e^{-2sd_0}}{(e^{s\ell_a}-1)(e^{s\ell_b}-1)} &\le Z_a(s)Z_b(s)<1.
\end{aligned}
\]
Letting $s\downarrow\delta$ and then taking the infimum over $p,q$, we get
\[
D\ge\frac1{2\delta} \log\frac4{(e^{\delta\ell_a}-1)(e^{\delta\ell_b}-1)} >T(\ell_a,\ell_b)>0,
\]
so the axes have a common perpendicular, and Lemma \ref{cuts} applies.
\end{proof}

\begin{prop}\label{scaled}
Suppose $a,b$ freely generate a convex cocompact Kleinian group $\G$. If
\[
0<\delta_\G\le\frac14,\qquad \delta_\G\max\{\ell_a,\ell_b\}\le\log2,
\]
then the pair $a,b$ admits a classical domain.
\end{prop}
\begin{proof}
We order the generators so that $x=\delta_\G\ell_a\le y=\delta_\G\ell_b\le\log2$, and denote by $D$
the axis distance. For $p$ on the first axis and $q$ on the second, by (\ref{pair-pressure}) and the
triangle inequality we have
\[
\begin{aligned}
\delta_\G\dis(p,bp)&\ge F(x):=\log\frac{e^x+3}{e^x-1},\\
\dis(p,bp)&\le2\dis(p,q)+\ell_b,\\
D&\ge\frac{F(x)-y}{2\delta_\G}.
\end{aligned}
\]
Here the last line is obtained by taking the infimum over $p,q$. With $Q=(e^x+3)/(2(e^x-1))$ and
$z=e^x\in(1,2]$,
\[
D\ge\frac{\log Q}{2\delta_\G}\ge2\log Q,\qquad \frac{Q}{\coth x}
=\frac{(z+3)(z+1)}{2(z^2+1)}\ge\frac32,
\]
where the last inequality is equivalent to $z\le2$. Since $\ell_a/4=x/(4\delta_\G)\ge x$ and
$\ell_b/4\ge x$, and since $\coth$ decreases,
\[
T(\ell_a,\ell_b)\le2\log\coth x \le D-2\log(3/2)<D,
\]
and Lemma \ref{cuts} applies.
\end{proof}
\begin{lem}\label{F}
For $y>0$, set $F(y)=\log((e^y+3)/(e^y-1))$. Then
\[
F(y)\ge\log\frac{13}{4y}\qquad(0<y\le5/4).
\]
\end{lem}
\begin{proof}
The inequality is equivalent to $(13-4y)e^y\le13+12y$, that is, to $h\ge0$ on $(0,5/4]$ for
\[
h(y)=\log\frac{13+12y}{13-4y}-y, \qquad h'(y)=\frac{48y^2-104y+39}{(13+12y)(13-4y)}.
\]
Of the roots $(13\pm2\sqrt{13})/12$ of the numerator only the smaller lies in $(0,5/4)$, which
implies $h$ first increases and then decreases on $[0,5/4]$ and attains its minimum there at an
endpoint, where $h(0)=0$ and $h(5/4)=\log(7/2)-5/4>0$, the latter since the exponential series shows
that
\[
e<\frac52+\frac16\sum_{j\ge0}4^{-j}=\frac{49}{18},
\]
so that $e^5<(49/18)^5<(7/2)^4$.
\end{proof}

\begin{lem}\label{bounded}
Suppose $a,b$ freely generate a convex cocompact Kleinian group $\G$ with
\[
0<\delta_\G\le1/3,\qquad 0<\ell_a,\ell_b\le2.
\]
Then this pair admits a classical domain.
\end{lem}
\begin{proof}
For $s,t>0$ set
\[
f_s(t)=\coth(t/4)^{2s}(e^{st}-1).
\]
For each fixed $t>0$,
\[
\frac{\partial}{\partial s}\log f_s(t) =2\log\coth(t/4)+\frac{t}{1-e^{-st}}>0,
\]
so $f_s(t)$ is strictly increasing in $s$, and for $0<s\le1/2$ it is also increasing in $t$, since
\[
\frac{\partial}{\partial t}\log f_s(t) =s\left(\frac{1}{1-e^{-st}}-\frac{1}{\sinh(t/2)}\right)>0
\]
since
\[
1-e^{-st}<st\le t/2<\sinh(t/2).
\]
Hence, for $0<t\le2$,
\[
f_{\delta_\G}(t)\le f_{1/3}(2) =\coth(1/2)^{2/3}(e^{2/3}-1)<2.
\]
Here the last inequality holds since $e>5/2$ implies $\coth(1/2)=(e+1)/(e-1)<7/3$, and
$e<49/18<11/4$ implies $e^2<8$ and hence $e^{2/3}-1<1$. Therefore
\[
e^{2\delta_\G T(\ell_a,\ell_b)} (e^{\delta_\G\ell_a}-1)(e^{\delta_\G\ell_b}-1)
=f_{\delta_\G}(\ell_a)f_{\delta_\G}(\ell_b)<4,
\]
which is the strict condition of Theorem \ref{scalar}, so the pair $a,b$ admits a classical domain.
\end{proof}
\begin{prop}\label{nielsen}
Suppose $a,b$ freely generate a convex cocompact group and both translation lengths are at least
$L_0=37/10$. If the axes are disjoint, let $p,q$ be the endpoints of their common perpendicular on
the axes of $a,b$, respectively, and if the axes intersect, let $p=q$ be their intersection point.
We cut the axis of $a$ perpendicularly at signed distances $\pm\ell(a)/2$ from $p$, and the axis of
$b$ at signed distances $\pm\ell(b)/2$ from $q$, and at each cutting point take the closed halfspace
whose intersection with the corresponding axis is the ray directed away from $p$ or $q$,
respectively. Its ideal boundary is the outward ideal ball. Then these four outward ideal balls give
a classical domain, or a Nielsen transformation decreases the sum of the translation lengths by more
than
\[
c_0=\frac{37}{10}-2\operatorname{arcosh}(61/19)>0.
\]
Repeating these replacements leads, after a number of steps bounded in terms of the initial sum of
lengths, to a classical basis or to a basis with one translation length smaller than $L_0$.
\end{prop}
\begin{proof}
Set $A_0=L_0/2=37/20$. Since the Taylor terms are positive we have
\[
\begin{aligned}
\sinh A_0 &>A_0+A_0^3/6+A_0^5/120>\frac{37}{12},\\
\cosh A_0 &>1+A_0^2/2+A_0^4/24+A_0^6/720>\frac{13}{4}>\frac{61}{19},
\end{aligned}
\]
so that $c_0>0$, and since $(37/12)^2-19/2=1/144$, also $\sinh^2 A_0>19/2$. After interchanging
$a,b$ if necessary, $A=\ell_a/2\le B=\ell_b/2$. By the endpoint-stabilizer argument following
Proposition \ref{tubes}, two loxodromic elements sharing an endpoint lie in one cyclic subgroup, so
the axes of the basis elements $a,b$ cannot share an endpoint. Take the common-perpendicular points
$p,q$, or $p=q$ at the intersection, denote by $D$ the axis distance, and, after inverting $b$ if
necessary, choose unit tangents $u,v$ pointing in the translation directions of $a,b$ with
$c=\langle u,v\rangle\ge0$. The normals of Lemma \ref{cuts} also apply when $p=q$. If
\[
\langle n_a^\sigma,n_b^\tau\rangle<-1 \qquad(\sigma,\tau\in\{+,-\}),
\]
then, as the remaining pairs already satisfy
\[
\langle n_a^+,n_a^-\rangle=-\cosh\ell_a<-1, \qquad \langle n_b^+,n_b^-\rangle=-\cosh\ell_b<-1,
\]
every pair of distinct normals has inner product less than $-1$, their sum $2\sinh A\,p+2\sinh B\,q$
is future timelike, and the cap-separation argument in Lemma \ref{cuts} makes all four closed ideal
balls disjoint. So if these balls are not pairwise disjoint, there are $\sigma,\tau\in\{+,-\}$ with
\[
\langle n_a^\sigma,n_b^\tau\rangle\ge-1,
\]
and hence
\[
\begin{aligned}
\sinh A\sinh B\cosh D&\le1+c\cosh A\cosh B,\\
\cosh D-c&\le\frac{1+\cosh(B-A)}{\sinh A\sinh B}.
\end{aligned}
\]
Parametrizing the axes by signed arclength from $p,q$ in the directions $u,v$, let $d_\pm$ be the
distances between their matching points at parameters $\pm A$, so that
\[
\cosh d_\pm=\cosh D+(\cosh D-c)\sinh^2A.
\]
Since $0\le c\le1$, $B\ge A\ge A_0=37/20$ and $\sinh B\ge\sinh A\cosh(B-A)$,
\[
\begin{aligned}
\cosh D &\le\frac{1+c\cosh A\cosh B}{\sinh A\sinh B}\\
&\le\coth A\coth B+\operatorname{csch}A\operatorname{csch}B\\
&\le1+2\operatorname{csch}^2A<\frac{23}{19},\\
(\cosh D-c)\sinh^2A &\le\frac{\sinh A}{\sinh B}\bigl(1+\cosh(B-A)\bigr)\le2,
\end{aligned}
\]
and therefore $d_\pm<\operatorname{arcosh}(61/19)<A_0$. For $y,z$ at parameter $-A$ on the first and
second axes and $w$ at parameter $A$ on the second, $ay$ has parameter $A$ and $bz$ has parameter
$\ell_b-A$, so
\[
\begin{aligned}
\ell(a^{-1}b)&\le\dis(ay,by)\\
&\le\dis(ay,w)+\dis(w,bz)+\dis(bz,by)\\
&=d_++(\ell_b-\ell_a)+d_-\\
&<\ell_b-\ell_a+2\operatorname{arcosh}(61/19).
\end{aligned}
\]
Only translations along the axes and the invariance of distance under isometries were used, so
arbitrary rotations do not affect these estimates. For the Nielsen replacement $b\mapsto a^{-1}b$,
\[
(\ell_a+\ell_b)-(\ell_a+\ell(a^{-1}b)) >\ell_a-2\operatorname{arcosh}(61/19)\ge c_0.
\]
Since inversion preserves lengths, we repeat this step while both lengths are at least $L_0$ and the
four outward ideal balls are not pairwise disjoint. The sum of lengths remains positive, and each
step decreases it by more than $c_0$, hence the number of steps is bounded by the initial sum
divided by $c_0$, and the process terminates.
\end{proof}

For the short-generator case, following Hou's use of two-sided Nielsen changes $a^kba^l$ to control
circle centers \cite{Hou2020} (Section 4 and Lemma 11.1), we place the poles of the second generator
in a radial annulus. For a three-dimensional M\"obius transformation represented by the complex
matrix $\left(\begin{smallmatrix}a_b&b_b\\c_b&d_b\end{smallmatrix}\right)$ with $c_b\ne0$, we
identify the poles $u=b(\infty)$ and $v=b^{-1}(\infty)$ with Hou's centers $\zeta_b=a_b/c_b$ and
$\eta_b=-d_b/c_b$ \cite{Hou2020} (Section 2). We call the points
\[
u=b(\infty),\qquad v=b^{-1}(\infty)
\]
the poles of $b$ in these coordinates. These poles need not be the fixed points of $b$, see
\cite{Hou2020} (Section 4 and Lemma 11.1). The complex matrix estimates are replaced below by a
computation allowing any orthogonal part.

\begin{lem}\label{pole}
We normalize a loxodromic generator $a$ of length $L>0$ by $a(x,t)=(e^L Rx,e^Lt)$ in the upper
half-space model, where $x\in\mathbb{R}^{N-1}$, $t>0$, and $R\in O(N-1)$. Suppose $a,b$ are distinct
elements of a free basis of a convex cocompact Kleinian group $\G$. Then $b$ has the form
\begin{equation}\label{pole-form}
b(x,t)=\left( u+\frac{KQ(x-v)}{|x-v|^2+t^2}, \frac{Kt}{|x-v|^2+t^2}\right),
\end{equation}
where $Q\in O(N-1)$, $K>0$, and the poles $u=b(\infty)$ and $v=b^{-1}(\infty)$ are nonzero vectors
in $\mathbb{R}^{N-1}$. Set $\nu=\sqrt{K/(|u||v|)}>0$. For $k,l\in\mathbb{Z}$, the basis change
$b'=a^kba^l$, which fixes the other basis generators, has
\[
u'=e^{kL}R^ku,\qquad v'=e^{-lL}R^{-l}v,\qquad K'=e^{(k-l)L}K,
\]
so $\nu$ is unchanged, and the integers $k,l$ and a number $r_0>0$ can be chosen such that
\begin{equation}\label{annulus-poles}
r_0e^{L/4}\le|u'|,|v'|\le r_0e^{3L/4}.
\end{equation}
For this choice of the second generator, at $p=(0,\sqrt{K'})$,
\begin{equation}\label{pole-displacement}
\cosh\dis(p,b'p) \le1+\frac{1+\cosh L}{\nu^2}+\frac{1}{2\nu^4}.
\end{equation}
\end{lem}
\begin{proof}
By the endpoint-stabilizer argument following Proposition \ref{tubes}, an element fixing either
endpoint of $a$ fixes both, and their stabilizer is cyclic, hence $b$ fixes neither endpoint, since
$a,b$ freely generate a nonabelian subgroup.

If $b(\infty)=0$, the element $bab^{-1}$ fixes $0$, hence also $\infty$, and as its fixed-point set
is $\{b(0),b(\infty)\}$, $b(0)=\infty$ and $b$ interchanges the endpoints of the axis of $a$. Its
restriction to this axis reverses arclength and fixes one point, which implies that $b$ is elliptic,
and a nonidentity elliptic element cannot belong to the discrete torsion-free group $\G$. The same
argument applied to $b^{-1}$ shows $u,v\in\mathbb{R}^{N-1}\setminus\{0\}$. The hyperbolic isometry
\[
h_v(x,t)= \left(\frac{x-v}{|x-v|^2+t^2}, \frac{t}{|x-v|^2+t^2}\right)
\]
satisfies $h_v(v)=\infty$ and $h_v(\infty)=0$, which implies that $b\circ h_v^{-1}$ fixes infinity
and is therefore a Euclidean similarity preserving the upper half-space,
\[
(b\circ h_v^{-1})(x,t)=(KQx+u,Kt), \qquad K>0,\quad Q\in O(N-1),
\]
whose translation vector is $u$ since $(b\circ h_v^{-1})(0)=b(\infty)=u$, and composing with $h_v$
implies (\ref{pole-form}).

Let $Q'$ be the orthogonal factor of $b'$. Direct composition gives the stated $u',v',K'$ and
\[
Q'=R^kQR^l,\qquad \frac{K'}{|u'||v'|} =\frac{e^{(k-l)L}K}{e^{kL}|u|e^{-lL}|v|}=\nu^2.
\]
On the circle $\mathbb{R}/L\mathbb{Z}$, the two residues of $\log|u|$ and $\log|v|$ leave a
complementary arc of length at least $L/2$. If $s_0\in\mathbb{R}$ is a lift of the midpoint of such
an arc, both residues have circular distance at least $L/4$ from $s_0+L\mathbb{Z}$, so each residue
has a representative in $[s_0+L/4,s_0+3L/4]$. We choose $k,l\in\mathbb{Z}$ independently such that
\[
\log|u|+kL,\quad \log|v|-lL \in[s_0+L/4,s_0+3L/4],
\]
and set $r_0=e^{s_0}$. The formulas for $u',v'$ give (\ref{annulus-poles}). The inverse substitution
$b'\mapsto a^{-k}b'a^{-l}$ fixes every other basis generator, so the substitution preserves the free
basis and the group.

For the displacement, drop the primes and use the coordinates $(x/\sqrt K,t/\sqrt K)$. With
$U=u/\sqrt K$, $V=v/\sqrt K$, $\alpha=|U|^2$, $\beta=|V|^2$ and $\zeta=U\cdot QV$,
\[
\begin{aligned}
p&=(0,1),\qquad bp=\left(U-\frac{QV}{1+\beta},\frac1{1+\beta}\right),\\
\cosh\dis(p,bp)&=1+\frac12(\alpha+\beta+\alpha\beta-2\zeta).
\end{aligned}
\]
With $r=|u|/|v|$, it follows from (\ref{annulus-poles}) that
\[
\begin{gathered}
\alpha=r/\nu^2,\qquad \beta=1/(r\nu^2),\qquad \alpha\beta=\nu^{-4},\qquad |\zeta|\le\nu^{-2},\\
e^{-L}\le e^{-L/2}\le r\le e^{L/2}\le e^L,\qquad r+r^{-1}\le2\cosh L,
\end{gathered}
\]
and substituting these bounds implies (\ref{pole-displacement}).
\end{proof}
\begin{lem}\label{axis-minimum}
Let $\G<\Isom(\mathbb{H}^N)$ be discrete, free, and convex cocompact, and let $a,b$ be distinct
elements of a free basis of $\G$. In upper half-space, set
\[
a(x,t)=(e^LRx,e^Lt), \qquad L=\ell(a)>0,\qquad R\in O(N-1),
\]
and
\[
b(x,t)= \left( u+\frac{KQ(x-v)}{|x-v|^2+t^2}, \frac{Kt}{|x-v|^2+t^2} \right),
\]
where $K>0$, $Q\in O(N-1)$, and
\[
u=b(\infty),\qquad v=b^{-1}(\infty),\qquad u,v\in\mathbb{R}^{N-1}\setminus\{0\}.
\]
Set
\[
A_a=\{(0,t):t>0\},\qquad \nu=\sqrt{\frac{K}{|u||v|}},
\]
and
\[
\Delta=|u|^2|v|^2+K^2-2K\,u\cdot Qv.
\]
Then $\Delta>0$, and the displacement by $b$ restricted to $A_a$ attains its minimum at the unique
point $(0,\Delta^{1/4})$. Define
\[
D_a(b)=\min_{p\in A_a}\dis(p,bp).
\]
Then we have
\[
\cosh D_a(b) =\frac{|u|^2+|v|^2+2\sqrt{\Delta}}{2K}.
\]
If
\[
e^{-L/2}\le\frac{|u|}{|v|}\le e^{L/2},
\]
then
\begin{equation}\label{pole-upper}
\cosh D_a(b)\le 1+\frac{1+\cosh(L/2)}{\nu^2}.
\end{equation}
\end{lem}
\begin{proof}
The boundary extension satisfies
\[
b(0)=u-\frac{KQv}{|v|^2}.
\]
The endpoint-stabilizer argument following Proposition \ref{tubes} implies $b(0)\ne0$. Hence
\[
\Delta=|v|^2|b(0)|^2>0.
\]
By the upper half-space distance formula we have
\[
\cosh\dis\bigl((0,t),b(0,t)\bigr) =\frac{|u|^2+|v|^2+t^2+\Delta/t^2}{2K}.
\]
Since
\[
t^2+\frac{\Delta}{t^2}\ge2\sqrt{\Delta},
\]
with equality if and only if $t=\Delta^{1/4}$, the minimum formula follows. Furthermore,
\[
|u\cdot Qv|\le|u||v|, \qquad \sqrt{\Delta}\le|u||v|+K.
\]
Therefore we have
\[
\begin{aligned}
\cosh D_a(b) &\le1+\frac{|u|^2+|v|^2+2|u||v|}{2K}\\
&=1+\frac{ 1+\frac12\left(\frac{|u|}{|v|} +\frac{|v|}{|u|}\right)}{\nu^2}.
\end{aligned}
\]
The ratio assumption bounds the expression in parentheses by $2\cosh(L/2)$, proving
(\ref{pole-upper}).
\end{proof}
\begin{lem}\label{third-pole}
Let $\G$ be a convex cocompact Kleinian group, and let $a,b$ belong to a free basis of $\G$. Suppose
\[
0<\delta_\G\le\frac13,\qquad 0<L=\ell(a)<L_0=\frac{37}{10}.
\]
We choose upper half-space coordinates such that
\[
a(x,t)=(e^LRx,e^Lt),\qquad R\in O(N-1),
\]
and let $u,v,K,Q$ be the pole parameters of $b$ in (\ref{pole-form}). Assume
\[
e^{-L/2}\le\frac{|u|}{|v|}\le e^{L/2}, \qquad \nu=\sqrt{\frac{K}{|u||v|}}.
\]
Then we have
\[
\nu<\tanh(L/4).
\]
\end{lem}
\begin{proof}
Since $\dis(p,ap)=L$ at every point $p$ of the axis of $a$, it follows from (\ref{pair-pressure})
that
\[
\dis(p,bp)\ge\frac{F(\delta_\G L)}{\delta_\G}, \qquad F(y)=\log\frac{e^y+3}{e^y-1},
\]
in particular at the minimizing point in Lemma \ref{axis-minimum}, and since $\delta_\G
L<37/30<5/4$, Lemma \ref{F} gives
\begin{equation}\label{D-lower}
D_a(b)\ge \frac{1}{\delta_\G}\log\frac{13}{4\delta_\G L} \ge3\log\frac{39}{4L},
\end{equation}
where the logarithms are positive since $L<37/10$ and $\delta_\G\le1/3$. Set $x=L/4$ and
$x_0=37/40<1$. If $\nu\ge\tanh x$, then by (\ref{pole-upper}) we have
\[
\cosh D_a(b)\le1+\frac{1+\cosh(2x)}{\tanh^2x} =1+\frac{2\cosh^4x}{\sinh^2x},
\]
and as the function
\[
G(x)=x^3\left(1+\frac{2\cosh^4x}{\sinh^2x}\right)
\]
is increasing on $(0,1]$,
\[
\begin{aligned}
x^3\cosh D_a(b) &\le G(x)<G(x_0)\\
&<\left(\frac{37}{40}\right)^3 \left(1+\frac{2(117/80)^4}{(17/16)^2}\right),
\end{aligned}
\]
and by (\ref{D-lower}) and $x=L/4$ we have
\[
x^3\cosh D_a(b) \ge\frac{x^3}{2}\left(\frac{39}{4L}\right)^3 =\frac12\left(\frac{39}{16}\right)^3,
\]
which contradicts the preceding upper bound, hence $\nu<\tanh(L/4)$.
\end{proof}
\begin{rem}
Claude was used to help with computations of Lemma \ref{third-pole}.
\end{rem}
Following Hou's estimates for the centers of circles \cite{Hou2010}, we next give pole estimates for
two generators.
\begin{lem}\label{pole-estimate}
Let $\G$ be a convex cocompact Kleinian group, and suppose that $a,b$ belong to a free basis of
$\G$. Assume
\[
0<\delta_\G\le\frac18,\qquad 0<L=\ell(a)<4.
\]
We normalize $a(x,t)=(e^LRx,e^Lt)$, and denote by $u,v,K$ the pole parameters of $b$ in
(\ref{pole-form}). If
\[
e^{-L}\le\frac{|u|}{|v|}\le e^L, \qquad \nu=\sqrt{\frac{K}{|u||v|}},
\]
then $\nu<L/8$.
\end{lem}
\begin{proof}
Set $p=(0,\sqrt K)$, $D=\dis(p,bp)$. By the computation in Lemma \ref{pole}, which uses only the
stated bound on $|u|/|v|$ we have
\[
\cosh D \le1+\frac{1+\cosh L}{\nu^2} +\frac{1}{2\nu^4}.
\]
Since $p$ lies on the axis of $a$, $\dis(p,ap)=L$, and by (\ref{pair-pressure}) applied to $a,b$
with the exponent $\delta_\G$ we have
\[
D\ge\frac{1}{\delta_\G} \log\frac{e^{\delta_\G L}+3}{e^{\delta_\G L}-1}.
\]
Since $\delta_\G L<1/2$ and $e^{\delta_\G L}-1\le2\delta_\G L$, and since $\delta_\G^{-1}\ge8$ and
$\log(2/(\delta_\G L))\ge\log(16/L)>0$,
\begin{equation}\label{pole-lower}
D\ge\frac{1}{\delta_\G}\log\frac{2}{\delta_\G L} \ge8\log\frac{16}{L}.
\end{equation}
If $\nu\ge L/8$, then, since $L<4$ and $\cosh4<28$,
\[
\cosh D \le1+\frac{1856}{L^2}+\frac{2048}{L^4}<\frac12\left(\frac{16}{L}\right)^8,
\]
and by (\ref{pole-lower}) we have
\[
\cosh D\ge\frac12\left(\frac{16}{L}\right)^8,
\]
a contradiction, hence $\nu<L/8$.
\end{proof}
The pole bound is combined next with Hou's annular placement and changes of generators
\cite{Hou2020} (Section 4 and Proposition 5.6).
\begin{lem}\label{annulus-fit}
In the pole parameters of Lemma \ref{pole}, set $U=|u|$ and $V=|v|$, and suppose that $e^{-L/2}\le
U/V\le e^{L/2}$ and $\nu<\tanh(L/4)$. Then the balls centered at $u,v$ with radii $\nu U,\nu V$ lie
strictly inside one annulus $R_0<|x|<e^LR_0$ and are paired by $b$.
\end{lem}
\begin{proof}
With $m=\min\{U,V\}$ and $M=\max\{U,V\}$ we have $M/m\le e^{L/2}$ and
\[
\frac{1+\nu}{1-\nu}<e^{L/2},
\]
Hence $M(1+\nu)e^{-L}<m(1-\nu)$, which implies that $R_0$ can be chosen with
\[
M(1+\nu)e^{-L}<R_0<m(1-\nu),
\]
and then every point of either closed ball has norm strictly between $R_0$ and $e^LR_0$. The formula
$|b(x)-u|=K/|x-v|$ and $(\nu U)(\nu V)=K$ show that $b$ maps the closed exterior of the ball at $v$
onto the ball at $u$, and the two concentric balls
\[
\{|x|\le R_0\},\qquad \{|x|\ge e^LR_0\}\cup\{\infty\}
\]
are paired by $a$.
\end{proof}

\begin{lem}\label{overlap-lem}
If the two balls in Lemma \ref{annulus-fit} intersect, including tangency, then
\[
\ell(b)\le g(L):=\operatorname{arcosh}(2+\cosh(L/2)).
\]
The function $g$ is increasing and $g(2)<2$, and with $\varepsilon_0=2-g(2)>0$ we have $g(L)\le
L-\varepsilon_0$ for $L\ge2$.
\end{lem}
\begin{proof}
Intersection implies $|u-v|\le\nu(U+V)$, and since $b(v,\sqrt K)=(u,\sqrt K)$, by the distance
formula we have
\[
\begin{aligned}
\ell(b) &\le\operatorname{arcosh}\left(1+\frac{|u-v|^2}{2K}\right)\\
&\le\operatorname{arcosh}\left(2+\frac12(U/V+V/U)\right) \le g(L).
\end{aligned}
\]
Note that we have
\[
0<g'(L) =\frac{\sinh(L/2)}{2\sqrt{(2+\cosh(L/2))^2-1}}<\frac12,
\]
so $g$ and $L-g(L)$ are increasing, and $g(2)<2$ since $\cosh2-(2+\cosh1)>0$, as $\cosh1>3/2$, hence
$L-g(L)\ge2-g(2)=\varepsilon_0$ when $L\ge2$.
\end{proof}

\begin{thm}\label{short-annulus}
Suppose $a,b$ freely generate a convex cocompact Kleinian group $\G$ with $0<\delta_\G\le1/3$ and
$0<\ell(a)<L_0=37/10$. Then a terminating sequence of two-sided changes $b\mapsto a^kba^l$ and
exchanges of generators gives a free basis admitting a classical domain.
\end{thm}
\begin{proof}
Set $(a_1,b_1)=(a,b)$. Successive free bases $(a_j,b_j)$ of $\G$ with
\[
L_j=\ell(a_j),\qquad 0<L_j<L_0,\qquad j\ge1,
\]
are defined as follows, with the function $g$ and the constant $\varepsilon_0>0$ of Lemma
\ref{overlap-lem}. At stage $j$, choose upper half-space coordinates such that
\[
a_j(x,t)=(e^{L_j}R_jx,e^{L_j}t), \qquad R_j\in O(N-1),
\]
and, by Lemma \ref{pole} applied to $(a_j,b_j)$, choose $k_j,l_j\in\mathbb{Z}$ such that
\[
b_j'=a_j^{k_j}b_ja_j^{l_j}
\]
has pole parameters $u_j',v_j',K_j',Q_j'$ satisfying (\ref{annulus-poles}), in particular
\[
e^{-L_j/2}\le\frac{|u_j'|}{|v_j'|} \le e^{L_j/2}.
\]
With
\[
\nu_j=\sqrt{\frac{K_j'}{|u_j'||v_j'|}},
\]
By Lemma \ref{third-pole}, applied to the free basis $(a_j,b_j')$ of $\G$, we have
\[
\nu_j<\tanh(L_j/4).
\]
The closed Euclidean balls $B_j^-$ and $B_j^+$ centered at $v_j'$ and $u_j'$, with respective radii
$\nu_j|v_j'|$ and $\nu_j|u_j'|$, satisfy
\[
b_j'\bigl(\mathbb{S}^d\setminus\inte B_j^-\bigr)=B_j^+,
\]
and by Lemma \ref{annulus-fit} they lie strictly between two concentric spheres paired by $a_j$. If
$B_j^-\cap B_j^+=\emptyset$, the four paired balls give a classical domain for $(a_j,b_j')$. If
$B_j^-\cap B_j^+\ne\emptyset$, then by Lemma \ref{overlap-lem}, applied to $(a_j,b_j')$, we have
\[
\ell(b_j')\le g(L_j),
\]
and in the case $L_j<2$
\[
\ell(a_j)<2,\qquad \ell(b_j')\le g(L_j)<g(2)<2,
\]
so that Lemma \ref{bounded} gives a classical domain for $(a_j,b_j')$, and in the case $2\le
L_j<L_0$ the next basis
\[
(a_{j+1},b_{j+1})=(b_j',a_j)
\]
satisfies
\[
0<L_{j+1} =\ell(b_j') \le g(L_j) \le L_j-\varepsilon_0<L_0.
\]
Each exchange thus decreases the first translation length by at least $\varepsilon_0$, so since
$L_1<L_0$, after at most
\[
\left\lceil\frac{L_0-2}{\varepsilon_0}\right\rceil+1
\]
exchanges either disjoint balls have been obtained or the first translation length is below two, and
in both cases a classical domain has been found.
\end{proof}

\begin{proof}[Proof of Theorem \ref{main-2}]
Let $a,b$ be a generating pair of $\G$ and $o\in\mathbb{H}^N$. By the structural theorem, $\G$ is
free and convex cocompact, with $\delta_\G\le1/3$, and being non-elementary and generated by two
elements, it has free rank two. If $F_2$ is the free group with basis $\alpha,\beta$, the surjection
\[
F_2\longrightarrow\G\cong F_2,\qquad \alpha\longmapsto a,\quad\beta\longmapsto b
\]
is an isomorphism by Hopficity, which implies that $a,b$ is a free basis, and every nonidentity
element is loxodromic. With $M=\max\{\dis(o,ao),\dis(o,bo)\}>0$, each of the $4\cdot3^{k-1}$ reduced
words of length $k$ has displacement at most $kM$, so
\[
\begin{aligned}
P_\G(s,o)&\ge\sum_{k\ge1}4\cdot3^{k-1}e^{-skM} =\infty\quad(sM\le\log3),\\
\delta_\G&\ge\frac{\log3}{M}>0.
\end{aligned}
\]
If one generator has length below $L_0=37/10$, exchange the generators if necessary and apply
Theorem \ref{short-annulus}. Otherwise Proposition \ref{nielsen} leads to a classical basis or to a
basis with a generator of length below $L_0=37/10$, and in the latter case Theorem
\ref{short-annulus} applies.
\end{proof}

\begin{cor}\label{2-subgroups}
Suppose $\G<\Isom(\mathbb{H}^N)$ is discrete and torsion-free and satisfies $\delta_\G\le1/3$. Then
every non-elementary subgroup generated by two elements is classical. The conclusion also holds
under $\mathfrak{D}_\G\le1/3$.
\end{cor}
\begin{proof}
Such a subgroup $H$ is discrete and torsion-free, its Poincar\'e series is a subseries of that of
$\G$, and $\Lambda_H\subset\Lambda_\G$, so
\[
\delta_H\le\delta_\G,\qquad \dim_H\Lambda_H\le\dim_H\Lambda_\G,
\]
and the corresponding version of Theorem \ref{main-2} applies. If $\G$ is elementary, there is no
such subgroup $H$.
\end{proof}


\section{Carrier Reductions and Standard Form}
In this section the pole estimate is used for simultaneous annular placement, and we introduce the
carrier graphs that will be used in the proof of Theorem \ref{main}. Following the annular
normalization in Hou's standard form and generator-selection argument \cite{Hou2020} (Definition 6.2
and Lemma 11.1), one annulus is obtained for every pole pairs in higher dimensions. Mutual
disjointness of the other pairs needs a separate argument, just as Hou's Lemma 5.2 imposes separate
conditions for annular position and mutual separation.

\begin{thm}\label{annulus}
Let $m\ge1$ and let $\G$ be a convex cocompact Kleinian group freely generated by $a,b_1,...,b_m$,
with $0<\delta_\G\le1/8$ and $0<L=\ell(a)<4$. We normalize $a(x,t)=(e^LRx,e^Lt)$, where $R\in
O(N-1)$, let $u_i=b_i(\infty),v_i=b_i^{-1}(\infty)\in\mathbb{R}^{N-1}\setminus\{0\}$ and $K_i>0$ be
the pole parameters of $b_i$ given by Lemma \ref{pole} as in (\ref{pole-form}), and set
$\nu_i=\sqrt{K_i/(|u_i||v_i|)}>0$. Then
\begin{equation}\label{sum-nu}
\sum_{i=1}^m\nu_i<\frac{L}{8}.
\end{equation}
Note that there are two-sided Nielsen changes $b_i'=a^{k_i}b_i a^{l_i}$, with
$k_i,l_i\in\mathbb{Z}$, and one $r_0>0$ such that, with $u_i',v_i',K_i'$ the pole parameters of
$b_i'$, all $2m$ paired Euclidean balls centered at $v_i',u_i'$, with respective radii $\nu_i|v_i'|$
and $\nu_i|u_i'|$, lie strictly inside the common annulus $r_0<|x|<r_0e^L$, called an $a$-annulus.
\end{thm}
\begin{proof}
We fix $R_0>0$ and apply independent changes $b_i\mapsto a^{k_i}b_i a^{l_i}$, which preserve
$\nu_i$, hence both pole norms lie in $[R_0,e^LR_0]$, and keep writing $b_i,u_i,v_i,K_i,Q_i$ for the
changed quantities, hence $e^{-L}\le|u_i|/|v_i|\le e^L$. Since $a,b_i$ freely generate a subgroup of
$\G$, (\ref{pair-pressure}) in Lemma \ref{pressure} gives, for all $p\in\mathbb{H}^N$,
\[
\bigl(e^{\delta_\G\dis(p,ap)}-1\bigr) \bigl(e^{\delta_\G\dis(p,b_ip)}-1\bigr)\ge4,
\]
and Lemma \ref{pole-estimate} applies to $a,b_i$, since
\[
e^{-L}\le\frac{|u_i|}{|v_i|}\le e^L,
\]
giving $\nu_i<L/8<1/2$ for all $i$.

At the common point $p=(0,R_0)$, set $D_i=\dis(p,b_i p)$. By the upper half-space distance formula,
\[
\cosh D_i= \frac{(|u_i|^2+R_0^2)(|v_i|^2+R_0^2)+K_i^2 -2K_i u_i\cdot Q_i v_i}{2K_iR_0^2},
\]
and since $U_i=|u_i|/R_0$ and $V_i=|v_i|/R_0$ lie in $[1,e^L]$, the right side is at most
\[
\frac{2\cosh^2L}{\nu_i^2} +\frac{\nu_i^2e^{2L}}2+e^{2L} <\frac{8192}{\nu_i^2},
\]
since $\cosh4<28$, $e^4<55$, $\nu_i<1$ and $2\cdot28^2+(3/2)\cdot55^2<8192$. It follows that
$D_i<14\log2-2\log\nu_i$.

For $s>\delta_\G$, consider the words
\[
a^{n_1}b_{i_1}\cdots a^{n_k}b_{i_k}, \qquad n_j\in\mathbb{Z},\quad 1\le i_j\le m.
\]
Since the $b_i$ letters are positive, they do not cancel, and reading each of them together with its
preceding power of $a$ gives a unique decomposition, also when $n_j=0$, which distinguishes words
with different values of $k$. With
\[
\rho(s)=\left(\sum_{n\in\mathbb{Z}}e^{-s|n|L}\right) \left(\sum_i e^{-sD_i}\right)
=\coth(sL/2)\sum_i e^{-sD_i},
\]
the triangle inequality gives the geometric lower bound
\[
P_\G(s,p)\ge\sum_{k\ge1}\rho(s)^k.
\]
Since the Poincar\'e series converges, $\rho(s)<1$, and letting $s\downarrow\delta_\G$,
\[
\sum_i\nu_i^{2\delta_\G} \le2^{14\delta_\G}\tanh(\delta_\G L/2) <2^{14\delta_\G}\delta_\G L/2.
\]
With $q=2\delta_\G\le1$ and $(\sum_i\nu_i)^q\le\sum_i\nu_i^q$ we get
\[
\sum_i\nu_i< 128(\delta_\G L/2)^{1/(2\delta_\G)} \le128(L/16)^4 =\frac{L^4}{512}<\frac{L}{8},
\]
since $\delta_\G L/2\le L/16<1$ and $1/(2\delta_\G)\ge4$, which is (\ref{sum-nu}).

On the circle $\mathbb{R}/L\mathbb{Z}$, exclude all points whose distance from the image of
$\log|u_i|$ or $\log|v_i|$ is at most $2\nu_i$ for some $i$. The resulting closed arcs have total
length at most $8\sum_i\nu_i<L$, so there is $s_0\in\mathbb{R}$ whose image in
$\mathbb{R}/L\mathbb{Z}$ lies outside them. Set $r_0=e^{s_0}$. After independent integer shifts,
denote by $b_i'$ the final generators and by $u_i',v_i',K_i'$ their pole parameters, hence each of
$|u_i'|,|v_i'|$ lies strictly between $r_0e^{2\nu_i}$ and $r_0e^{L-2\nu_i}$. Since $0<\nu_i<1/2$,
$1-e^{-2\nu_i}>\nu_i$ and $e^{2\nu_i}-1>\nu_i$, so every ball whose radius is $\nu_i$ times its
center norm lies strictly in the annulus, and the radius product $\nu_i^2|u_i'||v_i'|=K_i'$ shows
that each pairing is preserved. Note that every basis change fixes $a$ and every other $b_j$, so
these changes commute.
\end{proof}

\begin{rem}
The preceding argument places all balls between the two concentric spheres paired by $a$. Their
mutual disjointness requires a separate argument. In rank two, Theorem \ref{main-2} gives a
classical basis. Under $\mathfrak{D}_\G\le\lambda_{m+1}$, the global carrier method of Section $7$
gives a classical basis independently of these balls.
\end{rem}

For $x,y>0$, define $T(x,y)=\log\coth(x/4)+\log\coth(y/4)$. Suppose $r\ge2$ and $g_1,...,g_r$ are
loxodromic transformations whose axes $A_i$ are pairwise disjoint and no two share an ideal
endpoint. Each pair of axes then has a common perpendicular of positive length. Write
$\ell_i=\ell(g_i)$ and $D_{ij}=\dis(A_i,A_j)$, let $p_{ij}\in A_i$ be the foot of this perpendicular
for $j\ne i$, and let $W_i$ be the length of the smallest closed segment of $A_i$ containing all
$p_{ij}$, hence $W_i=0$ if these points coincide.

\begin{thm}\label{footspan}
Suppose that
\[
\ell_i>W_i,\qquad D_{ij}>T(\ell_i-W_i,\ell_j-W_j)\quad(i\ne j).
\]
Then $g_1,...,g_r$ admit a classical round-ball domain.
\end{thm}
\begin{proof}
We orient each $A_i$ in the translation direction of $g_i$, let $p_i$ be the midpoint of the segment
defining $W_i$, and cut $A_i$ perpendicularly at signed distances $\pm\ell_i/2$ from $p_i$. The two
outward ideal balls are disjoint, and they are paired by $g_i$, since it translates along $A_i$ by
$\ell_i$ and sends perpendicular hyperplanes to perpendicular hyperplanes under its rotation.

For $i\ne j$, let $p_i$ have signed coordinate $t$ from $p_{ij}$. Since $|t|\le W_i/2$,
\[
t+\ell_i/2\ge(\ell_i-W_i)/2,\qquad t-\ell_i/2\le-(\ell_i-W_i)/2,
\]
so by nesting along $A_i$ the actual outward halfspaces lie inside those cut at $\pm(\ell_i-W_i)/2$
around $p_{ij}$. Set
\[
A=\frac{\ell_i-W_i}{2},\qquad B=\frac{\ell_j-W_j}{2},\qquad p=p_{ij},\qquad q=p_{ji}.
\]
Let $u,v$ be the unit tangent vectors to $A_i,A_j$ at $p,q$, respectively, pointing in the
translation directions of $g_i,g_j$. The hyperplanes perpendicular to $A_i$ at signed distances $\pm
A$ from $p$, and to $A_j$ at signed distances $\pm B$ from $q$, have outward unit spacelike normals
\[
n_i^\sigma=\sinh A\,p+\sigma\cosh A\,u, \qquad n_j^\tau=\sinh B\,q+\tau\cosh B\,v, \qquad
\sigma,\tau\in\{-1,1\},
\]
whose positive halfspaces do not contain $p$ and $q$, respectively, and by Lemma \ref{cuts},
\[
\langle n_i^\sigma,n_j^\tau\rangle \le-\sinh A\sinh B\cosh D_{ij}+\cosh A\cosh B<-1
\qquad(\sigma,\tau\in\{-1,1\}),
\]
since $D_{ij}>T(2A,2B)$ and
\[
\cosh T(2A,2B) =\frac{1+\cosh A\cosh B}{\sinh A\sinh B}.
\]
For each choice of signs $\sigma,\tau$, the sum $n_i^\sigma+n_j^\tau$ is timelike, with
\[
\langle p,n_i^\sigma+n_j^\tau\rangle =-\sinh A-\sinh B\cosh D_{ij}<0,
\]
since $\langle p,v\rangle=0$, so the sum is future timelike. A future null vector in both positive
ideal balls would have nonnegative inner product with this sum, which is impossible, since the inner
product of a future null vector with a future timelike vector is negative. It follows that the
closed ideal balls determined by $n_i^\sigma,n_j^\tau$ are disjoint, and so are the original outward
balls, which lie in them. Applied to every pair, this gives all $2r$ disjoint balls with the
original pairings.
\end{proof}

\begin{cor}\label{footspan-cor}
Suppose $g_1,...,g_r$ form a free basis of a convex cocompact Kleinian group $\G$. Assume
\[
0<\delta_\G\le\frac14,\qquad \delta_\G\max_i\ell_i\le\log2.
\]
Then any two axes $A_i,A_j$, $i\ne j$, are disjoint and have no common endpoint on
$\partial_\infty\mathbb{H}^N$, so they have a unique common perpendicular of positive length. If, in
addition, this basis satisfies
\[
W_i\le(1-4\delta_\G)\ell_i\quad(1\le i\le r),
\]
then the whole group is classical. In particular, in every rank it suffices that $\delta_\G\le1/8$,
$\delta_\G\max_i\ell_i\le\log2$ and $W_i\le\ell_i/2$.
\end{cor}
\begin{proof}
For each pair, set $x=\delta_\G\min\{\ell_i,\ell_j\}$ and
$y=\delta_\G\max\{\ell_i,\ell_j\}\le\log2$, and let $F(t)=\log((e^t+3)/(e^t-1))$ for $t>0$. Assume
$\ell_i\le\ell_j$ and choose $p\in A_i$. Since $\dis(p,g_ip)=\ell_i$ and $x=\delta_\G\ell_i$, it
follows from (\ref{pair-pressure}) that
\[
(e^x-1) \bigl(e^{\delta_\G\dis(p,g_jp)}-1\bigr)\ge4.
\]
Therefore we have
\[
\dis(p,g_jp) \ge\frac1{\delta_\G} \log\left(1+\frac4{e^x-1}\right) =\frac{F(x)}{\delta_\G}.
\]
For all $p\in A_i$ and $q\in A_j$, we have
\[
\frac{F(x)}{\delta_\G} \le\dis(p,g_jp) \le2\dis(p,q)+\ell_j =2\dis(p,q)+\frac{y}{\delta_\G}.
\]
Taking the infimum over both variables, we get
\[
D_{ij} =\inf_{\substack{p\in A_i\\q\in A_j}}\dis(p,q) \ge\frac{F(x)-y}{2\delta_\G}>0,
\]
since $F(x)\ge\log5>y$. The estimate also applies before intersecting or asymptotic axes have been
excluded, hence the axes are disjoint, have no common ideal endpoint, and their common perpendicular
is attained. With
\[
Q(x)=\frac{e^x+3}{2(e^x-1)},
\]
the bounds $F(x)-y\ge\log Q(x)$ and $\delta_\G\le1/4$ give $D_{ij}\ge2\log Q(x)$, and for
$z=e^x\in(1,2]$
\[
\frac{Q(x)}{\coth x} =\frac{(z+3)(z+1)}{2(z^2+1)}>1.
\]
Since $\ell_i-W_i\ge4\delta_\G\ell_i$ and $\coth$ decreases,
\[
T(\ell_i-W_i,\ell_j-W_j) \le2\log\coth x <2\log Q(x)\le D_{ij},
\]
and Theorem \ref{footspan} applies.
\end{proof}

Let $G$ be a graph all of whose vertices have valence $<\infty$, the valence of a vertex being the
number of incident edge ends, counting each loop twice. Denote by $V(G)$ the vertex set and by
$E(G)$ the unoriented edge set. Represent each edge $e$ by an interval $[0,L_e]$, where $L_e>0$,
identify its endpoints with the corresponding vertices, and give $G$ the resulting quotient
topology. Denote an oriented edge from $v$ to $w$ by $e:v\to w$ and its reverse by $\bar e$, with
$L_{\bar e}=L_e$. On each connected component,
\[
d_G(x,y)=\inf_{P:x\to y}\leg(P),
\]
is a metric inducing the quotient topology on that component: at each vertex, sufficiently small
metric balls consist of initial segments of its incident edges, and since the valence at each vertex
is $<\infty$ and all $L_e>0$, these balls form a neighborhood basis for the quotient topology, and
at each interior edge point both topologies are the interval topology. For points $x,y$ in a tree
$K$, $[x,y]_K$ is the unique shortest path from $x$ to $y$, oriented in that direction, the constant
path when $x=y$. Paths are concatenated in their traversal order, and $P^{-1}$ is the reversed path.

For a free Kleinian group $\G<\Isom(\mathbb{H}^N)$, we call a connected graph $X$ with $V$ vertices
and $E$ edges, where $V,E$ are positive integers, a (metric) carrier graph if it has a continuous
mapping $\bar f:X\to\mathbb{H}^N/\G$ inducing an isomorphism on fundamental groups \cite{White2002}.
Let $\widetilde X$ be its universal cover, a tree, with its deck group identified with $\G$ through
this isomorphism, hence the lift $f:\widetilde X\to\mathbb{H}^N$ is $\G$-equivariant. For each edge
$e$, choose lifted endpoints $\widetilde v,\widetilde w$ and set
\[
L_e=\dis(f(\widetilde v),f(\widetilde w))>0,
\]
which by equivariance does not depend on the chosen lift. We give $X$ the corresponding path metric
and require $f$ to map each lifted edge isometrically onto $[f(\widetilde v),f(\widetilde w)]$.
\begin{thm}\label{carrier}
Let $\G<\Isom(\mathbb{H}^N)$ be a free Kleinian group, let $(X,\bar f)$ be a carrier graph of rank
$r\ge2$ with $|E(X)|<\infty$, and let $f:\widetilde X\to\mathbb{H}^N$ be its $\G$-equivariant lift,
mapping each lifted edge isometrically onto a geodesic segment of length $L_e>0$. Suppose that, for
all $\widetilde v\in V(\widetilde X)$ and any two distinct oriented edges $\widetilde e_1,\widetilde
e_2$ starting there,
\[
\angle_{f(\widetilde v)} \bigl(f(\widetilde e_1),f(\widetilde e_2)\bigr) \ge\theta_0, \qquad
0<\theta_0\le\pi.
\]
If
\[
L_e> 2\operatorname{arcosh}\bigl(\csc(\theta_0/2)\bigr) \qquad(e\in E(X)),
\]
then $\G$ is classical round-ball Schottky. For $\theta_0=2\pi/3$, the length bound equals
\[
2\operatorname{arcosh}(2/\sqrt3)=\log3.
\]
\end{thm}
\begin{proof}
For an oriented lifted edge $v\to w$, set $\ell(v,w)=\dis(f(v),f(w))$ and let $E(v,w)$ be the closed
halfspace beyond the perpendicular midpoint plane of $[f(v),f(w)]$ toward $f(w)$. By the length
assumption we have
\[
\cosh\bigl(\ell(v,w)/2\bigr)>\csc(\theta_0/2), \qquad \tanh\bigl(\ell(v,w)/2\bigr)>\cos(\theta_0/2).
\]
Seen from $f(v)$, the angular radius of its ideal cap is
\[
\arccos\bigl(\tanh(\ell(v,w)/2)\bigr)<\theta_0/2,
\]
so the outgoing closed ideal balls at a common lifted vertex are pairwise disjoint, and by the
argument in the proof of Proposition \ref{gap} so are their closed hyperbolic halfspaces. Along a
tree path $v_0,...,v_k$,
\[
E(v_j,v_{j+1})\subset\inte E(v_{j-1},v_j),
\]
since the new forward halfspace is disjoint from $E(v_j,v_{j-1})$, the opposite side of the
preceding midpoint plane. The inclusion is also strict on the ideal boundaries, and induction gives
the nesting.

Let $\pi_X:\widetilde X\to X$ be the universal covering mapping, choose a spanning tree $\mathcal
S\subset X$ and a component $S$ of $\pi_X^{-1}(\mathcal{S})$, hence $\pi_X|_S:S\to\mathcal{S}$ is an
isomorphism. There are $2r$ oriented edges leaving $S$, two for each edge of $X$ outside $\mathcal
S$. For an oriented lifted edge $e$, denote its ideal ball by $B_e$, and let $d_1, d_2$ leave $S$.
If their initial vertices agree, local disjointness applies. Otherwise, let their initial vertices
be $\tau_0,\tau_k$ and let
\[
e_j:\tau_{j-1}\to\tau_j,\qquad 1\le j\le k,
\]
be the path between them in $S$. Since $d_1$ leaves $S$ whereas $e_1$ lies in $S$, these are
distinct edges starting at $\tau_0$, so
\[
B_{d_1}\cap B_{e_1}=\emptyset.
\]
At $\tau_j$, by local disjointness and the opposite-side identity we have
\[
B_{e_{j+1}}\cap B_{\bar e_j}=\emptyset, \qquad B_{\bar e_j}=\mathbb{S}^d\setminus\inte B_{e_j},
\]
so we have
\[
B_{e_{j+1}}\subset\inte B_{e_j},
\]
and the same argument at $\tau_k$ implies $B_{d_2}\subset\inte B_{e_k}$. Therefore
\[
B_{d_2}\subset\inte B_{e_k} \subset\cdots\subset\inte B_{e_1},
\]
and consequently
\[
B_{d_1}\cap B_{d_2} \subset B_{d_1}\cap B_{e_1} =\emptyset.
\]
We orient an edge $e:i\to j$ outside $\mathcal{S}$, and let $\widetilde i,\widetilde j\in S$ be the
unique lifts of its endpoints. The lift of $e$ starting at $\widetilde i$ ends at $g_e\widetilde j$
for a unique $g_e\in\G$, and the corresponding lift of $\bar e$ is
\[
\widetilde e^+: \widetilde i\longrightarrow g_e\widetilde j, \qquad \widetilde e^-: \widetilde
j\longrightarrow g_e^{-1}\widetilde i.
\]
Hence we have
\[
g_e\widetilde e^-=\overline{\widetilde e^+}.
\]

Let $E_e^\pm$ be the closed halfspaces determined by their midpoint planes, directed towards their
terminal endpoint images, and set $B_e^\pm=\partial_\infty E_e^\pm$. By equivariance we have
\[
g_eE_e^-=\mathbb{H}^N\setminus\inte E_e^+, \qquad g_eB_e^-=\mathbb{S}^d\setminus\inte B_e^+.
\]
Hence we have
\[
g_e(\mathbb{S}^d\setminus\inte B_e^-)=B_e^+.
\]

We choose a vertex $o\in X$ and, for each vertex $v$, let $\sigma_v$ be the path in $\mathcal{S}$
from $o$ to $v$. The loops
\[
\beta_e=\sigma_i e\sigma_j^{-1}, \qquad e\in E(X)\setminus E(\mathcal{S}),
\]
form a free basis of $\pi_1(X,o)$, and their lifts starting at $\widetilde o\in S$ end at
$g_e\widetilde o$, so the carrier isomorphism sends $\beta_e$ to $g_e$ and these elements freely
generate $\G$.

For one such $g$, let $E_+$ be the halfspace of the first lift and $E_-$ that of the second. Since
their closed ideal balls are disjoint,
\[
\alpha_g=\dis(\partial E_-,\partial E_+)>0.
\]
By the pairing, $gE_-=\mathbb{H}^N\setminus\inte E_+$ and $gE_+\subset\inte E_+$. We choose $o$
strictly outside both and set $\beta_g=\min\{\dis(o,\partial E_-),\dis(o,\partial E_+)\}>0$. The
boundary planes of $E_+,gE_+,...,g^{n-1}E_+$ are nested, consecutive planes having gap $\alpha_g$
since $g^{-1}\partial E_+=\partial E_-$, and the segment from $o$ to $g^no$ crosses them in order
with initial and terminal distances at least $\beta_g$. The last plane is $g^{n-1}\partial E_+$, and
\[
\begin{aligned}
\dis(g^no,g^{n-1}\partial E_+) &=\dis(o,g^{-1}\partial E_+)\\
&=\dis(o,\partial E_-)\ge\beta_g.
\end{aligned}
\]
Therefore we have
\[
\dis(o,g^no)\ge2\beta_g+(n-1)\alpha_g\qquad(n\ge1).
\]
By subadditivity the stable translation length $\lim_n \dis(o,g^no)/n$ exists, and by this lower
bound it is positive, hence each of the chosen generators is loxodromic. Together with the disjoint
paired round balls already obtained, this proves that $\G$ is classical.
\end{proof}

An element is called primitive if it belongs to a free basis. Theorem \ref{global} below replaces
the rank-two long-generator reduction in arbitrary rank.

\begin{lem}\label{minimum}
Let $\G$ be a free convex cocompact Kleinian group of rank $r\ge2$. Then it has a metric carrier
graph of minimum total length among all such graphs of rank $r$. After leaf edges are deleted,
vertices of valence two suppressed and edges of length zero collapsed, this graph has at most $3r-3$
edges, all of positive length.
\end{lem}
\begin{proof}
Nearest-point projection to the convex hull is equivariant, does not increase length, and is
equivariantly homotopic to the identity along joining geodesics, so it is sufficient to work in the
compact convex core. Deleting leaves, suppressing vertices of valence two and straightening each
resulting edge in its relative homotopy class preserves the rank and the group isomorphism without
increasing length, and zero-length edges are allowed until the end. Note that every carrier graph
$(X_0,\bar f_0)$ obtained after these reductions satisfies
\[
2|E(X_0)|\ge3|V(X_0)|,\qquad r=|E(X_0)|-|V(X_0)|+1.
\]
Therefore we have
\[
|E(X_0)|\le3r-3,\qquad |V(X_0)|\le2r-2.
\]
We fix one such carrier and set
\[
A_0=\leg(X_0) =\sum_{e\in E(X_0)}L_e.
\]
By convex cocompactness, choose a compact set $C\subset\Hu(\Lambda_\G)$ such that
\[
\Hu(\Lambda_\G) =\bigcup_{\g\in\G}\g C.
\]
For any carrier $(X,\bar f)$ obtained after the preceding reductions with $\leg(X)\le A_0$, set
$v_X=|V(X)|$, index its vertices by $i=1,...,v_X$, choose a lift $p_i\in C$ of each vertex image,
choose one orientation of each edge $e\in E(X)$, and denote by $i(e),j(e)$ its initial and terminal
vertex indices, hence its lifted image joins $p_{i(e)}$ to $\g_ep_{j(e)}$, the reverse orientation
having a lifted image joining $p_{j(e)}$ to $\g_e^{-1}p_{i(e)}$, and lengths are summed once over
$e\in E(X)$. Set
\[
\mathcal{L}(A_0)= \left\{\g\in\G: \inf_{p,q\in C}\dis(p,\g q)\le A_0\right\}.
\]
The set
\[
\mathcal{N}_{A_0}(C) =\{z\in\mathbb{H}^N:\dis(z,C)\le A_0\}
\]
is compact, since $C$ is, and the infimum defining $\mathcal{L}(A_0)$ is attained on $C\times C$, so
\[
\mathcal{L}(A_0) =\{\g\in\G: \g C\cap\mathcal{N}_{A_0}(C)\ne\emptyset\},
\]
which is finite since the action of $\G$ on $\mathbb{H}^N$ is proper and $C$ and $\mathcal
N_{A_0}(C)$ are compact. If the total length is at most $A_0$, then $\g_e\in\mathcal{L}(A_0)$ for
all edge, so with the vertices and the oriented edges indexed, the number of graph and label choices
is at most
\[
\sum_{n_V=1}^{2r-2}\sum_{n_E=1}^{3r-3} n_V^{2n_E}|\mathcal{L}(A_0)|^{n_E}<\infty.
\]
For each fixed graph $X$, the space
\[
C^{v_X}:=\prod_{i=1}^{v_X}C =\{(p_1,...,p_{v_X}):p_i\in C\}
\]
is compact, and for each choice of edge labels inducing an isomorphism, the function
\[
\Phi_X(p_1,...,p_{v_X}) :=\sum_{e\in E(X)}\dis(p_{i(e)},\g_ep_{j(e)})
\]
satisfies, for $\mathbf p=(p_1,...,p_{v_X})$ and $\mathbf q=(q_1,...,q_{v_X})$ in $C^{v_X}$,
\[
\begin{aligned}
|\Phi_X(\mathbf p)-\Phi_X(\mathbf q)| &\le\sum_{e\in E(X)} \bigl(\dis(p_{i(e)},q_{i(e)})
+\dis(p_{j(e)},q_{j(e)})\bigr)\\ &=\sum_{i=1}^{v_X}\deg(i)\dis(p_i,q_i).
\end{aligned}
\]
Here $\deg(i)$ is the valence of vertex $i$, which implies that $\Phi_X$ is continuous and attains
its minimum on $C^{v_X}$. For a loop $e_1^{\varepsilon_1}\cdots e_k^{\varepsilon_k}$ starting and
ending at a chosen vertex, the induced homomorphism assigns
\[
[e_1^{\varepsilon_1}\cdots e_k^{\varepsilon_k}] \longmapsto \g_{e_1}^{\varepsilon_1}\cdots
\g_{e_k}^{\varepsilon_k}, \qquad \varepsilon_j\in\{-1,1\},
\]
a product depending only on the edge labels, so moving the vertices and joining the corresponding
lifted endpoints by geodesics leaves the homomorphism $\pi_1(X)\to\G$ unchanged, and it remains an
isomorphism even when some edges map to a point. Let $A_*$ be the smallest attained value. Note that
every carrier of length at most $A_0$ reduces to one of these choices, so
\[
A_*\le A_0,\qquad A_*\le\leg(Y)\quad \text{for every carrier }Y\text{ with }\leg(Y)\le A_0.
\]
Here $\leg(Y)$ is its total length, and if $\leg(Y)>A_0$, then $\leg(Y)>A_*$, hence $A_*$ is the
minimum over all carriers.

If the zero-length edges contained a circle, its lift would have a constant image fixed by its
nonidentity group label. Since a point stabilizer is compact and the group is discrete, two powers
of each stabilizing element coincide, so such an element has a positive power equal to the identity,
and torsion-freeness makes the stabilizer trivial, contrary to injectivity, hence each connected
component of the subgraph formed by the zero-length edges is a tree, and collapsing each such tree
preserves rank, total length and the induced isomorphism on fundamental groups, leaving only edges
of positive length. Collapsing a tree with $k$ vertices produces a vertex of valence at least
$3k-2(k-1)=k+2\ge3$, so the resulting positive-length carrier has at most $3r-3$ edges and attains
the minimum among all carriers.
\end{proof}

\begin{lem}\label{midpoint}
Let $u,v,x,y\in\mathbb{H}^N$, with $x,y$ at distances $R_x,R_y>0$ on geodesic rays from $u,v$,
respectively. Assume that the outward closed ideal balls of the perpendicular hyperplanes at $x,y$
meet. Then
\[
\dis(x,y)<|R_x-R_y|+\dis(u,v)+2.
\]
\end{lem}
\begin{proof}
Let $\xi$ be a common ideal point. The direction from $u$ to $\xi$ makes angle at most
$\arccos(\tanh R_x)$ with the direction of $x$. Let $x_\xi$ be the point at distance $R_x$ on the
ray $u\xi$, and $y_\xi$ the point at distance $R_y$ on $v\xi$. By the cosine rule,
\[
\cosh\dis(x,x_\xi) \le1+\sinh^2R_x(1-\tanh R_x)<3/2,
\]
so $\dis(x,x_\xi)<1$, and likewise $\dis(y,y_\xi)<1$. We send $\xi$ to infinity in upper half-space.
From $u=(z_u,t_u)$ and $v=(z_v,t_v)$, the unit-speed rays are
\[
\rho_u(s)=(z_u,e^s t_u),\qquad \rho_v(s)=(z_v,e^s t_v).
\]
For $s\ge0$,
\[
\cosh\dis(\rho_u(s),\rho_v(s)) =1+\frac{(t_u-t_v)^2}{2t_ut_v} +e^{-2s}\frac{|z_u-z_v|^2}{2t_ut_v}
\le\cosh\dis(u,v).
\]
Comparing the rays at $s=\min\{R_x,R_y\}$ and adding the remaining segment on the longer ray, we get
\[
\dis(x_\xi,y_\xi)\le\dis(u,v)+|R_x-R_y|,
\]
and since each of the other two distances is below $1$, the inequality follows.
\end{proof}

\begin{thm}\label{global}
Let $\G$ be a free convex cocompact Kleinian group of rank $r\ge2$, and set
\[
m=3r-3,\qquad B=2m+5=6r-1,\qquad C_r=mB^m+2.
\]
Then $\G$ is classical, or it has a primitive element $a$ such that $\ell(a)<C_r$. In particular, if
every primitive element has length at least $C_r$, then $\G$ is classical.
\end{thm}
\begin{proof}
Assume that
\[
\ell(a)\ge C_r\qquad\text{for every primitive }a\in\G
\]
and let $\bar f:X\to\mathbb{H}^N/\G$ be a shortest geodesic carrier, which by Lemma \ref{minimum}
exists and has at most $m$ edges, and let $f:\widetilde X\to\mathbb{H}^N$ be its $\G$-equivariant
lift and $p:\widetilde X\to X$ the covering projection. For an edge $e$ of $X$ and any lift
$\widetilde e$, set
\[
L_e=\leg(f(\widetilde e)),\qquad L(X)=\sum_{e\in E(X)}L_e,
\]
where each edge is given the length $L_e$, which does not depend on its orientation, and $L(Y)$
denotes the total length of any carrier $Y$. Since $B>1$, the intervals
\[
(B^j,B^{j+1}),\qquad j=0,...,m,
\]
are pairwise disjoint, each edge length belongs to at most one of them, and since $|E(X)|\le m$, at
least one of these $m+1$ intervals contains no edge length. Set $T=B^j$ and $S=mT$, so that
\[
1\le B^j \le B^m,\qquad mB^j\le C_r-2,\qquad L_e\le B^j\ \text{or}\ L_e\ge B^{j+1}.
\]
An edge is called short if $L_e\le B^j$ and long if $L_e\ge B^{j+1}$, and $F\subset X$ denotes the
subgraph consisting of all vertices and the short edges. If $F$ contains an embedded circle $c$,
choose an edge $e_0$ of $c$. Since $c\setminus\inte e_0$ is a path, it is contained in a spanning
tree $\mathcal{S}$ of $X$, and $e_0$ lies outside $\mathcal{S}$. We choose an orientation $e_0:v\to
w$ and orient $c$ so that one traversal follows $e_0$ and then the path $[w,v]_{\mathcal{S}}$. If
$\sigma$ is the path in the spanning tree from the base vertex to the initial vertex of $c$, the
corresponding basis element is represented by $\sigma c\sigma^{-1}$, with $c$ traversed once from
that vertex. We choose a lift $\widetilde z$ of a point on $c$ and let $a\in\G$ be the deck element
such that the lift of $c$ from $\widetilde z$ ends at $a\widetilde z$. Then $a$ is primitive and
\[
C_r\le\ell(a) \le\dis(f(\widetilde z),af(\widetilde z)) \le\sum_{e\in E(c)}L_e \le mB^j\le C_r-2,
\]
a contradiction, hence every component of $F$ is a tree, and for every component $F_0$ and every
lift $\widetilde F_0$,
\[
\sum_{e\in E(F_0)}L_e\le mB^j, \qquad \diam f(\widetilde F_0)\le mB^j.
\]

For an oriented long edge $\widetilde e:u\to w$ with midpoint $x$, set $R_{\widetilde
e}=L_{p(\widetilde e)}/2$, let $\mathcal{P}_{\widetilde e}$ be the hyperplane perpendicular to
$[f(u),f(w)]$ at $f(x)$, $E_{\widetilde e}$ its closed halfspace toward $f(w)$, and $B_{\widetilde
e}=\partial_\infty E_{\widetilde e}\subset\mathbb{S}^d$. For the reverse edge $\overline{\widetilde
e}$,
\[
E_{\overline{\widetilde e}}=\mathbb{H}^N\setminus\inte E_{\widetilde e}, \qquad
B_{\overline{\widetilde e}}=\mathbb{S}^d\setminus\inte B_{\widetilde e}.
\]
Since
\[
R_{\widetilde e}-S-2 \ge\frac{B^{j+1}}{2}-mB^j-2 =\frac52B^j-2\ge\frac12,
\]
the image of the short component containing $u$ is disjoint from $E_{\widetilde e}$, since for $z$
in that lifted component
\[
\dis(f(z),f(u))\le mB^j<R_{\widetilde e} =\dis(f(u),\mathcal{P}_{\widetilde e}).
\]

We fix a component $F_0$ of $F$ and one lift $\widetilde F_0$. The balls of distinct long edges
leaving $\widetilde F_0$ are disjoint, i.e. $B_{\widetilde e}\cap B_{\widetilde e'}=\emptyset$.
Suppose not, and let the two edges be $\widetilde e_1:u\to u_1$ and $\widetilde e_2:v\to v_1$, with
midpoints $x,y$ and half-lengths $R_x\le R_y$. Since $u,v$ lie in one lifted short tree,
\[
\dis(f(u),f(v))\le mB^j,
\]
and by Lemma \ref{midpoint} we have
\[
D:=\dis(f(x),f(y)) <R_y-R_x+mB^j+2 =R_y-\eta,\qquad \eta:=R_x-mB^j-2\ge\frac12.
\]

Let us first assume that the two edges project to different unoriented edges of $X$. We insert a
vertex at $p(x)$ in the edge $p(\widetilde e_1)$, replacing this edge by two edges, and denote the
resulting graph by $X_*$. The path
\[
P=p\bigl([v,u]_{\widetilde F_0}[u,x]_{\widetilde X}\bigr):p(v)\longrightarrow p(x),
\]
avoids $e_2=p(\widetilde e_2)$. We replace $e_2:p(v)\to p(v_1)$ by $e':p(x)\to p(v_1)$, denote the
new graph by $X'$, and define mappings $q:X'\to X_*$ and $t:X_*\to X'$ that fix all vertices and
every retained edge and satisfy
\[
q(e')=P^{-1}e_2,\qquad t(e_2)=Pe'.
\]
Since $q,t$ fix $P$, relative to endpoints
\[
qt(e_2)=PP^{-1}e_2\simeq e_2,\qquad tq(e')=P^{-1}Pe'\simeq e',
\]
so $q,t$ are homotopy inverses. Regarding $\bar f$ as a mapping on the subdivided graph $X_*$, take
the carrier $\bar f\circ q$ and straighten only $e'$. Its chosen lift joins $f(x)$ to $f(v_1)$, so
\[
\begin{aligned}
L_{e'}&=\dis(f(x),f(v_1))\\
&\le\dis(f(x),f(y))+\dis(f(y),f(v_1))\\
&=D+R_y<2R_y-\eta=L_{e_2}-\eta.
\end{aligned}
\]
Subdivision preserves the total length of the first edge and all other edge lengths are unchanged,
so
\[
L(X')=L(X)-L_{e_2}+L_{e'} <L(X)-\eta\le L(X)-\frac12.
\]
By the injectivity argument in Lemma \ref{minimum}, each connected component of the subgraph formed
by zero-length edges is a tree, and collapsing these trees preserves rank, total length and the
induced isomorphism on fundamental groups. This gives a carrier graph of rank $r$ with positive edge
lengths and total length $L(X')<L(X)$, contradicting the minimality of $X$ among all carrier graphs
of rank $r$.

Now assume that the two edges project to the same unoriented edge. Since $p|_{\widetilde F_0}$ is an
isomorphism onto $F_0$, each oriented long edge with its initial vertex in $F_0$ has one lift
starting in $\widetilde F_0$, so the two distinct leaving edges have opposite orientations. Let
their quotient edge be $e:i\to j$, with $i,j\in F_0$, and choose its endpoint lifts
$u,v\in\widetilde F_0$. For some $\g\in\G$, the leaving edges are
\[
\widetilde e_1:u\longrightarrow\g v, \qquad \widetilde e_2:v\longrightarrow\g^{-1}u.
\]
Here $\g\ne1$, since otherwise $\widetilde e_1$ and the short-tree path from $v$ to $u$ would form a
circle in $\widetilde X$, and
\[
\g\widetilde e_2=\overline{\widetilde e_1},\qquad R_x=R_y,\qquad \g y=x.
\]
Now combine equivariance with the overlap estimate to get
\[
\ell(\g) \le\dis(f(y),\g f(y)) =\dis(f(y),f(x))<mB^j+2\le C_r.
\]
Let $c=e[j,i]_{F_0}$, with a constant short path when $i=j$. Its lift
\[
\widetilde c=\widetilde e_1\,\g[v,u]_{\widetilde F_0}: u\longrightarrow\g v\longrightarrow\g u,
\]
shows that $c$ represents $\g$ up to conjugacy. We extend $F$ to a spanning tree of $X$. Since this
tree contains $[j,i]_{F_0}$ and does not contain $e$, its graph basis contains the class of $c$ up
to conjugacy, which implies that $\g$ is primitive, which is a contradiction. So we have proved
\begin{equation}\label{disjoint}
B_{\widetilde e}\cap B_{\widetilde e'}=\emptyset \quad\text{for distinct long edges leaving the same
lifted short component.}
\end{equation}
We collapse each component of $F$ in $X$ and each lifted component in $\widetilde X$, obtaining
$H=X/F$ and its universal covering tree $\mathcal{U}$. Since $V(F)=V(X)$ and each connected
component of $F$ is a tree, the number of these components is $|V(X)|-|E(F)|$, so by collapsing each
component to a vertex we get
\[
|V(H)|=|V(X)|-|E(F)|,\qquad |E(H)|=|E(X)|-|E(F)|,
\]
and
\[
|E(H)|-|V(H)|+1 =|E(X)|-|V(X)|+1=r.
\]
We identify the edges of $\mathcal{U}$ with the corresponding long edges of $\widetilde X$, keeping
their original midpoint balls. For consecutive oriented edges $\widetilde e,\widetilde e'$ in a
reduced path of $\mathcal{U}$, the edges $\overline{\widetilde e},\widetilde e'$ leave the same
vertex and are distinct, so
\[
B_{\widetilde e'}\cap B_{\overline{\widetilde e}}=\emptyset, \qquad B_{\widetilde e'}\subset
\mathbb{S}^d\setminus B_{\overline{\widetilde e}} =\inte B_{\widetilde e}.
\]

We lift a spanning tree $\mathcal{T}\subset H$ to $\widetilde{\mathcal{T}}\subset\mathcal{U}$. The
number of oriented edges leaving $\widetilde{\mathcal{T}}$ is
\[
2\bigl(|E(H)|-(|V(H)|-1)\bigr)=2r.
\]
For two such edges $\widetilde a,\widetilde b$ with the same initial vertex, the balls are disjoint
by (\ref{disjoint}). Otherwise, if $\widetilde e_1\cdots\widetilde e_k$ is the path in
$\widetilde{\mathcal{T}}$ between their initial vertices, by the nesting above we have
\[
B_{\widetilde b}\subset\inte B_{\widetilde e_k} \subset\cdots\subset\inte B_{\widetilde e_1}, \qquad
B_{\widetilde a}\cap B_{\widetilde e_1}=\emptyset.
\]
Hence all $2r$ leaving balls are pairwise disjoint.

We choose one orientation $e_1,...,e_r$ for the edges of $H$ outside $\mathcal{T}$ and, for each
$e_i$, its endpoint lifts $v_i,w_i\in\widetilde{\mathcal{T}}$, hence the two leaving lifts are
\[
\widetilde e_i^+:v_i\longrightarrow g_iw_i, \qquad \widetilde e_i^-:w_i\longrightarrow g_i^{-1}v_i,
\qquad 1\le i\le r.
\]
Set $B_i^\pm=B_{\widetilde e_i^\pm}$. Since $g_i\widetilde e_i^-=\overline{\widetilde e_i^+}$, by
equivariance we have
\[
g_iB_i^-=B_{\overline{\widetilde e_i^+}},\qquad g_i\bigl(\mathbb{S}^d\setminus\inte B_i^-\bigr)
=\mathbb{S}^d\setminus\inte B_{\overline{\widetilde e_i^+}} =B_i^+.
\]
Since collapsing $F$ preserves the carrier isomorphism, the spanning-tree loops give a free basis
$g_1,...,g_r$ of $\G$, each $g_i$ being nontrivial and hence loxodromic by convex cocompactness. The
disjoint closed round balls $B_i^\pm$ satisfy (\ref{pairing}), which implies that $\G$ is classical.
\end{proof}

The next estimate follows Hou's annular generator selection \cite{Hou2020} (Propositions 5.6--5.7
and Definition 6.2), with an arbitrary upper bound $C$ on the distinguished translation length.

\begin{thm}\label{bounded-annulus}
Let $C\ge4$ and $r\ge2$, and let $\G$ be a convex cocompact Kleinian group with free basis
$a,b_1,...,b_{r-1}$ such that
\[
0<L=\ell(a)\le C,\qquad 0<\delta_\G\le\frac{1}{32(C+1)}.
\]
We normalize $a$ and define the positive pole parameters $\nu_i=\sqrt{K_i/(|u_i||v_i|)}>0$ of the
$b_i$ as in Theorem \ref{annulus}. Then $\sum_i\nu_i<L/8$, and there are independent two-sided
changes by powers of $a$ after which every paired pole ball lies strictly inside one $a$-annulus.
\end{thm}
\begin{proof}
We fix $R_0>0$ and apply independent two-sided changes by powers of $a$ so that every pole norm lies
in $[R_0,e^LR_0]$. The changed quantities are written without primes, and $\nu_i$ is unchanged. Set
\[
\widehat p_i=(0,\sqrt{K_i}),\qquad \widehat D_i=\dis(\widehat p_i,b_i\widehat p_i).
\]
Since $\widehat p_i$ lies on the axis of $a$, by the pole bound we have
\[
\cosh\widehat D_i\le1+\frac{1+\cosh L}{\nu_i^2} +\frac{1}{2\nu_i^4}.
\]
If $\nu_i\ge1/2$, then $\cosh\widehat D_i\le13+4\cosh L\le17e^L$, which implies $\widehat
D_i<L+\log34<C+4$. For $t>0$ let
\[
F(t)=\log\frac{e^t+3}{e^t-1}.
\]
Using the full group exponent in (\ref{pair-pressure}), we get $\widehat D_i\ge F(\delta_\G
L)/\delta_\G$. Since $\delta_\G L\le1/32$ and $e^t-1\le2t$ in this range,
\[
F(\delta_\G L)\ge\log\frac{2}{\delta_\G L}\ge\log64,
\]
so $\widehat D_i\ge32(C+1)\log64>C+4$, which is a contradiction. It follows that each $\nu_i<1/2$.

At the common point $p_0=(0,R_0)$, set $D_i=\dis(p_0,b_i p_0)$, so that
\[
\cosh D_i \le\frac{2\cosh^2L}{\nu_i^2} +\frac{\nu_i^2e^{2L}}2+e^{2L} <\frac{4e^{2C}}{\nu_i^2},
\qquad D_i<\log8+2C-2\log\nu_i.
\]
The count of the distinct words
\[
a^{n_1}b_{i_1}\cdots a^{n_k}b_{i_k}, \qquad n_j\in\mathbb{Z},
\]
used in the proof of Theorem \ref{annulus} shows that
\[
\sum_i e^{-\delta_\G D_i}\le\tanh(\delta_\G L/2).
\]
Since $D_i<\log8+2C-2\log\nu_i$,
\[
\begin{aligned}
\sum_i\nu_i^{2\delta_\G} &<e^{\delta_\G(\log8+2C)} \sum_i e^{-\delta_\G D_i}\\
&\le e^{\delta_\G(\log8+2C)}\tanh(\delta_\G L/2)\\
&<e^{\delta_\G(\log8+2C)}\,\delta_\G L/2,
\end{aligned}
\]
and for $p=1/(2\delta_\G)\ge16(C+1)$, by subadditivity we have
\[
\begin{split}
\sum_i\nu_i &<\sqrt8\,e^C(\delta_\G L/2)^p\\
&=\sqrt8\,e^C(L/C)^p(\delta_\G C/2)^p\\
&\le\sqrt8\,e^C64^{-16(C+1)}\,L/C\\
&<\tfrac12 L/C\le L/8,
\end{split}
\]
using $L/C\le1$ and $\delta_\G C/2<1/64$, together with
\[
\sqrt8\,e^C64^{-16(C+1)}<\frac12,
\]
which holds for $C\ge4$ since
\[
\log\bigl(\sqrt8\,e^C64^{-16(C+1)}\bigr) =\frac12\log8+C-16(C+1)\log64
\]
has derivative $1-16\log64<0$ with respect to $C$ and value $\tfrac12\log8+4-80\log64<-\log2$ at
$C=4$.

Next, $s_0\in\mathbb{R}$ is to be chosen such that, for all $i$ and $k\in\mathbb{Z}$,
\[
|s_0-\log|u_i|-kL|>2\nu_i,\qquad |s_0-\log|v_i|-kL|>2\nu_i.
\]
For each $i$, the choices violating these inequalities form two closed arcs on
$\mathbb{R}/L\mathbb{Z}$, each of length $4\nu_i$. The union of these arcs over all $i$ has length
at most $8\sum_i\nu_i<L$. Since the circle has length $L$, we can choose $s_0\in\mathbb{R}$ whose
image lies outside these arcs. By Lemma \ref{pole}, we can choose $k_i,l_i\in\mathbb{Z}$ such that,
for $b_i'=a^{k_i}b_i a^{l_i}$,
\[
\{\log|u_i'|,\log|v_i'|\} \subset(s_0+2\nu_i,s_0+L-2\nu_i),
\]
and $\nu_i$ is unchanged by these choices. Since $0<\nu_i<1/2$,
\[
-\log(1-\nu_i)<2\nu_i,\qquad \log(1+\nu_i)<2\nu_i.
\]
Hence, for $z=u_i'$ or $z=v_i'$,
\[
e^{s_0}<(1-\nu_i)|z| \le(1+\nu_i)|z|<e^{s_0+L},
\]
so every closed ball centered at $z$ with radius $\nu_i|z|$ lies strictly inside the common annulus
\[
e^{s_0}<|x|<e^{s_0+L}.
\]
Finally,
\[
(\nu_i|v_i'|)(\nu_i|u_i'|) =\nu_i^2|u_i'||v_i'|=K_i',
\]
so each pair satisfies the pairing condition for $b_i'$.
\end{proof}

\begin{cor}\label{reduction}
Let $r\ge2$, set
\[
C_r=(3r-3)(6r-1)^{3r-3}+2,\qquad \varepsilon_r=\frac{1}{32(C_r+1)}.
\]
Let $\G<\Isom(\mathbb{H}^N)$ be discrete, torsion-free, and non-elementary, and suppose its least
number of generators is $r$. Assume $\delta_\G\le\varepsilon_r$. Then $\G$ is classical, or there is
a free basis $a,b_1,...,b_{r-1}$ with $\ell(a)<C_r$ for which all other paired pole balls lie
strictly in one $a$-annulus after two-sided changes, and their positive parameters $\nu_i$, as in
Theorem \ref{bounded-annulus}, satisfy $\sum_i\nu_i<\ell(a)/8$. In particular, the conclusion holds
under $\mathfrak{D}_\G\le\varepsilon_r$.
\end{cor}
\begin{proof}
Since $\varepsilon_r<1/2$, the structural theorem makes $\G$ free and convex cocompact with
$\delta_\G=\mathfrak{D}_\G$, and its free rank is its least generating number $r$. By Theorem
\ref{global}, either $\G$ is classical or it has a primitive element $a$ with $\ell(a)<C_r$. In the
latter case choose a free basis $a,b_1,...,b_{r-1}$ and apply Theorem \ref{bounded-annulus} with
$C=C_r$. The Hausdorff-dimension conclusion follows since $\delta_\G\le \mathfrak{D}_\G$.
\end{proof}

\begin{rem}
In rank two, $C_2=3995$ and $\varepsilon_2=1/(32\cdot3996)<1/3$, so Theorem \ref{main-2} already
gives classical Schottky. Theorem \ref{bounded-annulus} places the balls in one annulus but says
nothing about their mutual disjointness.
\end{rem}

\begin{prop}\label{short}
Let $X$ be a carrier graph with mapping $\bar f:X\to\mathbb{H}^N/\G$ and $E=|E(X)|$ edges, and fix
$T>0$. If
\[
\delta_\G<\frac{\log3}{2ET},
\]
every connected component formed by its edges of length at most $T$ is a tree or has fundamental
group of rank one. If, in addition, $X$ has rank $r\ge2$, $E\le3r-3$, and $T=4$, the sufficient
bound is $\delta_\G<\log3/(24(r-1))$.
\end{prop}
\begin{proof}
Let $f:\widetilde X\to\mathbb{H}^N$ be the equivariant lift, and suppose that a connected component
$Y$ of the subgraph consisting of all vertices of $X$ and its edges with $L_e\le T$ has
$\operatorname{rank}\pi_1(Y)\ge2$. We choose a root $z\in V(Y)$ and a spanning tree $K$ of $Y$, and
extend $K$ to a spanning tree of $X$. The associated basis contains the loops given by the edges of
$Y$ outside $K$, joined to $z$ within $K$. With a lift $\widetilde z$ of $z$ used to identify
$\pi_1(X,z)$ with $\G$, inclusion induces an injection $\pi_1(Y,z)\to\pi_1(X,z)\cong\G$. We choose
two such edges, whose free generators map under the carrier isomorphism to noncommuting elements
$a,b\in\G$. Each loop follows two tree paths and one further edge, so with $E_Y=|E(Y)|\le E$ it has
length at most $(2E_Y-1)T\le2ET$, and at $p=f(\widetilde z)$
\[
\dis(p,ap),\dis(p,bp)\le2ET.
\]
By Hou's two-generator pressure inequality in Lemma \ref{pressure} we have
\[
1\ge \frac{2}{1+e^{\delta_\G\dis(p,ap)}}+ \frac{2}{1+e^{\delta_\G\dis(p,bp)}}
\ge\frac{4}{1+e^{2\delta_\G ET}}.
\]
Hence $\delta_\G\ge\log3/(2ET)$, which is a contradiction.
\end{proof}


\section{Classical in Arbitrary Rank}
In this section we prove Theorem \ref{main}. As in Hou's strategy, a small exponent is used to
choose generators and paired regions \cite{Hou2020} (Section 1.1). In three dimensions, Hou resolves
circle degenerations and derives a rank-independent bound from the fixed-rank case \cite{Hou2020}
(Sections 6--15). In higher dimensions, we use a minimizing carrier and enlarge its rank-one
components to get compatible round balls, with a constant depending on $r$.

For $r\ge2$, set
\[
m=3r-3,\qquad B=2m+5,\qquad C_0=mB^m,\qquad A=2\cdot203^m,\qquad \lambda_r=\frac{1}{1000AC_0}.
\]
Carrier graphs are as in Section $6$. A carrier is geodesic if each edge maps to the geodesic in its
relative homotopy class, and each loop contributes two to the valence of its vertex. The hyperboloid
model carries the Lorentz form of signature $(-,+,...,+)$.

\begin{thm}\label{rank-r}
Let $\G<\Isom(\mathbb{H}^N)$ be discrete, free, and convex cocompact of rank $r\ge2$. If its
critical exponent satisfies $\delta_\G\le\lambda_r$, then $\G$ is classical Schottky.
\end{thm}
The main steps of the proof are the following. Let $G$ be a carrier graph of minimum total length,
with at most $m$ edges. The gap between short and long edge lengths gives a subgraph $G_0$ whose
components have length at most $C_0$, and the displacement estimate (\ref{pairs}) excludes a
component of rank at least two, since such a component would imply
\[
\log 3\le 2\delta_{\G} C_0 \le \frac{2}{1000A}<\log 3,
\]
so each component of $G_0$ is a tree or has rank one. For a rank-one component with loop generator
$a_i$, axis $A_i$ and $L_i=\ell(a_i)$, by the estimate (\ref{anchor}) we have
\[
\dis(p_i,A_i)\le U_i, \qquad U_i=C_0+\log^+(1/L_i).
\]
These components are enlarged by adding edges and tree components. The growth bound (\ref{growth})
and the axis displacement estimate (\ref{cycle}) prevent both the creation of a second loop and the
joining of two rank-one components, so the resulting components remain disjoint and satisfy
\[
Z_i=\leg(K_i)+U_i\le AU_i, \qquad L_e>200Z_i
\]
for every remaining edge to $K_i$.

For a lifted remaining edge with endpoint images $p,q$, let $\mathcal{B}$ be the closed ideal ball
on the $q$ side of the perpendicular midpoint plane $\mathcal{P}$ of $[p,q]$, and set
$2R=\dis(p,q)>200Z_i$. By minimality we have
\[
\dis(a_i^kp,q)\ge\dis(p,q) \qquad(k\in\mathbb{Z}),
\]
which bounds the distance along $A_i$ between the nearest points to $p$ and $q$, and together with
the distance bound $Z_i$ from the initial endpoint to $A_i$ we get
\[
D=\dis(A_i,\mathcal{P})>98U_i-2,
\]
so that, in coordinates with axis $(0,\infty)$,
\[
\{\log|\zeta|:\zeta\in\mathcal{B}\} \subset I_{\mathcal{B}} =[\tau-\omega(D),\tau+\omega(D)], \qquad
\omega(D)\le4e^{-D}<\frac{L_i}{8m}.
\]
Since $a_i$ multiplies the boundary radius by $e^{L_i}$, the intervals for $\mathcal{B}$ and
$a_i^k\mathcal{B}$ differ by $kL_i$, and as their lengths are less than $L_i$, distinct translates
of the same ball are disjoint. If balls from different edge orbits intersect, the interval bounds
and (\ref{phase}) give $\dis(p_1,p_2)<4Z_i+3$, and by the overlap estimate (\ref{overlap}) we have
\[
\dis(x_1,x_2) <|R_1-R_2|+4Z_i+5 <\max\{R_1,R_2\}.
\]
Here $x_1,x_2$ are the edge midpoints. For distinct edges this permits an edge replacement that
shortens the carrier, and for opposite orientations of one edge, retaining one half and attaching a
shorter loop gives the same contradiction, both replacements preserving the induced isomorphism,
hence all these closed balls and their translates are pairwise disjoint. There are at most $2m$ edge
orbits, and their intervals have total length less than $L_i/2$, so there is $s_0$ such that, after
translation, the chosen balls all lie in
\[
e^{s_0}<|\zeta|<e^{s_0+L_i},
\]
and are consequently disjoint from the two additional closed balls
\[
\{|\zeta|\le e^{s_0}\}, \qquad \{|\zeta|\ge e^{s_0+L_i}\}\cup\{\infty\},
\]
which are paired by $a_i$. For each remaining tree component, the simpler bound
\[
R\ge BT/2>mT+2
\]
allows the same overlap estimate and shortening argument to prove disjointness of its associated
balls.

The resulting regions are arranged along the $\G$-equivariant tree $\mathcal{T}'$, adjacent regions
having matching boundary planes and complementary halfspaces. Disjointness of the balls associated
with each region therefore gives strict nesting along paths in $\mathcal{T}'$, so the balls
corresponding to edges leaving any connected subtree are pairwise disjoint. A lift of a spanning
tree of $X=\G\backslash\mathcal{T}'$ has
\[
2|E(X)|-2(|V(X)|-1)=2r
\]
such balls, paired by free generators $g_1,...,g_r$ of $\G$ with
\[
g_i(\mathbb{S}^d\setminus\inte B_i^-)=B_i^+,
\]
which are the required disjoint round balls.

\begin{proof}
The first estimate is Hou's two-generator displacement inequality \cite{Hou2001}, equivalently
\cite{Hou2020} (Corollary 3.2), whose proof is repeated here in the form used later for components
whose fundamental groups have rank one. Let $a,b\in\G$ freely generate a subgroup of rank two, and
for $x\in\mathbb{H}^N$ set $d_a=\dis(x,ax)$ and $d_b=\dis(x,bx)$. The words
\[
a^{n_1}b^{k_1}\cdots a^{n_j}b^{k_j},\qquad j\ge1,\quad n_i,k_i\in\mathbb{Z}\setminus\{0\},
\]
are distinct and have displacement at most $\sum_i(|n_i|d_a+|k_i|d_b)$, so for $t>\delta_\G$ the
convergent Poincar\'e series dominates
\[
\sum_{j\ge1}\left[ \left(\sum_{n\ne0}e^{-t|n|d_a}\right)
\left(\sum_{k\ne0}e^{-t|k|d_b}\right)\right]^j =\sum_{j\ge1}\left[
\frac{4}{(e^{td_a}-1)(e^{td_b}-1)}\right]^j,
\]
whose ratio is therefore less than one, and letting $t\downarrow\delta_\G$ gives
\begin{equation}\label{pairs}
(e^{\delta_\G d_a}-1)(e^{\delta_\G d_b}-1)\ge4.
\end{equation}
In particular $\delta_\G>0$. When $x$ lies on the axis of $a$, set $L=\ell(a)>0$ and consider the
words
\[
a^{n_1}b\cdots a^{n_j}b,\qquad j\ge1,\quad n_i\in\mathbb{Z}.
\]
Expanding each $b^q$, $q\ge1$, in the reduced word into $q$ letters $b$, their number determines
$j$, and the power of $a$ immediately before the $i$th letter $b$ determines $n_i$, with $n_i=0$
when that power is absent, hence the reduced word determines $(j,n_1,...,n_j)$, and these words are
distinct. Their displacement is at most $L\sum_i|n_i|+j d_b$, and the geometric-series ratio is
\[
e^{-td_b}\sum_{n\in\mathbb{Z}}e^{-t|n|L} =e^{-td_b}\coth(tL/2)<1\qquad(t>\delta_\G).
\]
Letting $t\downarrow\delta_\G$ and using $\coth u\ge1/u$,
\begin{equation}\label{cycle}
\delta_\G \dis(x,bx)\ge\log\coth(\delta_\G L/2) \ge\log\frac{2}{\delta_\G L}.
\end{equation}

{\bf Step 1:} The minimal carrier. For a graph $Y$, $V(Y)$ is its vertex set and $E(Y)$ its
unoriented edge set. A carrier is rectifiable if each edge mapping is rectifiable. With $L_e$ the
arclength of the mapping on $e$, counted with multiplicity, set
\[
\leg(Y)=\sum_{e\in E(Y)}L_e,
\]
where each unoriented edge is counted once, and for a geodesic edge, $L_e$ is the distance between
its lifted endpoints. By the projection and reduction in Lemma \ref{minimum}, it is sufficient to
minimize over geodesic carriers in $\Hu(\Lambda_\G)/\G$ with all valences at least three,
temporarily allowing zero lengths. For such a graph $X$, with $v_X=|V(X)|$ and $e_X=|E(X)|$,
\[
r=e_X-v_X+1,\qquad 3v_X\le2e_X \quad\Longrightarrow\quad e_X\le m,\qquad v_X\le2r-2.
\]
The compact parametrization of that lemma is as follows. We fix a geodesic carrier of length $D_0$
and a compact $Q\subset\Hu(\Lambda_\G)$ with $\G Q=\Hu(\Lambda_\G)$, index the vertices of $X$ by
$1,...,v_X$, and choose lifts $p_i\in Q$ of their images. For one orientation of each edge $e$, let
$i(e),j(e)$ be its endpoints and $\g_e\in\G$ its label, hence its lifted image joins $p_{i(e)}$ to
$\g_ep_{j(e)}$, and the reverse label is $\g_e^{-1}$. By properness,
\[
\mathcal{L}(D_0)=\{\g\in\G:\dis(Q,\g Q)\le D_0\}, \qquad \#\mathcal{L}(D_0)<\infty.
\]
If $\leg(X)\le D_0$, all labels lie in $\mathcal{L}(D_0)$, and by indexing the edges by $1,...,e_X$
and choosing their ordered endpoints we get
\[
\#\{\text{graph and label choices up to isomorphism}\}
\le\sum_{v=1}^{2r-2}\sum_{e=1}^{m}v^{2e}\bigl(\#\mathcal{L}(D_0)\bigr)^e<\infty.
\]
For each graph, the Cartesian product
\[
Q^{v_X}:=\prod_{i=1}^{v_X}Q =\{(p_1,...,p_{v_X}):p_i\in Q\text{ for }1\le i\le v_X\},
\]
is compact, so for each label choice inducing an isomorphism on fundamental groups, the continuous
function
\[
\Phi_X:Q^{v_X}\longrightarrow\mathbb{R},\qquad \Phi_X(p_1,...,p_{v_X}) =\sum_{e\in
E(X)}\dis(p_{i(e)},\g_ep_{j(e)})
\]
with the sum over the chosen orientations, attains its minimum. For a based edge loop
$e_1^{\varepsilon_1}\cdots e_k^{\varepsilon_k}$, with $\varepsilon_j\in\{-1,1\}$, the group label
\[
\g_{e_1}^{\varepsilon_1}\cdots\g_{e_k}^{\varepsilon_k},
\]
does not depend on the $p_i$, so for fixed edge labels, varying $p_1,...,p_{v_X}\in Q$ preserves the
induced isomorphism $\pi_1(X)\to\G$, and every such choice gives a carrier, possibly with
zero-length edges. Let $M_0$ be the least of these minima. Note that every rectifiable carrier of
length at most $D_0$ projects, reduces and straightens to a choice above without increasing length,
and $M_0\le D_0$, and no carrier longer than $D_0$ can improve $M_0$. By the zero-edge argument in
Lemma \ref{minimum}, the zero edges form a forest, whose collapse preserves rank, total length and
the group isomorphism. This gives a positive-length geodesic carrier $G$ with at most $m$ edges and
\[
\leg(G)=M_0 =\inf\{\leg(Y):Y\text{ is a rectifiable carrier of rank }r\},
\]
so the same minimum applies to subdivided carriers and to carriers before straightening.

Let $\pi_G:\widetilde G\to G$ be the universal covering map, with deck group identified with $\G$,
and let
\[
F:\widetilde G\longrightarrow\Hu(\Lambda_\G)
\]
be the equivariant lift of $\bar F:G\to\Hu(\Lambda_\G)/\G$. For a connected subgraph $K\subset G$
and a component $\widetilde K$ of $\pi_G^{-1}(K)$, a spanning tree of $K$ extends to one of $G$, so
by the graph basis we have
\[
\pi_1(K)\hookrightarrow\pi_1(G)\cong\G, \qquad \operatorname{Stab}_\G(\widetilde K) \cong\pi_1(K),
\]
with the conjugate determined by $\widetilde K$. Note that this stabilizer is trivial for a tree and
generated by a loxodromic element for rank one.

Let $X$ be a subdivision of $G$ and $e:v\to w$ an oriented edge. We delete $e$, retaining all
vertices, to get $J$, choose a path $\sigma:v'\to v$ in $J$, and replace $e$ by $e':v'\to w$,
represented by $\sigma e$ in $X$, to get $X'$. The mappings
\[
\begin{aligned}
\phi:X'\to X,&\qquad \phi|_J=\operatorname{id}, &\phi(e')&=\sigma e,\\
\psi:X\to X',&\qquad \psi|_J=\operatorname{id}, &\psi(e)&=\sigma^{-1}e'
\end{aligned}
\]
satisfy, relative to endpoints,
\[
\phi\psi(e)=\sigma^{-1}\sigma e\simeq e, \qquad \psi\phi(e')=\sigma\sigma^{-1}e'\simeq e',
\]
so $\phi,\psi$ are homotopy inverses. Viewing $\bar F$ on $X$, composing it with $\phi$ and
straightening only $e'$ leaves all retained lengths unchanged, so
\[
\leg(X') =\leg(G)-L_e+L_{e'}.
\]
Each connected component of the subgraph formed by zero-length edges is a tree and may be collapsed
by Lemma \ref{minimum}, and reduction cannot increase length. Therefore
\[
\leg(X')<\leg(G)
\]
is impossible, including after subdivisions.

{\bf Step 2:} The overlap estimate. For $u,q\in\mathbb{H}^N$ with $\dis(u,q)=2R>0$, the midpoint cap
of $[u,q]$ is the closed ideal ball beyond its perpendicular midpoint hyperplane toward $q$. On the
hyperboloid, with $\nu$ the unit tangent at $u$ toward $q$, its outward unit normal is
\[
n=\sinh R\,u+\cosh R\,\nu.
\]
An ideal point has a null representative $u+w$, with $w\in T_u\mathbb{H}^N$ and $\langle
w,w\rangle=1$, and with $\cos\theta=\langle\nu,w\rangle$ its cap condition is
\[
\langle u+w,n\rangle =-\sinh R+\cosh R\cos\theta\ge0 \quad\Longleftrightarrow\quad\cos\theta\ge\tanh
R.
\]
Also,
\[
\sinh^2R(1-\tanh R) =\frac{(1-e^{-2R})^2}{2(1+e^{-2R})}<\frac12.
\]
If two such caps meet, with initial points $u,v$, midpoints $x,y$ and half-lengths $R_1,R_2>0$, then
\begin{equation}\label{overlap}
\dis(x,y)<|R_1-R_2|+\dis(u,v)+2.
\end{equation}
To see this, let $\xi$ be a common ideal point and $x_\xi,y_\xi$ the points at distances $R_1,R_2$
on $u\xi,v\xi$, respectively. The cap angle at $u$ has cosine at least $\tanh R_1$, so
\[
\cosh\dis(x,x_\xi) \le1+\sinh^2R_1(1-\tanh R_1)<\frac32<\cosh1,
\]
so $\dis(x,x_\xi)<1$, and likewise $\dis(y,y_\xi)<1$. With $\xi=\infty$ in upper half-space, write
\[
\begin{gathered}
u=(z_u,h_u),\qquad v=(z_v,h_v),\qquad z_u,z_v\in\mathbb{R}^{N-1},\quad h_u,h_v>0,\\
c_u(t)=(z_u,e^th_u),\qquad c_v(t)=(z_v,e^th_v)\qquad(t\ge0).
\end{gathered}
\]
Then we have
\[
\cosh\dis(c_u(t),c_v(t)) =\frac{h_u^2+h_v^2+e^{-2t}|z_u-z_v|^2}{2h_uh_v} \le\cosh\dis(u,v).
\]
Comparing at $t=\min\{R_1,R_2\}$, we get
\[
\dis(x_\xi,y_\xi)\le\dis(u,v)+|R_1-R_2|,
\]
and the triangle inequality completes the proof of (\ref{overlap}).

{\bf Step 3:} Short components. Since $|E(G)|\le m$, there is $j\in\{0,...,m\}$ such that no edge
length lies in $(B^j,B^{j+1})$. Set $T=B^j$, so that
\[
1\le T\le B^m,\qquad L_e\le T\ \text{or}\ L_e\ge BT.
\]
Such edges are called short and long, respectively. Let $G_0$ consist of all vertices of $G$ and its
short edges. Note that every component $K\subset G_0$ satisfies
\[
\leg(K)\le mT\le C_0.
\]
If $\operatorname{rank}\pi_1(K)\ge2$, extend a spanning tree of $K$ to one of $G$. Two associated
loops based at a vertex of $K$ give free generators $a,b$ of a subgroup, each represented by a loop
of length at most $2mT$. At the image $z$ of a common lifted base point,
\[
\dis(z,az),\dis(z,bz)\le2mT.
\]
By (\ref{pairs}),
\[
(e^{2\delta_\G mT}-1)^2\ge4 \quad\Longrightarrow\quad \log3\le2\delta_\G mT\le2\delta_\G C_0
\le\frac{2}{1000A}<\log3,
\]
a contradiction. Each component is therefore a tree or has rank one.

Let $K_i$, $i\in\mathcal{I}$, be the components of $G_0$ with $\operatorname{rank}\pi_1(K_i)=1$,
where $\#\mathcal{I}\le m$. Each contains a unique embedded circuit. We choose a vertex on this
circuit, a lift $\widetilde p_i\in\widetilde G$ of that vertex, and set $p_i=F(\widetilde p_i)$. We
choose an orientation of the circuit and lift one traversal beginning at $\widetilde p_i$. Its
endpoint is $a_i\widetilde p_i$ for a unique $a_i\in\G$,
\[
\widetilde p_i\longrightarrow a_i\widetilde p_i.
\]
Deleting one circuit edge and extending the remaining tree to a spanning tree of $G$ shows that
$a_i$ is primitive. Let $A_i$ be its axis, oriented by $a_i$, and set
\[
L_i=\ell(a_i)\in(0,C_0],\qquad U_i=C_0+\log^+(1/L_i),\qquad D_i=\dis(p_i,A_i).
\]
The circuit implies $\dis(p_i,a_ip_i)\le C_0$. We choose a unit-speed parametrization
$P_i:\mathbb{R}\to A_i$ such that $P_i(0)$ is the nearest point on $A_i$ to $p_i$ and
$a_iP_i(t)=P_i(t+L_i)$, and set
\[
\begin{gathered}
V_i=\operatorname{span}\{P_i(t):t\in\mathbb{R}\},\qquad R_i=a_i|_{V_i^\perp},\\
p_i=\cosh D_i\,P_i(0)+\sinh D_i\,v_i,\qquad v_i\in V_i^\perp,\quad\langle v_i,v_i\rangle=1,
\end{gathered}
\]
with any unit $v_i\in V_i^\perp$ if $D_i=0$. Since $R_i$ is orthogonal on the positive-definite
space $V_i^\perp$,
\[
\begin{aligned}
\cosh\dis(p_i,a_ip_i) &=\cosh^2D_i\cosh L_i-\sinh^2D_i\langle v_i,R_iv_i\rangle\\
&\ge1+2\sinh^2(L_i/2)\cosh^2D_i,\\
\sinh(C_0/2)&\ge\sinh(L_i/2)\cosh D_i,
\end{aligned}
\]
and $\operatorname{arcosh}y\le\log(2y)$ for $y\ge1$ gives
\begin{equation}\label{anchor}
\begin{aligned}
D_i&\le\log\frac{2\sinh(C_0/2)}{\sinh(L_i/2)}\\
&\le\frac{C_0-L_i}{2}+\log\frac{2}{1-e^{-L_i}}\\
&\le\frac{C_0}{2}+\log4+\log^+(1/L_i)\le U_i,
\end{aligned}
\end{equation}
using $1-e^{-L}\ge L/2$ for $0<L\le1$, $1-e^{-L}>1/2$ for $L\ge1$, and $C_0\ge3\cdot11^3>2\log4$.

{\bf Step 4:} Enlarging the components. We enlarge each $K_i$, keeping the components disjoint and
of rank one, hence $Z_i\le AU_i$ and every remaining edge to $K_i$ satisfies $L_e>200Z_i$. For each
$i$, in nonincreasing order of $U_i$, we repeatedly apply the edge rule below to $K_i$. We pass to
the next index when no edge satisfies the rule, and stop the process if either the second or third
case below occurs. The points $\widetilde p_i,p_i$ and the elements $a_i$ remain fixed, and for each
enlarged $K_i$ the component of $\pi_G^{-1}(K_i)$ containing $\widetilde p_i$ is used. Initially,
let $\mathcal{F}$ be the set of tree components of $G_0$. At each stage, set
\[
\mathcal{K}=\{K_i:i\in\mathcal{I}\}\cup\mathcal{F},\qquad E_{\rm rem}
=E(G)\setminus\bigcup_{K\in\mathcal{K}}E(K),
\]
and
\[
W_i=\leg(K_i),\qquad Z_i=W_i+U_i,
\]
so that the sets $V(K)$, $K\in\mathcal{K}$, partition $V(G)$. For the current index $i$, choose
$e\in E_{\rm rem}$ incident to $K_i$ with $L_e\le200Z_i$. If no such edge exists, proceed to the
next index, and otherwise orient $e:v\to w$ with $v\in V(K_i)$.
\begin{itemize}
\item If $w\in V(J)$ for $J\in\mathcal{F}$, replace the vertex and edge sets of $K_i$ by
\[
V(K_i)\cup V(J),\qquad E(K_i)\cup\{e\}\cup E(J),
\]
respectively, remove $J$ from $\mathcal{F}$, redefine $\mathcal{K},E_{\rm rem},W_i,Z_i$ as above,
and repeat for the same index $i$.

\item If $w\in V(K_i)$, stop the process with the \emph{new-circuit}
\[
H=K_i\cup e.
\]

\item If $w\in V(K_j)$ for $j\ne i$, stop the process and set
\[
H=K_j\cup e\cup K_i, \qquad\mbox{\emph{merger}}.
\]
\end{itemize}

When the edge $e$ joins $K_i$ to a tree component $J$, the graph enlarged by adjoining $e$ and $J$
deformation retracts onto the previous $K_i$ by collapsing $J\cup e$, fixing that graph pointwise,
so its fundamental group remains generated by the original circuit representing $a_i$, and the
component lift containing $\widetilde p_i$ still has stabilizer $<a_i>$. The added tree $J$ has
length at most $C_0\le U_i$, so the new value $Z_i'$ satisfies
\[
Z_i'\le Z_i+200Z_i+C_0 \le201Z_i+U_i\le202Z_i<203Z_i.
\]
Initially $Z_i\le2U_i$, and since each enlargement adds a distinct joining edge and $|E(G)|\le m$,
there are at most $m$ enlargements, so
\begin{equation}\label{growth}
Z_i\le2\cdot203^mU_i=AU_i.
\end{equation}
In the other two cases,
\[
\operatorname{rank}\pi_1(H)=
\begin{cases}
\operatorname{rank}\pi_1(K_i)+1=2,&\text{new circuit},\\
\operatorname{rank}\pi_1(K_i)+\operatorname{rank}\pi_1(K_j)=2,&\text{merger}.
\end{cases}
\]
If $K_j$ has not been enlarged,
\[
W_j\le C_0\le U_j\le AU_j,
\]
and otherwise by (\ref{growth}) we have
\[
W_j\le Z_j\le AU_j,
\]
so we have
\[
\leg(H)\le
\begin{cases}
201AU_i,&\text{new circuit},\\
201AU_i+AU_j,&\text{merger}.
\end{cases}
\]
With $U_*=U_i$ in the first case and $U_* =\max\{U_i,U_j\}$ in the second,
\[
\leg(H)\le202AU_*.
\]
We choose $k$ such that $K_k\subset H$ and $U_k=U_*$, and set
\[
(a_*,\widetilde p_*,p_*,A_*,L_*) =(a_k,\widetilde p_k,p_k,A_k,L_k), \qquad p_*=F(\widetilde p_*).
\]
Base loops in $H$ at $\pi_G(\widetilde p_*)$, a loop representing $\g\in\G$ if its lift from
$\widetilde p_*$ ends at $\g\widetilde p_*$. A spanning tree of $H$ containing all but one edge of
this starting circuit, extended to a spanning tree of $G$, gives a basis containing $a_*$ and a
second element $b\in\G$ with
\[
\dis(p_*,bp_*)\le2\leg(H)\le404AU_*.
\]
In particular $a_*,b$ freely generate a subgroup of rank two. For the point $z_*\in A_*$ nearest to
$p_*$, it follows from (\ref{anchor}) that
\[
\dis(z_*,bz_*) \le\dis(p_*,bp_*)+2\dis(p_*,z_*) \le(404A+2)U_*<1000AU_*,
\]
as $A\ge2$. On the other hand, by $0<L_*\le C_0$, $C_0\ge1$ and $\delta_\G C_0\le1/(1000A)$ we have
\[
\begin{aligned}
\log\frac{2}{\delta_\G L_*} &\ge\log\frac{2}{\delta_\G C_0}+\log^+(1/L_*)\\
&\ge\log(2000A)+\log^+(1/L_*)\\
&>1+\frac{\log^+(1/L_*)}{C_0}=\frac{U_*}{C_0}.
\end{aligned}
\]
Here the first inequality holds since for $L_*\ge1$ we have $\log^+(1/L_*)=0$ and $L_*\le C_0$, and
for $0<L_*<1$
\[
\log\frac{2}{\delta_\G L_*} =\log\frac{2}{\delta_\G C_0} +\log C_0+\log(1/L_*) \ge
\log\frac{2}{\delta_\G C_0}+\log^+(1/L_*),
\]
since $C_0\ge1$.

Since $z_*\in A_*$, it follows from (\ref{cycle}) that
\[
\dis(z_*,bz_*) \ge\frac{1}{\delta_\G}\log\frac{2}{\delta_\G L_*} >\frac{U_*}{\delta_\G
C_0}\ge1000AU_*,
\]
a contradiction, hence neither the second nor the third case occurs, and the process ends with
pairwise disjoint connected graphs $K_i$, each containing a unique circuit, satisfying
\[
Z_i\le AU_i,\qquad L_e>200Z_i \quad\text{for every remaining edge to }K_i.
\]
The graphs remaining in $\mathcal{F}$ are the tree components of $G_0$ that were not added to any
$K_i$.

{\bf Step 5:} The outward ideal balls of the leaving edges have small variation of $\log|\zeta|$: We
fix an enlarged component and set
\[
K=K_i,\quad a=a_i,\quad A_a=A_i,\quad L=L_i,\quad U=C_0+\log^+(1/L),\quad W=\leg(K),\quad Z=W+U,
\]
and let $\widetilde K$ be its lift with stabilizer $<a>$. We will show that, for every edge leaving
$\widetilde K$, the function $\zeta\mapsto\log|\zeta|$ varies by less than $L/(4m)$ on the outward
ideal ball of its perpendicular midpoint plane.

Deleting one circuit edge gives a spanning tree, whose lift containing $\widetilde p_i$ gives vertex
representatives $\widetilde v$ satisfying
\[
\dis(F\widetilde v,A_a)\le W+D_i\le Z, \qquad \dis(F(a^k\widetilde v),A_a) =\dis(F\widetilde
v,A_a)\le Z \quad(k\in\mathbb{Z}).
\]
Let $\widetilde e:\widetilde p\to\widetilde q$ leave $\widetilde K$, and set
\[
p=F\widetilde p,\qquad q=F\widetilde q,\qquad 2R=\dis(p,q)>200Z.
\]
For $k\in\mathbb{Z}$, let $\sigma_k$ be the projection to $K$ of a path in $\widetilde K$ from
$a^k\widetilde p$ to $\widetilde p$, which avoids the projection $e$ of $\widetilde e$ to $G$, since
that edge lies outside $K$. Replacing $e$, through the mutually inverse mappings in Step~1, by an
edge represented by the path $\sigma_ke$, whose lift joins $a^k\widetilde p$ to $\widetilde q$, and
straightening it to the segment $[a^kp,q]$ preserves the induced isomorphism on fundamental groups
and leaves all retained edge lengths unchanged, so by minimality of $G$ we have
\[
\leg(G) \le\leg(G)-\dis(p,q)+\dis(a^kp,q),
\]
that is,
\begin{equation}\label{winding}
\dis(a^kp,q)\ge\dis(p,q)=2R \qquad(k\in\mathbb{Z}).
\end{equation}

We choose Lorentz coordinates such that
\[
\begin{gathered}
P(t)=(\cosh t,\sinh t,0),\qquad aP(t)=P(t+L),\\
V_a=\operatorname{span}\{P(t):t\in\mathbb{R}\},\qquad \|v\|_a^2=\langle v,v\rangle\quad(v\in
V_a^\perp),\\ D_p=\dis(p,A_a),\qquad D_q=\dis(q,A_a),\\ p=\cosh D_p\,P(0)+\sinh D_p\,v_p,\qquad
q=\cosh D_q\,P(z)+\sinh D_q\,v_q,\\ v_p,v_q\in V_a^\perp,\qquad \|v_p\|_a=\|v_q\|_a=1, \qquad
R_a=a|_{V_a^\perp},
\end{gathered}
\]
with any unit vectors $v_p,v_q\in V_a^\perp$ if $D_p=0$ or $D_q=0$. For $k\in\mathbb{Z}$,
\[
\cosh\dis(a^kp,q) =\cosh D_p\cosh D_q\cosh(kL-z) -\sinh D_p\sinh D_q\langle R_a^kv_p,v_q\rangle,
\]
and for $k$ with $|kL-z|\le L/2$, by (\ref{winding}),
\[
\begin{aligned}
&\cosh D_p\cosh D_q \bigl(\cosh z-\cosh(kL-z)\bigr)\\
&\qquad\le\sinh D_p\sinh D_q \bigl(\langle v_p,v_q\rangle -\langle R_a^kv_p,v_q\rangle\bigr)
\le2\sinh D_p\sinh D_q,\\ &\cosh z\le\cosh(L/2)+2\tanh D_p\tanh D_q \le\cosh(L/2)+2.
\end{aligned}
\]
For $t\ge0$,
\[
\cosh(t+\log6)-\cosh t =\frac52e^t-\frac5{12}e^{-t}\ge\frac{25}{12}>2.
\]
Hence
\begin{equation}\label{phase}
|z|\le M:=\operatorname{arcosh}(\cosh(L/2)+2) \le L/2+\log6\le C_0\le U,
\end{equation}
where the last bound uses $L\le C_0$ and $C_0/2\ge\log6$, and by the path through the axis feet we
have
\[
2R\le D_p+|z|+D_q,\qquad D_p\le Z,\qquad D_q\ge2R-Z-M.
\]

Let $\mathcal{P}$ be the perpendicular midpoint plane of $[p,q]$. Its outward normal $n$, directed
towards $q$, satisfies
\[
\begin{gathered}
n=\frac{q-p}{2\sinh R},\qquad \langle q-p,q-p\rangle=2(\cosh2R-1)=4\sinh^2R,\\
\langle n,n\rangle=1,\qquad n_A=\alpha P(z)-\beta P(0),\qquad \alpha=\frac{\cosh D_q}{2\sinh
R},\quad \beta=\frac{\cosh D_p}{2\sinh R}.
\end{gathered}
\]
Here $n_A$ is the Lorentz projection onto $V_a$. Since $R>100Z$ and $M\le U\le Z$,
\[
\frac{\alpha}{\beta} \ge\frac12e^{D_q-D_p} \ge\frac12e^{2R-2Z-M},\qquad \frac{\alpha}{\beta}e^{-M}
>\frac12e^{196Z}>2.
\]
Hence we have
\[
\begin{gathered}
n_A^0\pm n_A^1=\alpha e^{\pm z}-\beta>0,\qquad 0<\frac\beta\alpha e^{\pm z}<\frac12,\\
s_A:=\sqrt{-\langle n_A,n_A\rangle}>0,\qquad s_A^2=-\langle n_A,n_A\rangle\\
=(\alpha-\beta e^z)(\alpha-\beta e^{-z})>\frac{\alpha^2}{4},\\
n_A=s_AP(\tau),\qquad \tau=\frac12\log\frac{\alpha e^z-\beta}{\alpha e^{-z}-\beta},\\
\langle P(t),n\rangle=-s_A\cosh(t-\tau)<0.
\end{gathered}
\]
In particular, $n_A$ is future timelike and the axis lies outside the outward halfspace. By the
point-to-plane distance formula we have
\[
D:=\dis(A_a,\mathcal{P}) =\min_t\operatorname{arsinh}\bigl(s_A\cosh(t-\tau)\bigr)
=\operatorname{arsinh}s_A.
\]
Moreover,
\[
\alpha\ge\frac12e^{D_q-R},\qquad \operatorname{arsinh}(\alpha/2)
=\log\!\left(\frac\alpha2+\sqrt{1+\frac{\alpha^2}{4}}\right) \ge\log\alpha.
\]
Hence,
\begin{equation}\label{cut-axis}
\begin{aligned}
D&>\operatorname{arsinh}(\alpha/2)\ge\log\alpha\\
&\ge D_q-R-\log2\ge R-Z-M-\log2\\
&>99Z-M-\log2\ge98U-\log2>98U-2.
\end{aligned}
\end{equation}
Also
\[
\tau-z=\frac12\log \frac{1-(\beta/\alpha)e^{-z}}{1-(\beta/\alpha)e^z},\qquad
|\tau-z|<\frac12\log2<1.
\]

In upper-half-space coordinates with axis $(0,\infty)$ and $P(0)$ at height one, oriented towards
infinity, the ideal boundary point $\zeta\in\mathbb{R}^{N-1}\setminus\{0\}$ is represented by the
vector
\[
\iota(\zeta) =\left(\frac{|\zeta|^2+1}{2}, \frac{|\zeta|^2-1}{2},\zeta\right)
=|\zeta|\left[P(\log|\zeta|) +\left(0,0,\frac{\zeta}{|\zeta|}\right)\right].
\]
Set $n_\perp=n-n_A$, so that
\[
\mathcal{B}=B_n=\partial_\infty E_n,\qquad \|n_\perp\|_a^2=1+s_A^2=\cosh^2D,
\]
and for $\zeta\in\mathcal{B}\setminus\{0,\infty\}$, with
\[
t=\log|\zeta|,\qquad v=\left(0,0,\frac{\zeta}{|\zeta|}\right),
\]
we have $\|v\|_a=1$ and $\iota(\zeta)=|\zeta|(P(t)+v)$, so that
\[
\begin{aligned}
0 &\le\frac{\langle\iota(\zeta),n\rangle}{|\zeta|}=\langle P(t)+v,n\rangle\\
&\le -\sinh D\cosh(t-\tau)+\langle v,n_\perp\rangle\\
&\le-\sinh D\cosh(t-\tau) +\|v\|_a\|n_\perp\|_a\\
&\le-\sinh D\cosh(t-\tau)+\cosh D.
\end{aligned}
\]
The two axis endpoints are excluded since
\[
\langle(1,\pm1,0),n\rangle=-s_Ae^{\mp\tau}<0.
\]
Hence we have
\[
\begin{gathered}
\{\log|\zeta|:\zeta\in\mathcal{B}\} \subset I_{\mathcal{B}}:=[\tau-\omega(D),\tau+\omega(D)],\\
\omega(D)=\operatorname{arcosh}\coth D =\log\coth(D/2) =\log\frac{1+e^{-D}}{1-e^{-D}}.
\end{gathered}
\]
For $D\ge\log2$,
\[
\omega(D)=2\int_0^{e^{-D}}\frac{du}{1-u^2} \le\frac83e^{-D}\le4e^{-D}.
\]
Since $U\ge C_0\ge m\ge3$,
\[
\log(32m)+2-98C_0 <m+5-98m=5-97m<0,
\]
and
\[
\frac{32m e^{2-98U}}{L} =32m e^{2-98C_0}
\begin{cases}
L^{-1},&L\ge1,\\
L^{97},&0<L<1
\end{cases}
<1,
\]
and together with (\ref{cut-axis}) we get
\[\omega(D)\le4e^{-D}<4e^{2-98U}<\frac{L}{8m}.
\]

{\bf Step 6:} Disjointness at a cyclic component. The ideal balls determined by the midpoint planes
of all edges leaving $\widetilde K$ are pairwise disjoint, and one ball from each $\langle
a\rangle$-orbit can be placed inside an annulus $e^{s_0}<|\zeta|<e^{s_0+L}$. Consider the lifts of
oriented remaining edges whose initial vertices lie in $\widetilde K$, both orientations included if
both endpoints of an edge in $G$ lie in $K$, and choose one such lift from each $\langle
a\rangle$-orbit. Since $a$ multiplies $|\zeta|$ by $e^L$, with $I_{a^k\mathcal
B}:=[\tau+kL-\omega(D),\,\tau+kL+\omega(D)]=I_{\mathcal{B}}+kL$,
\[
|I_{\mathcal{B}}|=2\omega(D)<L,\qquad \mathcal{B}\cap a^k\mathcal{B}=\emptyset \quad(k\ne0).
\]
Suppose the ideal balls associated with two lifted oriented edges leaving $\widetilde K$ from
distinct $<a>$-orbits intersect. Let $\widetilde p_j,\widetilde q_j$ be their endpoint vertices and
$\widetilde x_j$ their midpoint vertices after subdivision, and set
\[
p_j=F\widetilde p_j,\qquad q_j=F\widetilde q_j,\qquad x_j=F\widetilde x_j,\qquad
2R_j=\dis(p_j,q_j)\quad(j=1,2).
\]
With the same parametrization $P:\mathbb{R}\to A_a$, define $t_j,\tau_j$ by
\[
\dis(p_j,P(t_j))=\dis(p_j,A_a),\qquad \dis(P(\tau_j),\mathcal{P}_j)=\dis(A_a,\mathcal{P}_j)
\quad(j=1,2),
\]
and set $D_j=\dis(A_a,\mathcal{P}_j)$ and $\omega_j=\omega(D_j)$. The calculations of Step 5 show
that
\[
\begin{gathered}
|\tau_j-t_j|<M+1,\qquad \dis(p_j,A_a)\le Z,\\
\omega_1+\omega_2<8e^{2-98U}<1,\qquad |\tau_1-\tau_2|\le\omega_1+\omega_2<1,\\
|t_1-t_2|<2M+3,\qquad \dis(p_1,p_2)<2Z+2M+3\le4Z+3.
\end{gathered}
\]
By (\ref{overlap}),
\[\dis(x_1,x_2)<|R_1-R_2|+4Z+5<\max\{R_1,R_2\},
\]
since $\min\{R_1,R_2\}>100Z>4Z+5$.

For distinct quotient edges, order $R_1\le R_2$ and subdivide the first edge at its midpoint. The
path
\[
\widetilde x_1\longrightarrow\widetilde p_1 \longrightarrow\widetilde p_2
\]
uses the first half-edge in reverse and a component path, so it avoids the second quotient edge.
Replacing the second edge by an edge whose chosen lift joins $\widetilde x_1$ to $\widetilde q_2$,
mapped to the geodesic segment $[x_1,q_2]$, we get
\[
\begin{aligned}
\leg(G_{\rm new})-\leg(G) &=\dis(x_1,q_2)-2R_2\\
&\le\dis(x_1,x_2)-R_2<0,
\end{aligned}
\]
and after collapsing any edges of length zero this is a carrier shorter than $G$, contradicting the
minimality of $G$, hence the two ideal balls are disjoint in this case.

For opposite orientations of one quotient edge $e$, $R_1=R_2=R$. Choose
\[
\widetilde e:\widetilde p\longrightarrow b_e\widetilde q, \qquad \widetilde p,\widetilde
q\in\widetilde K, \qquad x_0=F(\text{midpoint of }\widetilde e).
\]
The lift of the reversed edge starting at $\widetilde q$ ends at $b_e^{-1}\widetilde p$, and by
equivariance of $F$ its midpoint maps to $b_e^{-1}x_0$. Hence, for some $j_1,j_2\in\mathbb{Z}$,
\[
\begin{gathered}
x_1=a^{j_1}x_0,\qquad x_2=a^{j_2}b_e^{-1}x_0,\\
\g=a^{j_2}b_e^{-1}a^{-j_1},\qquad \g x_1=x_2,\\
\g_0=a^{-j_1}\g a^{j_1}=a^{j_2-j_1}b_e^{-1},\qquad \dis(x_0,\g_0x_0)=\dis(x_1,x_2)<R.
\end{gathered}
\]
The endpoints of $e$ are joined in $K$, which implies that $G-e$ is connected, and at the initial
vertex $o_e$ of $e$
\[
H_e=\pi_1(G-e,o_e),\qquad \pi_1(G,o_e)=H_e*<t_e>,
\]
for a spanning tree of $G-e$ and the corresponding oriented generator $t_e$ for $e$, the chosen lift
at $o_e$ identifying this group with $\G$. A component path joins $\widetilde q$ to $\widetilde p$,
so the loop that traverses $e$ and then this path represents $b_e$, and hence
\[
a\in H_e,\qquad b_e=h_1t_eh_2\quad(h_1,h_2\in H_e), \qquad \g_0=
\underbrace{a^{j_2-j_1}h_2^{-1}}_{u}\, t_e^{-1}\underbrace{h_1^{-1}}_{v}.
\]
We subdivide $e$ at its midpoint, retain its first half, delete its second half, and attach a loop
at the midpoint vertex whose lifted image is the geodesic segment $[x_0,\g_0x_0]$. The based loop
has lifted image
\[
F\widetilde p\longrightarrow x_0 \longrightarrow\g_0x_0 \longrightarrow\g_0F\widetilde p,
\]
so its class is $ut_e^{-1}v$. The homomorphisms $\phi,\psi$ defined by
\[
\phi|_{H_e}=\psi|_{H_e}=\mathrm{id},\qquad \phi(t_e)=ut_e^{-1}v,\qquad \psi(t_e)=vt_e^{-1}u,
\]
satisfy
\[
\phi\psi(t_e)=v(ut_e^{-1}v)^{-1}u=t_e,\qquad \psi\phi(t_e)=u(vt_e^{-1}u)^{-1}v=t_e,
\]
so the new carrier induces an isomorphism. Its added loop is nontrivial, and
\[
\leg(G_{\rm new})-\leg(G) =R+\dis(x_0,\g_0x_0)-2R<0,
\]
which again contradicts minimality, hence the ideal balls determined by the midpoint planes of all
lifted oriented remaining edges leaving $\widetilde K$ are pairwise disjoint.

Let $q_K$ be the number of $<a>$-orbits of lifted oriented remaining edges leaving $\widetilde K$.
Since the component $\widetilde K$ is the universal cover of $K$, its stabilizer $<a>$ acts
transitively on the lifts of each vertex of $K$, and unique lifting of an oriented edge from a
chosen initial vertex gives one orbit for each oriented remaining-edge incidence at $K$, of which
there are at most $2|E(G)|\le2m$ and, as $G$ is connected and $\operatorname{rank}\pi_1(K)=1<r$, at
least one, which implies $1\le q_K\le2m$. For their chosen intervals $I_j$,
\[
\sum_{j=1}^{q_K}|I_j|<\frac{q_KL}{4m}\le\frac L2,
\]
so with $\mathbb{T}_L=\mathbb{R}/L\mathbb{Z}$, $p_L:\mathbb{R}\to\mathbb{T}_L$ the quotient mapping
and $\mu_L$ arclength measure on this circle,
\[
\mu_L\!\left(\bigcup_{j=1}^{q_K}p_L(I_j)\right) \le\sum_{j=1}^{q_K}|I_j|<L/2<L=\mu_L(\mathbb{T}_L).
\]
We choose $s_0\in\mathbb{R}$ such that $p_L(s_0)$ avoids this union. Each connected interval $I_j$
avoids $s_0+L\mathbb{Z}$, so it lies in one component of $\mathbb{R}\setminus(s_0+L\mathbb{Z})$, and
there is a unique $k_j\in\mathbb{Z}$ such that
\[
I_j+k_jL\subset(s_0,s_0+L).
\]
Let $n_j$ be the outward unit normal of the image under $a^{k_j}$ of the perpendicular midpoint
plane corresponding to $I_j$, and set
\[
\begin{gathered}
n_-=-P'(s_0),\qquad n_+=P'(s_0+L),\qquad an_-=-n_+,\\
C=\{X\in\mathbb{H}^N: \langle X,n_-\rangle\le0, \ \langle X,n_+\rangle\le0, \ \langle
X,n_j\rangle\le0\ (1\le j\le q_K)\}.
\end{gathered}
\]
The ideal balls defined by $n_-$ and $n_+$ are, respectively,
\[
\{|\zeta|\le e^{s_0}\},\qquad \{|\zeta|\ge e^{s_0+L}\}\cup\{\infty\},
\]
and the ideal balls defined by $n_j$, $1\le j\le q_K$, lie in the open annulus
\[
e^{s_0}<|\zeta|<e^{s_0+L}.
\]
For $s_0<t<s_0+L$,
\[
\langle P(t),n_-\rangle=-\sinh(t-s_0)<0,\qquad \langle P(t),n_+\rangle=\sinh(t-s_0-L)<0,\qquad
\langle P(t),n_j\rangle<0,
\]
so $P(t)\in\inte C$, and all closed outward ideal balls of $C$ are disjoint.

{\bf Step 7:} Tree components. Each remaining tree component $K$ also has a region $C$ with
$F(\widetilde K)\subset\inte C$ and pairwise disjoint ideal balls. Let $K$ be a tree component of
$G_0$ that has not been added to any $K_i$. Since $K$ is a tree,
\[
\operatorname{Stab}_\G(\widetilde K) \cong\pi_1(K)=\{1\}.
\]
Its edges have length at most $T$, so
\[
\leg(K) \le |E(K)|T\le mT.
\]
Note that every short edge incident to $K$ belongs to $K$, since $K$ is a component of $G_0$, so an
edge leaving $\widetilde K$ projects to a long edge, of length $2R\ge BT$, and since $B=2m+5$ and
$T\ge1$,
\[
R\ge\frac{BT}{2} =mT+\frac52T \ge mT+\frac52>mT+2.
\]
Let $\widetilde e_j:\widetilde p_j\to\widetilde q_j$, $j=1,2$, be distinct edges leaving $\widetilde
K$ whose outgoing closed ideal balls intersect. Set
\[
p_j=F\widetilde p_j,\qquad q_j=F\widetilde q_j,\qquad x_j=F(\text{midpoint of }\widetilde
e_j),\qquad 2R_j=\dis(p_j,q_j).
\]
Since $\dis(p_1,p_2)\le mT$ and $R_j>mT+2$, it follows from (\ref{overlap}) that
\[
\dis(x_1,x_2)<|R_1-R_2|+mT+2<\max\{R_1,R_2\}.
\]
If the lifted edges project to distinct edges of $G$, assume $R_1\le R_2$, subdivide the first edge
at its midpoint and replace the second edge by one whose lifted image is $[x_1,q_2]$. As above, this
preserves the carrier property, while
\[
\begin{aligned}
\leg(G_{\rm new})-\leg(G) &=\dis(x_1,q_2)-2R_2\\
&\le\dis(x_1,x_2)-R_2<0,
\end{aligned}
\]
contradicting the minimality of $G$. If instead the two lifted edges project to opposite
orientations of one edge $e:v\to w$, then $v,w\in V(K)$, since both lifted edges start in
$\widetilde K$, and since $K$ is a tree, the restriction $\pi_G|_{\widetilde K}:\widetilde
K\longrightarrow K$ is an isomorphism, so there are unique vertices $\widetilde p,\widetilde q$ of
$\widetilde K$ lying over $v,w$, respectively. Denote the first lifted edge by $\widetilde e$, and
let $b_e\in\G$ be the deck transformation, unique since the deck group acts freely and transitively
on the vertices lying over $w$, for which the terminal vertex of $\widetilde e$ is $b_e\widetilde
q$, hence $\widetilde e:\widetilde p\longrightarrow b_e\widetilde q$ and the other edge is
$b_e^{-1}\overline{\widetilde e}:\widetilde q\longrightarrow b_e^{-1}\widetilde p$. With
$x_0=F(\text{midpoint of }\widetilde e)$ and $2R=L_e$, the midpoint images of the two edges are
$x_0$ and $b_e^{-1}x_0$, so
\[
\dis(x_0,b_e^{-1}x_0)<mT+2<R.
\]
The endpoints of $e$ are joined in $K$, which implies that $G-e$ is connected, and for a graph basis
based at $v$ and built from a spanning tree of $G-e$ containing $K$, the generator associated with
$e$ is represented by $e[w,v]_K$, whose lift ends at $b_e\widetilde p$, so this generator maps to
$b_e$. Subdividing $e$ at its midpoint, retaining its first half, deleting its second half, and
attaching a loop at the midpoint whose lifted image is $[x_0,b_e^{-1}x_0]$ gives a based loop with
lifted image
\[
F\widetilde p\longrightarrow x_0 \longrightarrow b_e^{-1}x_0 \longrightarrow b_e^{-1}F\widetilde p,
\]
so the new mapping replaces $b_e$ by $b_e^{-1}$ and fixes every other basis generator, a
substitution which is its own inverse, and the new graph remains a carrier, while
\[
\leg(G_{\rm new}) -\leg(G) =R+\dis(x_0,b_e^{-1}x_0)-2R<0,
\]
contradicting minimality. Since each oriented edge starting in $K$ has one lift starting in
$\widetilde K$, the two cases prove that all outgoing closed ideal balls are pairwise disjoint.

Let $\widetilde e:\widetilde p\to\widetilde q$ leave $\widetilde K$, and set
\[
u=F(\widetilde p),\qquad 2R=\dis(u,F(\widetilde q)).
\]
With $\nu$ the unit tangent at $u$ towards $F(\widetilde q)$, the unit normal to the midpoint plane,
directed towards $F(\widetilde q)$, is
\[
n=\sinh R\,u+\cosh R\,\nu.
\]
For $w\in F(\widetilde K)$, set $d=\dis(u,w)\le mT$ and choose $\eta\in T_u\mathbb{H}^N$ with
$w=\cosh d\,u+\sinh d\,\eta$ and $\langle\eta,\eta\rangle=1$, so that
\[
\begin{aligned}
\langle w,n\rangle &=-\cosh d\sinh R+\sinh d\cosh R\langle\eta,\nu\rangle\\
&\le-\sinh(R-d)\le-\sinh(R-mT)<0.
\end{aligned}
\]
With
\[
C=\bigcap_e\{X\in\mathbb{H}^N:\langle X,n_e\rangle\le0\},
\]
where $e$ ranges over the lifted oriented remaining edges leaving $\widetilde K$, we get
\[
F(\widetilde K)\subset\inte C.
\]
The closed ideal balls
\[
\partial_\infty\{X\in\mathbb{H}^N:\langle X,n_e\rangle\ge0\}
\]
are pairwise disjoint.

{\bf Step 8:} The tree $\mathcal{T}'$. This step constructs a tree $\mathcal{T}'$ with a free
$\G$-action and quotient graph of fundamental group rank $r$, with regions assigned so that adjacent
regions have a common boundary plane with opposite halfspaces. We collapse each enlarged component
$K_i$ and each tree component of $G_0$ not added to any $K_i$ to a vertex, obtaining a connected
graph $Y$, and let $K_v$ be the component collapsed to $v\in V(Y)$. Set
\[
n=\#V(Y),\qquad c=\#E(Y),\qquad s=\#\{v\in V(Y):K_v\text{ is cyclic}\}.
\]
Since $\#E(K_v)-\#V(K_v)$ is zero for a cyclic component and $-1$ for a tree,
\[
r=\#E(G)-\#V(G)+1 =c+\sum_v\bigl(\#E(K_v)-\#V(K_v)\bigr)+1 =c+s-n+1.
\]
Each component lift is a connected subtree of $\widetilde G$, and $q_0:\widetilde G\to\mathcal{T}$
denotes the map collapsing these subtrees. For an edge $\widetilde f$ outside them, let $\mathcal
U_1,\mathcal{U}_2$ be the two components, so
\[
\widetilde G\setminus\inte\widetilde f =\mathcal{U}_1\sqcup\mathcal{U}_2.
\]
Each component lift lies in one $\mathcal{U}_j$, so, for $f=q_0(\widetilde f)$,
\[
\mathcal{T}\setminus\inte f =q_0(\mathcal{U}_1)\sqcup q_0(\mathcal{U}_2).
\]
Hence every edge of the connected graph $\mathcal{T}$ is a bridge, which implies that $\mathcal{T}$
is a tree. For each $v\in V(Y)$, apply Step 6 or Step 7 to $K_v$, and denote the chosen component of
$\pi_G^{-1}(K_v)$ by $\widetilde K_v^0$ and the corresponding region in $\mathbb{H}^N$ by $C_v^0$.
We choose the component lift $\widetilde K_v$ and $\eta_v\in\G$ such that
\[
\widetilde K_v=\eta_v\widetilde K_v^0,\qquad C_v=\eta_vC_v^0.
\]
At a cyclic vertex, let $a_v^0,A_v^0$ be the generator and axis, and set
\[
a_v=\eta_va_v^0\eta_v^{-1},\qquad A_v=\eta_vA_v^0,
\]
and
\[
E_v^\pm =\{\eta_vX:X\in\mathbb{H}^N,\ \langle X,n_\pm\rangle\ge0\},
\]
where $n_-,n_+$ are the normals used to define $C_v^0$, and if $K_v$ is a tree, set $a_v=1$. Hence
\[
H_v=\operatorname{Stab}_\G(\widetilde K_v) =\begin{cases} \langle a_v\rangle,&K_v\text{ cyclic},\\
\{1\},&K_v\text{ a tree}.
\end{cases}
\]
At a cyclic vertex the lower and upper outgoing halfspaces satisfy
\[
a_vE_v^-=\mathbb{H}^N\setminus\inte E_v^+,
\]
so that $a_vP_v^-=P_v^+$ for $P_v^\pm=\partial E_v^\pm$. Since the closed ideal balls of the outward
halfspaces of $C_v$ are pairwise disjoint, so are these closed halfspaces, as in Proposition
\ref{gap}, and no other outward halfspace removes a point of $P_v^+$ or $P_v^-$ from $C_v$, so
\[
P_v^+\subset C_v,\qquad P_v^+=a_vP_v^-\subset a_vC_v,
\]
while
\[
C_v\subset\mathbb{H}^N\setminus\inte E_v^+, \qquad a_vC_v\subset E_v^+,
\]
and therefore
\[
C_v\cap a_vC_v =E_v^+\cap a_vE_v^- =P_v^+, \qquad E_v^+\cup a_vE_v^-=\mathbb{H}^N.
\]

We orient each edge of $Y$ and write $e:v\to w$, where $v$ and $w$ are its initial and terminal
vertices, and a lift of the corresponding edge of $G$ of the form
\[
\widetilde e:\widetilde p_e\longrightarrow h_e\widetilde q_e, \qquad h_e\in\G,\qquad \widetilde
p_e\in\widetilde K_v,\quad \widetilde q_e\in\widetilde K_w.
\]
Its forward and reverse midpoint halfspaces satisfy
\[
h_eE_{e,-}=\mathbb{H}^N\setminus\inte E_{e,+}.
\]
Let $k_e,l_e\in\mathbb{Z}$ so that the selected outgoing halfspaces of $C_v,C_w$ are
\[
E'_{e,+}=a_v^{k_e}E_{e,+},\qquad E'_{e,-}=a_w^{l_e}E_{e,-},
\]
with $k_e=0$ if $K_v$ is a tree and $l_e=0$ if $K_w$ is a tree, and define
\[
h'_e=a_v^{k_e}h_ea_w^{-l_e}, \qquad h_e=a_v^{-k_e}h'_ea_w^{l_e}.
\]
For the reverse orientation $\bar e$, use
\[
h'_{\bar e}=a_w^{l_e}h_e^{-1}a_v^{-k_e}=(h'_e)^{-1}.
\]
The selected midpoint planes match, since
\[
\begin{aligned}
h'_eE'_{e,-} &=a_v^{k_e}h_ea_w^{-l_e}a_w^{l_e}E_{e,-}\\
&=a_v^{k_e}\bigl(\mathbb{H}^N\setminus\inte E_{e,+}\bigr) =\mathbb{H}^N\setminus\inte E'_{e,+}.
\end{aligned}
\]

Let $\mathcal{E}(Y)$ be the set of orientations chosen above, one for each unoriented edge of $Y$,
and let $\mathcal{T}'$ be the graph with vertices $(\g,v)$, where $\g\in\G$ and $v\in V(Y)$, and
with edges
\[
\begin{array}{lll}
\varepsilon_{\g,v}:& (\g,v)\longrightarrow(\g a_v,v), &v\text{ cyclic},\\[2pt]
\eta_{\g,e}:& (\g,v)\longrightarrow(\g h'_e,w), &e\in\mathcal{E}(Y),\quad e:v\to w,
\end{array}
\]
Here $\g$ ranges over $\G$, and where reversing a displayed edge interchanges its initial and
terminal vertices without adding a new unoriented edge, so that
\[
\overline{\varepsilon_{\g,v}}: (\g a_v,v)\longrightarrow(\g,v)
\]
and
\[
\overline{\eta_{\g,e}}: (\g h'_e,w)\longrightarrow(\g,v), \qquad h'_{\bar e}=(h'_e)^{-1}.
\]
We define a map $\pi_{\mathcal{T}}:\mathcal{T}'\to\mathcal{T}$ as follows. Each edge $(\g,v)\to(\g
a_v,v)$ is mapped to the vertex of $\mathcal{T}$ obtained by collapsing $\g\widetilde K_v$. Let
$t_v\in V(\mathcal{T})$ be the vertex obtained by collapsing $\widetilde K_v$. Since $a_v\widetilde
K_v=\widetilde K_v$, we have $a_vt_v=t_v$. Hence
\[
\pi_{\mathcal{T}}(\g,v)=\g t_v,\qquad \pi_{\mathcal{T}}(\g a_v,v) =\g a_vt_v=\g t_v.
\]
We therefore map the entire edge $(\g,v)\to(\g a_v,v)$ to the single vertex $\g t_v$. For each
chosen $e:v\to w$, the edge
\[
(\g,v)\longrightarrow(\g h'_e,w)
\]
is mapped to
\[
q_0\bigl(\g a_v^{k_e}\widetilde e\bigr).
\]
The identities
\[
a_v^{k_e}\widetilde K_v=\widetilde K_v,\qquad a_v^{k_e}h_e\widetilde K_w=h'_e\widetilde K_w
\]
show that the edge images have the required endpoint images, which implies that $\pi_{\mathcal{T}}$
is well defined. Since the deck action on $\widetilde G$ is free on points,
\[
\operatorname{Stab}_\G(\widetilde e)=\{1\},\qquad \g_1a_v^{k_e}\widetilde e =\g_2a_v^{k_e}\widetilde
e \quad\Longrightarrow\quad\g_1=\g_2.
\]
These translates are all lifts of $e$, hence each edge of $\mathcal{T}$ is the image of a unique
edge of $\mathcal{T}'$. The preimages have vertex sets
\[
V\!\left(\pi_{\mathcal{T}}^{-1}(\g t_v)\right)=
\begin{cases}
\{(\g a_v^j,v):j\in\mathbb{Z}\}, &K_v\text{ is cyclic}\\
\{(\g,v)\},&K_v\text{ is a tree},
\end{cases}
\]
with the edges forming a line in the cyclic case. For two vertices of $\mathcal{T}'$ with different
images, we follow the path between their images in $\mathcal{T}$ using the unique edges of $\mathcal
T'$ mapped onto the edges of that path. The endpoints of successive lifted edges lie in the same
preimage, where they are joined by edges or coincide when the preimage is a point. Joining also the
two original vertices to the first and last edge endpoints gives a path in $\mathcal{T}'$. If both
images agree, the path lies in their common preimage, hence $\mathcal{T}'$ is connected. For an edge
$f$, let $\mathcal{T}_1,\mathcal{T}_2$ be the two components of $\mathcal{T}\setminus\inte f$, with
the endpoints retained. Let $f'$ be the unique edge of $\mathcal{T}'$ mapped onto $f$. Since
\[
\pi_{\mathcal{T}}^{-1}(\inte f)=\inte f',
\]
we have
\[
\mathcal{T}'\setminus\inte f' =\pi_{\mathcal{T}}^{-1}(\mathcal{T}_1)
\sqcup\pi_{\mathcal{T}}^{-1}(\mathcal{T}_2).
\]
Hence every edge mapped onto an edge of $\mathcal{T}$ is a bridge. Any circuit would therefore lie
in the preimage of a vertex, which is a line or a point, hence $\mathcal{T}'$ is a tree.

For $g\in\G$, define
\[
\begin{aligned}
g(\g,v)&=(g\g,v),\\
g\varepsilon_{\g,v}&=\varepsilon_{g\g,v},\\
g\eta_{\g,e}&=\eta_{g\g,e},
\end{aligned}
\]
and for every oriented edge $f$ of $\mathcal{T}'$ set
\[
g\bar f=\overline{gf}.
\]
These formulas preserve initial and terminal vertices, and the action is free on vertices. On a
cyclic preimage we have
\[
(\g a_v^q\g^{-1})(\g a_v^j,v) =(\g a_v^{j+q},v),\qquad q,j\in\mathbb{Z}.
\]
An edge inversion would have square fixing a vertex, so torsion-freeness excludes inversions. Set
$X=\mathcal{T}'/\G$ and let $\pi_X$ be the quotient mapping. For a vertex $\tau$ of either graph,
let $\operatorname{St}(\tau)$ be the star of oriented edge incidences with initial vertex $\tau$.
For $\tau\in V(\mathcal{T}')$, every oriented edge of $X$ starting at $\pi_X(\tau)$ has a lift
starting at $\tau$, since any lift can be translated by an element of $\G$ so that its initial
vertex becomes $\tau$. If two oriented edges starting at $\tau$ have the same image in $X$, some
$g\in\G$ maps one to the other, which implies $g\tau=\tau$, $g=1$ and the edges coincide. Hence
\[
\operatorname{St}(\tau)\xrightarrow{\ \cong\ } \operatorname{St}(\pi_X\tau)
\]
is a bijection, which implies that $\pi_X$ is a graph covering. Since $\mathcal{T}'$ is a tree, this
covering is universal. A deck transformation agrees with some $g\in\G$ at one vertex and then along
every edge path, so
\[
\pi_X:\mathcal{T}'\longrightarrow X=\mathcal{T}'/\G, \qquad \operatorname{Deck}(\pi_X)=\G.
\]
The quotient graph $X$ is obtained from $Y$ by adding one loop at each vertex $v$ for which $K_v$ is
cyclic, so
\[
\#V(X)=n,\qquad \#E(X)=c+s,\qquad \operatorname{rank}\pi_1(X)=c+s-n+1=r.
\]

For $(\g,v)\in V(\mathcal{T}')$, define $C_{(\g,v)}=\g C_v$. For $g\in\G$, we have
\[
C_{g(\g,v)}=C_{(g\g,v)} =g\g C_v=gC_{(\g,v)},
\]
so the assignment is $\G$-equivariant.

{\bf Step 9:} The classical domain. In this step the regions are joined into $2r$ pairwise disjoint
closed round balls paired by a free basis of $\G$, using that disjoint closed ideal balls determine
disjoint closed hyperbolic halfspaces, as shown at the beginning of the proof of Proposition
\ref{gap}.

For adjacent vertices $\tau,\upsilon$, let $E_{\tau\upsilon}$ be the outgoing closed halfspace from
$C_\tau$ towards $C_\upsilon$, and set $B_{\tau\upsilon}=\partial_\infty E_{\tau\upsilon}$. By face
matching we have
\[
E_{\upsilon\tau}=\mathbb{H}^N\setminus\inte E_{\tau\upsilon},\qquad
B_{\upsilon\tau}=\mathbb{S}^d\setminus\inte B_{\tau\upsilon},\qquad C_\upsilon\subset
E_{\tau\upsilon}.
\]
For vertices $\omega,\upsilon\in V(\mathcal{T}')$, denote $\omega\sim\upsilon$ when they are joined
by an edge. If $\omega\sim\upsilon$ and $\omega\ne\tau$, the balls $B_{\upsilon\omega}$ and
$B_{\upsilon\tau}$ are disjoint, so their corresponding closed halfspaces are disjoint. By the
pairing identities we have
\[
\begin{aligned}
E_{\upsilon\omega} &\subset\mathbb{H}^N\setminus E_{\upsilon\tau} =\inte E_{\tau\upsilon},\\
B_{\upsilon\omega} &\subset\mathbb{S}^d\setminus B_{\upsilon\tau} =\inte B_{\tau\upsilon},
\end{aligned}
\]
and for any path $\tau_0,...,\tau_k$ with no repeated vertices and $k\ge2$, by successive
application we have
\[
E_{\tau_{k-1}\tau_k}\subset\inte E_{\tau_0\tau_1}, \qquad B_{\tau_{k-1}\tau_k}\subset\inte
B_{\tau_0\tau_1}.
\]

Let $\mathcal{P}\subset\mathcal{T}'$ be a connected subtree with $\#V(\mathcal{P})<\infty$, and
define
\[
\partial^+\mathcal{P} =\{(\tau,\upsilon):\tau\in V(\mathcal{P}),\ \upsilon\notin V(\mathcal{P}),\
\tau\sim\upsilon\}.
\]
Its outgoing balls are pairwise disjoint: for $(\tau,\upsilon),(\sigma,\omega)\in\partial^+\mathcal
P$ with $\tau\ne\sigma$, if $\tau=\tau_0,...,\tau_j=\sigma$ is their path in $\mathcal{P}$, by
nesting we have
\[
B_{\sigma\omega}\subset\inte B_{\tau_0\tau_1},\qquad B_{\tau\upsilon}\cap
B_{\tau_0\tau_1}=\emptyset, \qquad B_{\tau\upsilon}\cap B_{\sigma\omega}=\emptyset,
\]
and if $(\tau,\upsilon),(\tau,\omega)\in\partial^+\mathcal{P}$ are distinct, by the disjointness of
the ideal balls associated with $C_\tau$ we have
\[
B_{\tau\upsilon}\cap B_{\tau\omega}=\emptyset.
\]

The identity
\[
D_{\mathcal{P}}:=\bigcup_{\tau\in V(\mathcal{P})}C_\tau =\mathbb{H}^N\setminus
\bigcup_{(\tau,\upsilon)\in\partial^+\mathcal{P}} \inte E_{\tau\upsilon},
\]
is proved by induction on $\#V(\mathcal{P})$. If $V(\mathcal{P})=\{\tau\}$, it is the definition of
$C_\tau$. Let $\mathcal{P}_0$ be the current subtree, let $\upsilon\notin V(\mathcal{P}_0)$ be a
vertex joined to $\tau\in V(\mathcal{P}_0)$, let $\mathcal{P}$ be obtained from $\mathcal{P}_0$ by
adding $\upsilon$ and the edge $(\tau,\upsilon)$, and set $H=E_{\tau\upsilon}$. Let $P_i$ be the
halfspaces $E_{\alpha\beta}$ indexed by
\[
(\alpha,\beta)\in \partial^+\mathcal{P}_0\setminus\{(\tau,\upsilon)\},
\]
and let $Q_j$ be the halfspaces $E_{\upsilon\omega}$ with $\omega\sim\upsilon$ and $\omega\ne\tau$,
so that
\[
H\cap P_i=\emptyset,\qquad Q_j\subset\inte H.
\]
By the induction hypothesis we have
\[
D_{\mathcal{P}_0} =\mathbb{H}^N\setminus \left(\inte H\cup\bigcup_i\inte P_i\right),
\]
while
\[
C_\upsilon=H\setminus\bigcup_j\inte Q_j.
\]
Hence,
\[
\begin{aligned}
D_{\mathcal{P}_0}\cap C_\upsilon&=\partial H,\\
D_{\mathcal{P}} =D_{\mathcal{P}_0}\cup C_\upsilon &=\mathbb{H}^N\setminus \left(\bigcup_i\inte
P_i\cup\bigcup_j\inte Q_j\right),
\end{aligned}
\]
which is the induction step.

We choose a spanning tree $\mathcal{S}\subset X$, a root $v_0$, and a lift
$\mathcal{P}\subset\mathcal T'$ with root $\widetilde v_0$, containing one lift of each vertex of
$X$. The boundary-face count is
\[
\#\partial^+\mathcal{P} =2\#E(X)-2\#E(\mathcal{S}) =2(c+s)-2(n-1)=2r.
\]
We orient the $r$ edges of $X\setminus\mathcal{S}$ as $e_i:v_i\to w_i$. Their chosen lifts determine
unique $g_i\in\G$ by
\[
\widetilde e_i:\widetilde v_i\longrightarrow g_i\widetilde w_i, \qquad \widetilde v_i,\widetilde
w_i\in\mathcal{P},
\]
and set
\[
B_i^+=B_{\widetilde v_i,\,g_i\widetilde w_i}, \qquad B_i^-=B_{\widetilde w_i,\,g_i^{-1}\widetilde
v_i}.
\]
Then we have
\[
\{B_i^+,B_i^-:1\le i\le r\} =\{B_{\tau\upsilon}: (\tau,\upsilon)\in\partial^+\mathcal{P}\}.
\]
For distinct $(\tau,\upsilon),(\sigma,\omega)\in\partial^+\mathcal{P}$ with $\tau=\sigma$, by the
disjointness proved for $C_\tau$ we have
\[
B_{\tau\upsilon}\cap B_{\tau\omega}=\emptyset,
\]
and if $\tau\ne\sigma$ and $\tau=\tau_0,...,\tau_k=\sigma$ is their path in $\mathcal{P}$, then
$B_{\sigma\omega}\subset\inte B_{\tau_0\tau_1}$ and $B_{\tau\upsilon}\cap
B_{\tau_0\tau_1}=\emptyset$, so
\[
B_{\tau\upsilon}\cap B_{\sigma\omega} \subset B_{\tau\upsilon}\cap B_{\tau_0\tau_1} =\emptyset.
\]
Hence all $2r$ closed balls are pairwise disjoint, and
\[
\begin{aligned}
g_iB_i^-&=\mathbb{S}^d\setminus\inte B_i^+,\\
g_i(\mathbb{S}^d\setminus\inte B_i^-) &=\mathbb{S}^d\setminus\inte(\mathbb{S}^d\setminus\inte B_i^+)
=\overline{\inte B_i^+}=B_i^+.
\end{aligned}
\]
For each $v$, let $\sigma_v$ be the path in $\mathcal{S}$ from $v_0$ to $v$. The loops
\[
\beta_i=\sigma_{v_i}e_i\sigma_{w_i}^{-1}, \qquad 1\le i\le r,
\]
form a free basis of $\pi_1(X,v_0)$, and their lifts beginning at $\widetilde v_0$ satisfy
\[
\widetilde\beta_i(0)=\widetilde v_0,\qquad \widetilde\beta_i(1)=g_i\widetilde v_0.
\]
Since $X/\mathcal{S}$ is a bouquet of $r$ circles, by the covering identification we have
\[
\pi_1(X,v_0)=<\beta_1,...,\beta_r\mid\ >, \qquad \pi_1(X,v_0)\xrightarrow{\ \cong\ }\G, \qquad
\beta_i\longmapsto g_i.
\]
Hence $g_1,...,g_r$ freely generate $\G$. They are nontrivial and therefore loxodromic by
torsion-freeness and convex cocompactness. The pairing identities and the pairwise disjoint closed
round balls $B_i^\pm$, $1\le i\le r$, show that $\G$ is classical Schottky.
\end{proof}

Theorem \ref{main} follows in both formulations.
\begin{cor}\label{hausdorff}
Let $\G<\Isom(\mathbb{H}^N)$ be discrete, torsion-free, and non-elementary, and suppose it is
generated by $r\ge2$ elements. If $\delta_\G\le\lambda_r$, then $\G$ is classical, and the
conclusion also holds if $\mathfrak{D}_\G\le\lambda_r$.
\end{cor}
\begin{proof}
Since $\lambda_r<1/2$, Theorem \ref{structural} gives freeness and convex cocompactness with
$\delta_\G=\mathfrak{D}_\G$. By non-elementarity the free rank satisfies $k\ge2$, and since the
abelianization $\mathbb{Z}^k$ is generated by the images of the given $r$ generators, $k\le r$. The
quantities $m_k$, $B_k$, $C_{0,k}$ and $A_k$ increase with $k$, which implies
$\delta_\G\le\lambda_r\le\lambda_k$ and Theorem \ref{rank-r} applies at rank $k$. The
Hausdorff-dimension statement follows since $\delta_\G\le \mathfrak{D}_\G$.
\end{proof}
\begin{thm}\label{n=3}
Let $\G<\Isom(\mathbb{H}^3)$ be discrete, torsion-free, and non-elementary. Assume that either $\G$
is generated by $r\ge2$ elements and $\mathfrak{D}_\G\le\lambda_r$, or $\G$ is generated by two
elements and $\mathfrak{D}_\G\le1/3$. Then $\G$ is classical Schottky. The same conclusion holds
when the corresponding Hausdorff-dimension bound is replaced by the critical-exponent bound. The
elements of $\G$ may preserve or reverse orientation.
\end{thm}
\begin{proof}
We identify $\mathbb{H}^3$ with a totally geodesic subspace $Y\subset\mathbb{H}^4$ and, in Lorentz
coordinates, extend each isometry by
\[
A\longmapsto
\begin{pmatrix}
A&0\\
0&1 \end{pmatrix},
\]
which identifies $\G$ with a discrete subgroup of $\Isom(\mathbb{H}^4)$ acting on $Y$ as before. For
$o\in Y$,
\[
\dis_{\mathbb{H}^4}(o,\g o) =\dis_Y(o,\g o) \qquad(\g\in\G),
\]
so the critical exponent is unchanged, and the limit set lies in
\[
\partial_\infty Y\cong\mathbb{S}^2\subset\mathbb{S}^3
\]
and has the same Hausdorff dimension. The already proved case $N=4$ of the corresponding theorem
gives a free basis $g_1,...,g_k$ and pairwise disjoint closed round balls
$B_i^\pm\subset\mathbb{S}^3$ such that
\[
g_i\bigl(\mathbb{S}^3\setminus\inte B_i^-\bigr)=B_i^+,
\]
and by disjointness and the pairings we have
\[
g_iB_i^+\subset\inte B_i^+,\qquad g_i^{-1}B_i^-\subset\inte B_i^-.
\]
Each $g_i$ is loxodromic on $Y$, since
\[
\inf_{p\in Y}\dis(p,g_ip) \ge\inf_{p\in\mathbb{H}^4}\dis(p,g_ip)>0,
\]
and its attracting and repelling fixed points $\xi_i^\pm$ lie in $\partial_\infty Y$ as limits of
the forward and backward orbits of any point of $Y$. For $z\in B_i^+\setminus\{\xi_i^-\}$, the
iterates $g_i^nz$ lie in $B_i^+$ and converge to $\xi_i^+$, which implies $\xi_i^+\in B_i^+$ since
$B_i^+$ is closed, and then
\[
\xi_i^+=g_i\xi_i^+\in g_iB_i^+\subset\inte B_i^+,
\]
while the same argument applied to $g_i^{-1}$ implies $\xi_i^-\in\inte B_i^-$.

In these coordinates,
\[
\partial_\infty Y=\{(\xi_1,\xi_2,\xi_3,0): \xi_1^2+\xi_2^2+\xi_3^2=1\}\subset\mathbb{S}^3.
\]
For either ball, write
\[
B=\{\xi\in\mathbb{S}^3:u\cdot\xi\ge c\}, \qquad u=(u_1,u_2,u_3,u_4)\in\mathbb{S}^3,\quad -1<c<1,
\]
and set $u_{\partial_\infty Y}=(u_1,u_2,u_3,0)$, so that
\[
u\cdot\xi=u_{\partial_\infty Y}\cdot\xi \qquad(\xi\in \partial_\infty Y).
\]
One fixed point $\xi\in \partial_\infty Y$ lies in $\inte B$ and the other fixed point $\eta\in
\partial_\infty Y$ lies outside $B$, hence
\[
-|u_{\partial_\infty Y}|\le u_{\partial_\infty Y}\cdot\eta<c <u_{\partial_\infty
Y}\cdot\xi\le|u_{\partial_\infty Y}|,
\]
so that $u_{\partial_\infty Y}\ne0$ and $-1<c/|u_{\partial_\infty Y}|<1$, and
\[
B\cap \partial_\infty Y =\{\xi\in \partial_\infty Y: (u_{\partial_\infty Y}/|u_{\partial_\infty
Y}|)\cdot\xi\ge c/|u_{\partial_\infty Y}|\}
\]
is a proper closed round disk with nonempty relative interior. It follows that the pairwise disjoint
disks $D_i^\pm=B_i^\pm\cap \partial_\infty Y$ satisfy $\operatorname{int}_{\partial_\infty Y}
D_i^\pm=\partial_\infty Y\cap\inte B_i^\pm$, and since each $g_i$ preserves $\partial_\infty Y$, the
pairings restrict to
\[
g_i\bigl(\partial_\infty Y \setminus\operatorname{int}_{\partial_\infty Y}D_i^-\bigr) =D_i^+.
\]
Hence the original subgroup of $\Isom(\mathbb{H}^3)$ is classical, including when some generators
reverse orientation.
\end{proof}

\begin{thm}\label{n=2}
Let $\G<\Isom(\mathbb{H}^2)$ be free and convex cocompact of rank $r\ge1$. Then $\G$ is a classical
Schottky group.
\end{thm}
\begin{proof}
The case $r=1$ is trivial, so let $r\ge2$ and set
\[
\G^+=\G\cap\Isom^+(\mathbb{H}^2).
\]
Note that this subgroup has index at most two, is free and convex cocompact, and has the same limit
set as $\G$, and this limit set is a proper subset of $\mathbb{S}^1$, since otherwise convex
cocompactness would make $\mathbb{H}^2/\G^+$ a closed orientable hyperbolic surface, whose second
integral homology is $\mathbb{Z}$, and a free group has zero second integral homology. By
Patterson's theorem \cite{Patterson1976} (Theorem 4.1 and corollary of Theorem 7.1) we have
\[
\dim_H\Lambda_\G =\dim_H\Lambda_{\G^+} =\delta_{\G^+}<1.
\]
In the hyperboloid model, identify $\Isom(\mathbb{H}^2)$ with the Lorentz matrices preserving the
future hyperboloid and define
\[
\iota(A)=\operatorname{diag}(A,\det A),
\]
a faithful closed embedding into $\Isom^+(\mathbb{H}^3)$, since
\[
\det\iota(A)=(\det A)^2=1,
\]
which preserves a totally geodesic plane $Y\cong\mathbb{H}^2$ and restricts there to the original
action, hence $\iota(\G)$ is discrete and convex cocompact, with
\[
\dim_H\Lambda_{\iota(\G)} =\dim_H\Lambda_\G<1,
\]
and Hou's theorem \cite{Hou} gives a free basis $\iota(g_1),...,\iota(g_r)$ and paired disjoint
closed round disks $D_i^\pm\subset\partial_\infty\mathbb{H}^3$. Set $B_i^\pm=D_i^\pm\cap
\partial_\infty Y$. The disk pairings imply
\[
\iota(g_i)D_i^+\subset\inte D_i^+,\qquad \iota(g_i)^{-1}D_i^-\subset\inte D_i^-,
\]
so iteration places the attracting fixed point in $\inte D_i^+$ and the repelling fixed point in
$\inte D_i^-$, and as both fixed points lie on $\partial_\infty Y$, the sets $B_i^\pm$ are proper
closed round intervals with nonempty relative interiors. They are pairwise disjoint, and by
restriction of the disk pairings we have
\[
g_i\bigl(\partial_\infty Y\setminus\inte_{\partial_\infty Y} B_i^-\bigr)=B_i^+ \qquad(1\le i\le r).
\]
Hence $\G$ is classical, without an orientation assumption on its elements.
\end{proof}


\section{Lorentz Normals and Pole Radius}
In this section we give a necessary and sufficient condition for a round domain in terms of the
outward Lorentz normals. The criteria so far use ratio balls. The following condition on the outward
Lorentz normals is necessary and sufficient for a chosen tuple to admit a domain bounded by
arbitrary round spheres, and it also implies that the tuple freely generates a discrete group.

\begin{thm}\label{lorentz}
Let $r\ge1$ be an integer, let $G_1,...,G_r$ be Lorentz matrices on $\mathcal{V}$ preserving the
future hyperboloid, let $g_i$ be the represented isometries, and set $\G=\langle
g_1,...,g_r\rangle$. For $n_i\in\mathcal{V}$, $1\le i\le r$, let
\[
n_{i,-}=n_i,\qquad n_{i,+}=-G_i n_i,\qquad S=\sum_i(n_{i,-}+n_{i,+}),
\]
where $s,t$ range over $\mathcal{I}_r=\{(i,-),(i,+):1\le i\le r\}$, ordered by
$(1,-)<(1,+)<\cdots<(r,-)<(r,+)$, $n_s=n_{i,\sigma}$ for $s=(i,\sigma)$, and $S_0$ is the zeroth
coordinate of $S$. Then the tuple $g_1,...,g_r$ admits a classical round-ball domain if and only if
the system
\begin{equation}\label{normals}
\begin{aligned}
\langle n_i,n_i\rangle&=1 &&(1\le i\le r),\\
\langle n_s,n_t\rangle&<-1 &&(s\ne t),\\
S_0&>0&
\end{aligned}
\end{equation}
has a solution.
\end{thm}
\begin{proof}
For a unit spacelike normal $n$, let $E_n=\{X\in\mathbb{H}^N:\langle X,n\rangle\ge0\}$ and
$B_n=\partial_\infty E_n$. Assume (\ref{normals}) and set $m=2r$. Summing the pair inequalities, we
get
\[
\langle S,S\rangle =m+2\sum_{s<t}\langle n_s,n_t\rangle <m(2-m)\le0, \qquad \langle
S,n_s\rangle<2-m\le0.
\]
Since $S_0>0$, the point $o=S/\sqrt{-\langle S,S\rangle}$ is future timelike and lies strictly
outside every halfspace. We conjugate $o$ to $(1,0,...,0)$, so that
\[
n_s=(\sinh t_s,\cosh t_s\,u_s),\qquad t_s>0,\quad u_s\in\mathbb{S}^{N-1}.
\]
The ideal cap has center $u_s$ and angular radius $\theta_s=\arccos(\tanh t_s)<\pi/2$. The pair
inequality is therefore equivalent to
\[
u_s\cdot u_t< \tanh t_s\tanh t_t-\operatorname{sech}t_s\operatorname{sech}t_t
=\cos(\theta_s+\theta_t),
\]
so $\angle(u_s,u_t)>\theta_s+\theta_t$ and the closed balls are disjoint, with the pairing
\[
G_i B_{-n_{i,-}}=B_{-G_i n_{i,-}}=B_{n_{i,+}}.
\]

Each generator is loxodromic: with $E_{i,\pm}=E_{n_{i,\pm}}$, set
\[
\begin{aligned}
\alpha_i&=\dis(\partial E_{i,-},\partial E_{i,+})>0,\\
\beta_i&=\min\{\dis(o,\partial E_{i,-}),\dis(o,\partial E_{i,+})\}>0.
\end{aligned}
\]
Since $G_iE_{i,-}=\mathbb{H}^N\setminus\inte E_{i,+}$, the halfspaces
$E_{i,+},G_iE_{i,+},...,G_i^{k-1}E_{i,+}$ are strictly nested, and inverse powers of $G_i$ send
consecutive boundary planes to $\partial E_{i,-},\partial E_{i,+}$, so their gap is $\alpha_i$, and
the initial and terminal crossing distances are at least $\beta_i$, so
\[
\dis(o,G_i^ko)\ge2\beta_i+(k-1)\alpha_i\qquad(k\ge1).
\]
Hence we have
\[
\lim_{k\to\infty}\frac{\dis(o,G_i^ko)}{k}\ge\alpha_i>0,
\]
and by the classification of real hyperbolic isometries $G_i$ is loxodromic. We choose $z$ outside
all the ideal balls and write $B_{g_i}=B_{n_{i,+}}$ and $B_{g_i^{-1}}=B_{n_{i,-}}$. For every
nonempty reduced word $w=s_1\cdots s_k$, the pairing implies $wz\in\inte B_{s_1}$, which implies
$wz\ne z$, and the tuple is a free basis. Proposition \ref{gap} gives discreteness and convex
cocompactness, and the result follows.

For the converse, take the unit normals of a paired domain, hence $n_{i,+}=-G_i n_{i,-}$, and a
future point $x\in\mathbb{H}^N$ outside all halfspaces, obtained by taking it near an ideal point
outside all balls. After conjugating $x$ to the origin, the same cap calculation gives the strict
pair inequalities, and since
\[
\langle S,S\rangle<m(2-m)\le0, \qquad \langle S,x\rangle =\sum_s\langle n_s,x\rangle<0,
\]
the vector $S$ is future timelike, which implies $S_0>0$.
\end{proof}

\begin{cor}\label{limit}
The generating tuples admitting a round domain form an open set among Lorentz tuples. Note that
suppose that $G_i^{(k)}\to G_i$ as $k\to\infty$, for $1\le i\le r$, that for each $k$ the vectors
$n_{i,-}^{(k)}$ satisfy (\ref{normals}) for this tuple with bounded Euclidean norms, and set
$n_{i,+}^{(k)}=-G_i^{(k)}n_{i,-}^{(k)}$. If there is a common $\eta>0$ such that
\[
\langle n_s^{(k)},n_t^{(k)}\rangle\le-1-\eta \quad(s\ne t),
\]
then the limit tuple admits a round domain as well.
\end{cor}
\begin{proof}
We fix a matrix norm $\|\cdot\|$. For Lorentz matrices $G_1',...,G_r'$ satisfying
\[
\max_{1\le i\le r}\|G_i'-G_i\|<\varepsilon.
\]
Here $\varepsilon>0$ is sufficiently small, keep $n_{i,-}=n_i$ and set
\[
n_{i,+}'=-G_i'n_i,\qquad S'=\sum_{i=1}^r(n_i+n_{i,+}').
\]
The unit norms are preserved, and continuity gives the strict pair inequalities and $S_0'>0$ for
nearby tuples, which is openness.

For the limit statement, choose a convergent subsequence of the bounded normals $n_{i,-}^{(k)}$.
Since $G_i^{(k)}\to G_i$, the paired normals converge as well. The unit normalizations and uniform
pair inequalities pass to the limit, and the sum is timelike with $S_0\ge0$, which implies $S_0>0$,
and Theorem \ref{lorentz} applies.
\end{proof}

The centers are restricted next to the chosen poles, with unequal radii allowed. Let $a,b_1,...,b_q$
be a free basis of a convex cocompact Kleinian group $\G$, where $q\ge1$ and $r=q+1$, set
$L=\ell(a)>0$, and normalize
\[
a(x,t)=(e^LAx,e^Lt),\qquad A\in O(N-1).
\]
We fix $r_0>0$ and set
\[
\begin{aligned}
B_a^-&=\{x\in\mathbb{R}^{N-1}:|x|\le r_0\},\\
B_a^+&=\{x\in\mathbb{R}^{N-1}:|x|\ge r_0e^L\}\cup\{\infty\},\\
a\bigl(\mathbb{S}^d\setminus\inte B_a^-\bigr)&=B_a^+.
\end{aligned}
\]
For $1\le i\le q$, let $u_i=b_i(\infty)$, $v_i=b_i^{-1}(\infty)$ and $K_i>0$ be the pole parameters
in (\ref{pole-form}), and assume that all pole centers lie in $r_0<|x|<r_0e^L$. For the signed
labels $s=(i,-)$ and $s=(i,+)$, set
\[
z_{i,-}=v_i,\qquad z_{i,+}=u_i,\qquad c_s=\min\{|z_s|-r_0,r_0e^L-|z_s|\},\qquad d_{st}=|z_s-z_t|.
\]
We order the signed labels by $(1,-)<(1,+)<\cdots<(q,-)<(q,+)$, denote by $s<t$ unordered distinct
pairs, and for $s=(i,\sigma)$ write $\bar s=(i,-\sigma)$, where $-(+)=-$ and $-(-)=+$. Then $c_s>0$,
and $d_{st}>0$ for $s<t$: since $a$ belongs to this free basis, its endpoint stabilizer is $\langle
a\rangle$, and equal signed poles would give a nonempty reduced word involving the $b_i$ in this
stabilizer, contrary to freeness.

\begin{thm}\label{dual}
There are mutually disjoint paired balls centered at the poles and lying strictly inside the annulus
if and only if there are positive radii such that
\begin{equation}\label{primal}
r_{i,-}r_{i,+}\ge K_i,\qquad r_s<c_s,\qquad r_s+r_t<d_{st}\quad(s\ne t).
\end{equation}
For real numbers $\lambda_s\ge0$, one for each signed label $s$, and $\lambda_{st}\ge0$, one for
each pair $s<t$, set
\[
\lambda_{ts}=\lambda_{st},\qquad \lambda_{ss}=0, \qquad
\lambda=((\lambda_s)_s,(\lambda_{st})_{s<t}),
\]
so that $\lambda\ne0$ means
\[
\sum_s\lambda_s+\sum_{s<t}\lambda_{st}>0,
\]
and set
\[
\begin{split}
w_s&=\lambda_s+\sum_{t\ne s}\lambda_{st},\\
C(\lambda)&=\sum_s\lambda_sc_s+\sum_{s<t}\lambda_{st}d_{st},\\
D(\lambda)&=2\sum_i\sqrt{K_iw_{i,-}w_{i,+}}.
\end{split}
\]
The system (\ref{primal}) has a solution if and only if
\begin{equation}\label{dual-eq}
D(\lambda)<C(\lambda) \quad\text{for every nonzero nonnegative family }\lambda.
\end{equation}
\end{thm}
\begin{proof}
Each closed ball lies inside the annulus if and only if $r_s<c_s$, and the closed balls are pairwise
disjoint if and only if $r_s+r_t<d_{st}$ for all $s<t$. By (\ref{pole-form}),
\[
|b_i(x)-u_i|=\frac{K_i}{|x-v_i|},
\]
so $b_i$ maps the exterior of the ball centered at $v_i$ with radius $r_{i,-}$ onto the ball
centered at $u_i$ with radius $K_i/r_{i,-}$, and the required pairing holds if and only if
\[
r_{i,+}=\frac{K_i}{r_{i,-}}, \qquad\text{equivalently}\qquad r_{i,-}r_{i,+}=K_i.
\]
If this product exceeds $K_i$, multiplying both radii by
\[
\left(\frac{K_i}{r_{i,-}r_{i,+}}\right)^{1/2}\le1
\]
preserves every strict inequality and makes the product an equality. Hence (\ref{primal}) is
equivalent to the required domain.

We write $\mathbf r=(r_s)_s\in\mathbb{R}^{2q}$ and consider the set
\[
\mathcal{C}=\{\mathbf r\in(0,\infty)^{2q}:r_{i,-}r_{i,+}\ge K_i\}.
\]
Note that this set is convex, since each pair satisfies
\[
\left|(2\sqrt{K_i},r_{i,-}-r_{i,+})\right| \le r_{i,-}+r_{i,+},
\]
or equivalently $r_{i,+}\ge K_i/r_{i,-}$. By the arithmetic--geometric mean inequality we have
\[
\inf_{\mathbf r\in\mathcal{C}}\sum_s w_sr_s =2\sum_i\sqrt{K_iw_{i,-}w_{i,+}}=D(\lambda),
\]
where, if one of $w_{i,-},w_{i,+}$ vanishes and the other is positive, the radius with zero weight
is sent to infinity and its partner to zero, and if both vanish, the pair contributes zero.

Set $p=2q+\binom{2q}{2}$ and define the linear mapping $M:\mathbb{R}^{2q}\to\mathbb{R}^p$ and
vectors $b,\lambda\in\mathbb{R}^p$, with the same row order in all three expressions, by
\[
\begin{aligned}
M\mathbf r&=((r_s)_s,(r_s+r_t)_{s<t}),\\
b&=((c_s)_s,(d_{st})_{s<t}),\\
\lambda&=((\lambda_s)_s,(\lambda_{st})_{s<t}).
\end{aligned}
\]
Then $\lambda\cdot M\mathbf r=\sum_s w_sr_s$ and $\lambda\cdot b=C(\lambda)$. If $\mathbf
r\in\mathcal{C}$ satisfies $M\mathbf r<b$, then, for every nonzero $\lambda\ge0$,
\[
D(\lambda)\le\lambda\cdot M\mathbf r <\lambda\cdot b=C(\lambda).
\]
Conversely, suppose the convex set $M\mathcal{C}$ is disjoint from the open convex set
$\{z\in\mathbb{R}^p:z<b\}$. The separating hyperplane theorem implies
$\lambda\in\mathbb{R}^p\setminus\{0\}$ and $\alpha\in\mathbb{R}$ such that
\[
\lambda\cdot M\mathbf r\ge\alpha\ge\lambda\cdot z \qquad(\mathbf r\in\mathcal{C},\ z<b).
\]
Taking the infimum on the left and the supremum on the right, we get $\inf_{\mathbf r\in\mathcal
C}\lambda\cdot M\mathbf r\ge\sup_{z<b}\lambda\cdot z$. If $\lambda_j<0$, letting $z_j\to-\infty$
while keeping the other coordinates fixed below their corresponding $b_k$ implies $\lambda\cdot
z\to+\infty$, so a bounded supremum requires $\lambda\ge0$. For $\lambda\ge0$, $\lambda\cdot
z\le\lambda\cdot b$ whenever $z<b$, and by taking $z_j=b_j-\varepsilon$ for all $j$ we get
\[
\sup_{z<b}\lambda\cdot z =\lim_{\varepsilon\downarrow0} \sum_{j=1}^p\lambda_j(b_j-\varepsilon)
=\lambda\cdot b=C(\lambda).
\]
The infimum on the left is $D(\lambda)$, contrary to (\ref{dual-eq}).
\end{proof}

\begin{cor}\label{optimization}
Set
\[
\tau_* =\inf\left\{\tau\in\mathbb{R}:
\begin{array}{l}
\text{there exist }r_s>0\text{ such that}\\
r_{i,-}r_{i,+}\ge K_i\quad(1\le i\le q),\\
r_s\le\tau c_s\quad\text{for every }s,\\
r_s+r_t\le\tau d_{st}\quad(s<t) \end{array}\right\}.
\]
Then we have
\[
\tau_*=\max_{\lambda\ge0,\ C(\lambda)=1}D(\lambda).
\]
Both extrema are attained. The paired Euclidean balls can be chosen inside the fixed annulus to
form, together with the fixed balls $B_a^-,B_a^+$ paired by $a$, a classical domain for the whole
group if and only if $\tau_*<1$.
\end{cor}
\begin{proof}
We keep $M,b,\mathbf r,\mathcal{C}$ from the proof of Theorem \ref{dual}. Set
$r_{i,-}=r_{i,+}=\sqrt{K_i}$ and choose
\[
\tau_0>\max\left\{\max_s\frac{r_s}{c_s},\ \max_{s<t}\frac{r_s+r_t}{d_{st}}\right\}.
\]
These choices satisfy all the required inequalities. For $\tau\le\tau_0$, each radius lies between
$K_i/(\tau_0c_{\bar s})$ and $\tau_0c_s$, and $\tau\ge\sqrt{K_i/(c_{i,-}c_{i,+})}>0$. A minimizing
sequence within these bounds has a convergent subsequence by compactness, and since the lower radius
bounds are positive and the inequalities are closed, the limit satisfies all the required
inequalities.

For $\lambda\ge0$ with $C(\lambda)=1$, for every $\mathbf r\in\mathcal{C}$ satisfying $M\mathbf
r\le\tau b$ we have
\[
D(\lambda)\le\lambda\cdot M\mathbf r \le\tau\lambda\cdot b=\tau,
\]
and taking the minimum over such $(\mathbf r,\tau)$ implies $\tau_*\ge D(\lambda)$. For
$\tau<\tau_*$, no $\mathbf r\in\mathcal{C}$ satisfies $M\mathbf r<\tau b$, so the separating
hyperplane theorem provides a nonzero $\lambda\ge0$ such that
\[
D(\lambda)\ge\tau C(\lambda).
\]
Since $C(\lambda)>0$, replacing $\lambda$ by $\lambda/C(\lambda)$ implies $C(\lambda)=1$ and
$D(\lambda)\ge\tau$, and by letting $\tau\uparrow\tau_*$ we get
\[
\sup_{\lambda\ge0,\ C(\lambda)=1}D(\lambda)\ge\tau_*,
\]
which together with the opposite inequality gives equality. Since $\{\lambda\ge0:C(\lambda)=1\}$ is
compact and $D$ is continuous, the supremum is attained. The system (\ref{primal}) has a solution if
and only if $\tau_*<1$, and reducing the radii so that $r_{i,-}r_{i,+}=K_i$ and adding $B_a^-,B_a^+$
gives a classical domain for $\G$.
\end{proof}


\section{Invariant Metrics and Rectifiable Trees}
In this section, following Bowen's rectifiability of quasicircles \cite{Bowen} and its use in Hou's
classification \cite{Hou}, we study invariant metrics and rectifiable trees in higher dimensions.
Denote $\Omega_\G=\mathbb{S}^{N-1}\setminus\Lambda_\G$.
\begin{prop}\label{derivative}
Let $\G<\Isom(\mathbb{H}^N)$ be discrete with $P_\G(1,o)<\infty$, and in the ball model centered at
$o=0$ set $x_\g=\g^{-1}o$. Then the series
\[
\rho(\xi)=\sum_{\g\in\G} \frac{1-|x_\g|^2}{|\xi-x_\g|^2}
\]
is positive, converges with all derivatives locally on $\Omega_\G$, and satisfies
\[
\rho(g\xi)|Dg(\xi)|=\rho(\xi)\qquad(g\in\G,\ \xi\in\Omega_\G).
\]
Here $|Dg(\xi)|$ is the conformal dilation of the boundary mapping in $g_{\rm round}$, hence
$\rho^2g_{\rm round}$ is a smooth metric on $\Omega_\G$ invariant under the group. For every
piecewise smooth compact arc $c\subset\Omega_\G$,
\begin{equation}\label{orbit-length}
\sum_{\g\in\G}\leg_{\rm round}(\g c) =\int_c\rho\,ds<\infty.
\end{equation}
\end{prop}
\begin{proof}
The conformal derivative formula and convergence on $\Omega_\G$ follow from \cite{FMS2010} (Lemma
4.1). On each compact subset of $\Omega_\G$, every derivative of each summand is bounded by a
constant times $e^{-\dis(o,\g o)}$, so the series converges locally with all derivatives, and the
identity term implies $\rho\ge1$. By the chain rule and reindexing we have
\[
\rho(g\xi)|Dg(\xi)| =\sum_{\g\in\G}|D(\g g)(\xi)| =\rho(\xi),
\]
and Tonelli's theorem implies (\ref{orbit-length}), the integral being finite since $\rho$ is
continuous on the compact arc $c$.
\end{proof}
The length estimate is applied next to an embedded lift of a free basis, with the following tree
conventions. A geodesic metric space is a real tree if its geodesic segments are unique and every
geodesic triangle is a tripod. The valence of a point is the number of components of its complement,
and an endpoint has valence one. For a tree whose vertices have valence $<\infty$, an end is an
equivalence class of rays, where two rays are equivalent when they have a common tail.
\begin{thm}\label{tree}
Let $N\ge4$, and let $\G<\Isom(\mathbb{H}^N)$ be a discrete non-elementary torsion-free convex
cocompact group such that $\dim_H\Lambda_\G<1$. Then there is an embedded tree $T\subset\Omega_\G$,
invariant under the group, such that
\[
\overline T=T\cup\Lambda_\G,\qquad \mathcal{H}^1(T\cup\Lambda_\G)<\infty.
\]
The closure is a compact dendrite, that is, a compact connected locally connected space containing
no simple closed curve. Its endpoint set is $\Lambda_\G$.
\end{thm}
\begin{proof}
Set $C=\Hu(\Lambda_\G)$ and use the ball model centered at $o=0$. The action on the proper geodesic
space $C$ is proper and cocompact, which implies that $\G$ is finitely generated by the
\v{S}varc--Milnor lemma \cite{BH} (Proposition I.8.19), and free by Theorem \ref{structural}. We
choose a free basis $g_1,...,g_r$ with $r\ge2$. By Corollary \ref{simply-connected}, $\Omega_\G$ is
path connected and simply connected, and the action on $\Omega_\G$ is properly discontinuous, as
follows. Let $\mathcal{K}_o$ be the compact group fixing $o$, and in boundary coordinates with
endpoints $0,\infty$ of a fixed axis through $o$ set $a_D(z)=e^Dz$. Given $\g$, set $D=\dis(o,\g
o)$, choose $k\in\mathcal{K}_o$ with $ka_Do=\g o$, and set $l=a_{-D}k^{-1}\g\in\mathcal{K}_o$, so
that
\[
\g=k a_D l,\qquad k,l\in\mathcal{K}_o,\qquad D=\dis(o,\g o),\qquad a_D(z)=e^Dz.
\]
The set $\mathcal{K}_o\{a_D:0\le D\le R\}\mathcal{K}_o$ is compact for all $R>0$, so for distinct
$\g_j\in\G$ discreteness forces $D_j\to\infty$, and after passage to a subsequence $k_j\to k$ and
$l_j\to l$. Since $e^{D_j}z\to\infty$ uniformly on compact sets avoiding zero,
\[
\begin{gathered}
\g_jo\longrightarrow\xi=k(\infty)\in\Lambda_\G, \qquad
\g_j^{-1}o\longrightarrow\zeta=l^{-1}(0)\in\Lambda_\G,\\ \g_j\longrightarrow\xi \quad\hbox{locally
uniformly on }\mathbb{S}^{N-1}\setminus\{\zeta\}.
\end{gathered}
\]
For compact $K\subset\Omega_\G$, this convergence implies $\g_jK\cap K=\emptyset$ eventually, so
\[
\#\{\g\in\G:\g K\cap K\ne\emptyset\}<\infty.
\]
Note that every nonidentity element is loxodromic and has its fixed points in $\Lambda_\G$, so the
action on $\Omega_\G$ is also free, and $\Omega_\G\to M=\Omega_\G/\G$ is the universal covering. We
choose $m_0\in M$ and a lift $p\in\Omega_\G$, and identify $\pi_1(M,m_0)$ with $\G$ through this
lift.

We realize the basis by a bouquet of piecewise smooth loops $f_i:I=[0,1]\to M$ based at $m_0$, with
smooth open edges. After a homotopy fixing the vertex, the $2r$ edge segments meeting there are
mapped, in local coordinates centered at its image, to straight segments in pairwise distinct
directions. These edge segments are kept fixed, and $d=N-1\ge3$. Transversality for families of
mappings and for their values and first derivatives \cite{Hirsch} (Chapter 3, Theorems 2.7--2.8) is
used in the following relative form. In a coordinate chart, let $f:I\to\mathbb{R}^d$ be the edge
mapping. For distinct parameter values $s_1,...,s_q$ outside the fixed intervals, choose smooth bump
functions $\chi_\alpha$ supported away from those intervals with
$\chi_\alpha(s_\beta)=\delta_{\alpha\beta}$, and for small $\mathbf
v=(v_1,...,v_q)\in(\mathbb{R}^d)^q$ set
\[
f_{\mathbf v}(s) =f(s)+\sum_{\alpha=1}^q\chi_\alpha(s)v_\alpha.
\]
Then $f_{\mathbf v}(s_\alpha)=f(s_\alpha)+v_\alpha$, so these values vary independently, and at
$s=s_0$, replacing a bump $\chi_\alpha(s)$ with $(s-s_0)\chi_\alpha(s)$ also allows the first
derivative to vary while keeping the value fixed, hence the zero-derivative condition has
codimension $d>1$, so every open edge can be made an immersion. For parameters $s,t\in\inte I$
representing distinct points of the bouquet, arrange that
\[
(s,t)\longmapsto(f_i(s),f_j(t))
\]
is transverse to the diagonal $\Delta_M=\{(x,x):x\in M\}$ in $M\times M$, also when one parameter
lies in a fixed edge segment. Exhausting these parameter pairs by compact subsets and applying
Baire's theorem gives all these conditions simultaneously, and the fixed edge segments already have
disjoint interiors. An intersection $f_i(s)=f_j(t)=x$ would require
\[
df_i(T_sI)+df_j(T_tI)=T_xM,
\]
which is impossible since the left side has dimension at most two and $\dim M=d>2$. Immersion gives
injectivity on sufficiently short parameter intervals, and transversality of each $f_i$ to $\{m_0\}$
excludes $f_i(s)=m_0$ for $s\in\inte I$, since $1<d$, hence the bouquet is embedded and still
represents the chosen basis.

The induced mapping on fundamental groups is surjective, so the full lift is connected, and
injective, so the lift has no loop, hence the lift is an embedded Cayley tree $T$. Let $c_i$ be the
piecewise smooth lift of $f_i$ from $p$ to $g_ip$. The conical-limit-set equality implies
$\mathfrak{D}_\G=\delta_\G<1$, hence $P_\G(1,o)<\infty$, so Proposition \ref{derivative} applies and
by (\ref{orbit-length}) we have
\[
\mathcal{H}^1(T)\le \sum_{i=1}^r\sum_{\g\in\G} \leg_{\rm round}(\g c_i)<\infty.
\]
We give $T$ the arclength path metric $d_T$ and denote by $[x,y]$ the unique joining arc in $T$.
Since the total edge length is bounded, for $\eta>0$ there are finitely many edges whose complement
has length sum below $\eta$. Joining them by their tree paths produces a compact connected subtree
$F$, and each component of $T\setminus F$ attaches at one point of $F$ and has diameter and distance
from $F$ at most $\eta$. It follows that $T$ is totally bounded, and its completion $\widehat T$ is
compact.

The completion is a real tree. For $z,u,v\in T$, set
\[
a=d_T(z,u),\qquad b=d_T(z,v),\qquad c=\leg([z,u]\cap[z,v]).
\]
For the points $u_s\in[z,u]$ and $v_s\in[z,v]$ at the fraction $s\in[0,1]$ from $z$,
\[
\begin{aligned}
d_T(u_s,v_s)&=s(a+b)-2\min\{sa,sb,c\}\\
&\le s(a+b-2c)=s\,d_T(u,v),
\end{aligned}
\]
since $a,b\ge c$ implies $\min\{sa,sb,c\}\ge sc$. For constant-speed segments
$c_{xy},c_{x'y'},c_{x'y}:[0,1]\to T$, by the triangle inequality and this common-endpoint estimate
we have
\[
\begin{aligned}
d_T(c_{xy}(t),c_{x'y'}(t)) &\le d_T(c_{xy}(t),c_{x'y}(t)) +d_T(c_{x'y}(t),c_{x'y'}(t))\\
&\le(1-t)d_T(x,x')+t\,d_T(y,y').
\end{aligned}
\]
We keep the notation $d_T$ for the completed metric. For $x,y\in\widehat T$, choose $x_n,y_n\in T$
with $x_n\to x$, $y_n\to y$ and write $c_n=c_{x_ny_n}$, so that
\[
d_T(c_n(t),c_k(t))\le (1-t)d_T(x_n,x_k)+t d_T(y_n,y_k),\qquad 0\le t\le1,
\]
and these segments have a uniform limit $c$ with $d_T(c(s),c(t))=|s-t|d_T(x,y)$. For $x,y,z,w\in T$,
the four-point inequality
\[
d_T(x,y)+d_T(z,w)\le \max\{d_T(x,z)+d_T(y,w),\ d_T(x,w)+d_T(y,z)\},
\]
passes to the completion. If $u,v$ lie at parameter $t$ on two segments from $x$ to $y$, the
inequality applied to $u,v,x,y$ implies $d_T(u,v)+d_T(x,y)\le d_T(x,y)$, which implies $u=v$. For
the tripod property, take three points $x,y,z\in\widehat T$ and set
\[
L=d_T(x,y),\qquad A=d_T(x,z),\qquad B=d_T(y,z),\qquad a=\frac{L+A-B}{2}.
\]
By the triangle inequalities, $0\le a\le L$. With $[x,y]$ the unique geodesic segment in $\widehat
T$, choose $b\in[x,y]$ with $d_T(x,b)=a$. By the four-point inequality we have
\[
L+d_T(z,b) \le\max\{A+L-a,a+B\}=L+A-a,
\]
and the triangle inequality provides the opposite bound $d_T(z,b)\ge A-a$, so
\[
d_T(z,b)=A-a,\qquad d_T(y,b)+d_T(b,z)=L+A-2a=B.
\]
Hence $b$ lies on both $[x,z]$ and $[y,z]$, and since the segments are unique, their union with
$[x,y]$ is a tripod. It follows that $\widehat T$ is a real tree, its metric balls are convex and
path connected, it contains no simple closed curve, and being compact, it is a dendrite.

Since each vertex has $2r$ edges of positive length, a small metric ball consists of their initial
segments, and each interior edge point has a neighborhood in its edge, which implies that $d_T$
gives the original topology on $T$. We fix a root vertex and let $F_j$ be the subtree spanned by the
vertices at combinatorial distance at most $j$ from the root. These compact connected subtrees
exhaust $T$. For a Cauchy sequence $(z_n)$ having no limit in $T$, there is $\varepsilon_j>0$ such
that eventually $d_T(z_n,F_j)\ge\varepsilon_j$, since otherwise compactness would give a subsequence
converging in $F_j$, and the Cauchy property would force the whole sequence to converge there. If
two sufficiently large index terms of the sequence lie in different complementary components, their
distance is at least $2\varepsilon_j$, so the Cauchy property puts all sufficiently large $z_n$ in
one component of $T\setminus F_j$, and the nested components determine a unique combinatorial ray
from the root. Conversely, every ray is Cauchy, since the sum of its edge lengths converges, and
distinct rays have distinct limits, since their lengths beyond the branching vertex are positive,
hence $\widehat T\setminus T$ consists of the topological ends.

The inclusion into the round sphere is $1$-Lipschitz for $d_T$, since the joining arclength bounds
round distance, and since the sphere is complete, this inclusion extends continuously to $\widehat
T$.

We choose $o_0\in C$, set $\mathcal{S}=\{g_i,g_i^{-1}:1\le i\le r\}$, and denote by $|\g|_{\mathcal
S}$ word length. The orbit mapping from the chosen Cayley tree to $C$, extended geodesically along
edges, is a quasi-isometry by the \v{S}varc--Milnor lemma, so there are constants $\kappa\ge1$,
$\varepsilon_0\ge0$ and $R_C>0$ such that
\[
\kappa^{-1}|\g^{-1}\eta|_{\mathcal{S}}-\varepsilon_0 \le\dis(\g o_0,\eta o_0)
\le\kappa|\g^{-1}\eta|_{\mathcal{S}}+\varepsilon_0, \qquad \sup_{x\in C}\dis(x,\G o_0)\le R_C.
\]
The stability theorem for quasi-geodesics \cite{BH} (Theorem III.H.1.7) provides $R_{\rm M}>0$ such
that every image word segment has Hausdorff distance at most $R_{\rm M}$ from its joining hyperbolic
segment, and taking limits gives the same bound for rays, so each word ray has a unique endpoint in
$\partial_\infty C=\Lambda_\G$. If two image rays have the same endpoint, their joining geodesic
rays from $o_0$ agree, so their image vertex sets have bounded Hausdorff distance, since the image
edge lengths are uniformly bounded. The lower orbit inequality gives the same property for their
word vertex sets, and two rays from one root in a tree have this property only when they agree.

For surjectivity, fix a unit-speed ray $\sigma:[0,\infty)\to C$ from $o_0$ to $\xi\in\Lambda_\G$,
and for integers $n>R_C+1$ choose $\eta_n\in\G$ such that
\[
\dis(\eta_no_0,\sigma(n))\le R_C,\qquad D_n=\dis(o_0,\eta_no_0)\ge n-R_C.
\]
Let $\sigma_n:[0,D_n]\to C$ be the unit-speed segment to $\eta_no_0$, and let $\theta_n$ be its
angle with $\sigma$ at $o_0$. For $0\le t\le D_n$, by the hyperbolic cosine rule we have
\[
0\le1-\cos\theta_n \le\frac{\cosh R_C}{\sinh D_n\sinh n},\qquad \cosh\dis(\sigma_n(t),\sigma(t))
=1+\sinh^2t\,(1-\cos\theta_n),
\]
so $\sigma_n\to\sigma$ uniformly on bounded parameter intervals. The word segments from $1$ to
$\eta_n$ have lengths tending to infinity, and since the Cayley tree has bounded vertex degrees, a
subsequence agrees eventually at every fixed depth with a word ray $(\g_j)_{j\ge0}$. For fixed $j$,
stability gives points $z_n\in\sigma_n$ such that
\[
\dis(\g_jo_0,z_n)\le R_{\rm M},\qquad \dis(o_0,z_n)\le\dis(o_0,\g_jo_0)+R_{\rm M}.
\]
By local convergence of $\sigma_n$ we have $\dis(\g_jo_0,\sigma)\le R_{\rm M}$, and by the lower
orbit bound $\dis(o_0,\g_jo_0)\to\infty$, hence $\g_jo_0\to\xi$, hence the word ends correspond
bijectively to $\Lambda_\G$, and since
\[
\dis(\g o,\g o_0)=\dis(o,o_0)\qquad(\g\in\G),
\]
the two orbit mappings have the same ideal limits.

For the continuous extension of $T\hookrightarrow\mathbb{S}^{N-1}$ on an end, let $\g_jp$ be its
successive ray vertices, hence $\g_jo\to\xi$. The compact-factor argument above shows, along a
further subsequence of any subsequence,
\[
\g_j\longrightarrow\xi \quad\hbox{locally uniformly off some }\zeta\in\Lambda_\G.
\]
Since $p\in\Omega_\G$, $\g_jp\to\xi$ along this further subsequence, so the whole vertex sequence
converges to $\xi$, and the extension on the end has the same value.

Different ends map to different points of $\Lambda_\G$, and these images are disjoint from
$T\subset\Omega_\G$, so the extension is injective, and being continuous on a compact space, it is a
homeomorphism onto $\overline T=T\cup\Lambda_\G$. Since $\mathcal{H}^1(\Lambda_\G)=0$, the
one-dimensional measure of the closure equals that of $T$ and is bounded. The valence is $2r$ at
each original vertex and two in each edge interior. Let $e\in\widehat T\setminus T$. For two points
distinct from $e$, their rooted paths either terminate in $T$ or represent ends distinct from $e$,
so their joining arc omits $e$. If $\widehat T\setminus\{e\}$ had two components, the joining arc of
a point in each component would pass through $e$, which is a contradiction, hence this complement is
connected and $e$ has valence one. No point of $T$ is an endpoint, since its valence is two or
$2r\ge4$. Since $\overline T=T\cup\Lambda_\G$, the endpoint set of $\overline T$ is $\Lambda_\G$.
\end{proof}

\begin{rem}
Note that this invariant construction can be seen as a higher-dimensional generalization of Bowen's
rectifiable curve \cite{Bowen} and is used in \cite{Hou}.
\end{rem}
\centerline{ACKNOWLEDGEMENTS}\par The author would like to thank Peter Sarnak, Benson Farb, Dave
Gabai, and Peter Shalen for their unwavering support and many great mathematical inspirations.  The author owe debt of gratitude to  Ying Zhou for their invaluable support and guidance throughout this work. Claude is also used for grammar editing.\par
\text{}\\ E-mail: houyong@cuhk.edu.cn 

\pdfbookmark[1]{Reference}{Reference}
\bibliographystyle{plain}

\end{document}